\documentclass[a4paper,10pt, twoside,titlepage]{book}
\usepackage[T1]{fontenc} % imposta la codifica dei font
\usepackage{pifont}
\usepackage[utf8,latin1]{inputenc} % lettere accentate da tastiera
\usepackage[italian]{babel} 
\usepackage{layaureo} % imposta i margini di pagina
\usepackage{lipsum} %  genera testo fittizio
\usepackage{url} % per scrivere gli indirizzi Internet
\usepackage{amssymb}
\usepackage{amsmath}
\usepackage{amsfonts}
\usepackage{amsthm}
\usepackage{graphicx}
\usepackage{epstopdf}
\usepackage{enumitem}
\usepackage{tikz}
\usepackage{emptypage}
\usepackage{hyperref}
\usepackage{fancyhdr}
\usepackage[mathscr]{eucal}
\usepackage{cancel}
\usepackage{wrapfig}
\usepackage[noadvisor]{frontespizio}
\usepackage{hyperref}

\addto\captionsitalian{}

\numberwithin{equation}{chapter}
\newtheorem{teo}[equation]{Teorema}

\newtheorem{lem}[equation]{Lemma}

\newtheorem{es}[equation]{Esempio}

\title{\Huge``Curve algebriche piane e sghembe''\\ \vspace{0,5 cm}\Large Corso del Prof. G. Castelnuovo 1922-23 } 
 \date{}
 
\begin{document}
\pagestyle{plain}
\begin{center}
\vspace{6cm}
\huge{Guido Castelnuovo: \emph{Curve algebriche \\ piane e sghembe} (Roma, 1922--23)}\\
\vspace{1cm}
\Large{Ciro Ciliberto e Claudio Fontanari}
\end{center}

\vspace{2cm}
 {\fontsize{4mm}{4.8mm}\selectfont Il quaderno, che qui presentiamo nell'accurata trascrizione di Natascia Zangani, \`e stato 
ritrovato da Ciro Ciliberto nel fondo Franchetta e contiene gli appunti di un anonimo uditore 
del corso \emph{Curve algebriche piane e sghembe}, tenuto da Guido Castelnuovo sulla cattedra 
di Geometria Superio\-re dell'Universit\`a di Roma nell'a.a. 1922--23. Per ovvie ragioni anagrafiche, 
l'estensore degli appunti non pu\`o essere identificato con Alfredo Franchetta (1916--2011).

L'importanza del quaderno, che a nostro parere ne giustifica pienamente la trascrizione e la 
pubblicazione, \`e da un lato la completezza e la precisione del testo, che riproduce scrupolosamente i 
dettagli delle lezioni di Castel\-nuovo e ne riflette il ben noto stile piano, elegante e rigoroso, dall'altro il carattere eccezionale, dal punto di vista storico, del corso dell'a.a. 1922--23. 
Si tratta infatti dell'ultimo corso di Geometria Superiore tenuto da Castel\-nuovo, che a partire 
dall'anno accademico successivo passer\`a sulla cattedra di Matematiche Complementari per 
lasciare il posto a Federigo Enriques, che si trasferiva a Roma dall'Universit\`a di Bologna. 
Un ulteriore motivo di pregio, storico e scientifico, del presente quaderno consiste nel fatto che, dei vari corsi superiori di argomento algebrico geometrico tenuti da Castelnuovo nella sua carriera, questo dell'a.a. 1922--23 \`e l'unico di cui sia giunto un testo 
dettagliato in una forma, diremmo, pubblicabile. E ci\`o \`e tanto pi\`u rimarchevole se si pensi che Castelnuovo non pubblic\`o mai un trattato di argomento algebrico geometrico. 

Grazie al prezioso lavoro di digitalizzazione curato da Paola Gario (cfr. \cite{G2}), gli appunti 
autografi di Castelnuovo per i corsi di Geometria Superiore, che consistono in brevi note molto sintetiche ma al contempo esaustive riguardo agli  argomenti trattati, sono attualmente disponibili in 
rete, anche se purtroppo non ancora trascritti. Anche solo una rapida scorsa ai titoli permet\-te di apprezzare la ricchezza dei temi 
trattati nel ventennio 1903--1923: \medskip

1903--04: \emph{Indirizzi geometrici}

1904--05: \emph{Sulle curve algebriche}

1905--06: \emph{Geometria sulle curve -- Curve algebriche sghembe}

1906--07: \emph{Sistemi lineari di curve piane -- Superficie razionali}

1907--08: \emph{Superficie razionali}

1908--09: \emph{Funzioni algebriche e loro integrali}

1909--10: \emph{Integrali abeliani e loro inversione}

1910--11: \emph{Geometria non-euclidea}

1911--12: \emph{Geometria differenziale}

1912--13: \emph{Geometria sulle curve algebriche} 

1913--14: \emph{Matematica di precisione e Matematica di approssimazione}

1913--14: \emph{Integrali abeliani}

1914--15: \emph{Calcolo delle probabilit\`a}

1915--16: \emph{Indirizzi geometrici}

1916--17: \emph{Geometria differenziale}

1917--18: \emph{Curve algebriche piane}

1919--20: \emph{Geometria non euclidea}

1920--21: \emph{Funzioni algebriche e loro integrali}

1921--22: \emph{Funzioni abeliane}

1922--23: \emph{Curve algebriche piane e sghembe}\medskip

Colpisce in particolare l'argomento del corso di Geometria Superiore del\-l'a.a. 1914--15, \emph{Calcolo delle probabilit\`a}. 
In questo caso, gli appunti autografi di Castelnuovo iniziano con le seguenti precisazioni: 

\begin{quote}\emph{Indole dei corsi di matem. 
super. Inopportunit\`a di una rigida divisione tra i vari rami di matem., ed anche tra la stessa matem. e le scienze 
naturali.} \end{quote}

Qui viene espresso con rapide note manoscritte un punto di vista profondamente radicato nel pensiero 
di Castelnuovo. Infatti, gi\`a nell'introduzione autografa al corso dell'a.a. 1904--05, dedicato alla geometria 
delle curve algebriche, si legge: 
\begin{quote}
\emph{Inopportunit\`a di dividere nettamente tra loro i vari rami della mate\-matica. 
La geometria degli enti algebrici studia in realt\`a questioni di algebra facendosi guidare in parte nella 
posizione dei problemi e nella ricerca dei teoremi e loro dimostrazioni dall'intui\-zione geometrica.}
\end{quote}

Come evidenziato in \cite{CG}, si tratta della stessa intima convinzione di Fede\-rigo Enriques, che in \cite{Enr} dichiara: 

\begin{quote}\emph{[...] la Geometria astratta si pu\`o identificare coll'Analisi. Appunto per ci\`o: Le due scienze debbono coltivarsi 
insieme. Non soltanto ne deriver\`a alla Geometria il vantaggio di una gene\-ralit\`a e di una potenza di metodi sperimentati 
ormai da Des Cartes in poi; ma ugualmente l'Analisi potr\`a essere indirizzata alle pi\`u belle scoperte dalla feconda 
intuizione geometrica, come Monge, Clebsch, Klein e Lie hanno insegnato. Per le cose dette apparir\`a naturale che nel 
seguito (pur riferendoci di pre\-ferenza ad interpretazioni geometriche intuitive) consideriamo indifferentemente gli enti 
definiti geometricamente o analiticamente.} 
\end{quote}

La fonte di questa visione unitaria della matematica che, nelle parole di 
Enriques, ambisce a 

\begin{quote}\emph{[...] collegare in un tutto organico diversi indirizzi e profittando dei metodi svariati tendere di 
regola con tutte le forze al resultato, il quale appartiene non ad un ramo della Matematica ma alla Matematica intera, 
se pure non si voglia porre in relazione al progresso della scienza in generale [...]}\end{quote} 
\`e indicata dallo stesso Enriques 
nel lavoro \cite{Seg} di Corrado Segre. 
Nelle \emph{Osservazioni dirette ai miei studenti} contenute in \cite{Seg}, il maestro torinese proclama: 

\begin{quote}\emph{[...] nel 
rivolgermi ai miei studenti di Geometria io sento il dover di raccomandar loro col massimo calore lo studio dell'Analisi.
Un giovane che voglia oggid\`i coltivare la Geometria staccandola nettamente dall'Analisi, non tenendo conto dei progressi 
che questa ha fatto e va facendo, quel giovane, dico, per quanto grande abbia l'ingegno, non sar\`a mai un geometra completo}.\end{quote}
Nello stesso spirito, Segre poi conclude: 
\begin{quote}\emph{Per ogni ricerca si scelga liberamente il metodo che sembra pi\`u opportuno; spesso converr\`a alternare fra loro il metodo sintetico che appare pi\`u penetrante, pi\`u luminoso, e quello analitico che in molti casi 
\`e pi\`u potente, pi\`u generale, o pi\`u rigoroso; e parecchie volte accadr\`a pure che uno stesso argomento non sar\`a bene 
illuminato sotto ogni aspetto se non sar\`a trattato con ambo i metodi.}\end{quote}

In questa prospettiva, non deve dunque stupire che gli appunti di Castel\-nuovo qui riprodotti inizino con brevi cenni 
di richiamo sulle funzioni di varia\-bile reale e complessa e proseguano applicando alle curve piane il metodo algebrico 
dell'eliminazione; che descrivano le propriet\`a locali delle curve con tecniche differenziali e poi analizzino
le molteplicit\`a di intersezione mediante l'algebra dei polinomi. Anche nelle lezioni pi\`u avanzate, Castelnuovo  
mette a frutto con disinvolta libert\`a metodi algebrici, analitici e proiettivi iperspaziali per dimostrare prima il cosiddetto 
\emph{Fundamentalsatz} di Noether, quindi il teorema di Riemann-Roch (per una traduzione moderna delle dimostrazioni di 
Segre e Castelnuovo rimandiamo a \cite{Ser}) e infine la sua stima sul genere di una curva sghemba, nota in letteratura 
come \emph{Castelnuovo's bound}. 

Per ulteriori informazioni su Castelnuovo e i suoi corsi rimandiamo a \cite{G1}; ci piace piuttosto riportare diffusamente 
da \cite{Par} l'esperienza diretta di O. Zariski, allievo del corso di Geometria Superiore dell'anno accademico precedente, 
il 1921--22, dedicato alle funzioni abeliane: 

\begin{quote}\emph{Standing in front of the class with his long black beard and quiet hands,
Castelnuovo often reminded Zariski of the Moses of Miche\-langelo, although
Zariski would also remember ``the sweet smile that suddenly transformed his face.''
His lectures, which were on analytic geometry during Zariski's first term at the
university [in the fall of 1921], were tightly structured, reflecting the formality 
of his manner. Zariski enjoyed them so thoroughly that it was almost a month before 
he realized that he was wasting his time.}

\emph{One day after class he found the courage to introduce himself and was relieved
when Castelnuovo, as cordially as his forbidding manner would allow, said, ``Come
with me. I am going home.'' As they made their way through the narrow streets,
Zariski explained that he knew all of analytic geometry and more calculus than had
been taught in his college because his study had been based on the French textbooks
on integral calculus (which were more like the analysis courses given in the first year
of graduate school). He told Castelnuovo how he'd been forced to enroll as a student
of philosophy at Kiev because the mathematics department had been full, and how
he'd studied only mathematics. He described the books that he'd read, including
Salmon's algebra and Goursat's calculus, and how he'd studied alone at home and
taken examinations secretly after they had been abolished by the Communists.}

\emph{Castelnuovo, who'd been quiet during all this time, suddenly began to ask mathematical
questions, simple ones at first, and then more and more advanced and searching
ones. By the time they stood in front of his house he seemed to have reached a
decision. ``Well, Zariski, you go tomorrow to the registrar's office and tell them that I
sent you there, and that I suggested that you should change your application, which
has already been accepted, in the following way: instead of asking admission to the
first year, you ask to be admitted to the third year. Then you come to my course.''
As Zariski himself liked to put it, on that short walk he gained two years. He also
gained a thesis advisor who would encourage his independence. Once Zariski complained
that he needed to know more about the functions of complex variables in
order to understand the abelian functions that Castelnuovo was covering in his
third-year course on algebraic geometry and algebraic functions. ``Go to the
library,'' Castelnuovo said. ``There are books and you can read them.''}\end{quote}

Queste note meriterebbero certo un'analisi critica approfondita ed un paragone, altrettanto critico, con 
corsi di contenuto algebrico geometrico dati dallo stesso Castelnuovo e da altri autori suoi contemporanei, e infine con i famosi testi \cite {EC, S}. Non essendoci possibile impegnarci nell'immediato in questo lavoro, abbiamo preferito mettere subito a disposizione degli studiosi questo interessante materiale, porgendo cos\`i anche un nostro modestissimo omaggio alla memoria di Guido Castelnuovo, di cui il prossimo 14 Agosto ricorrono i centocinquanta anni dalla nascita.

\hspace{0.5cm}

\begin{center}
\begin{tabular}[t]{@{}l}
Ciro Ciliberto \\
Dipartimento di Matematica \\
Universit\`a di Roma "Tor Vergata" \\
Via della Ricerca Scientifica 1 \\
00133 Roma, Italy. \\
E-mail: cilibert@mat.uniroma2.it
\end{tabular}\hfill
\begin{tabular}[t]{@{}l}
Claudio Fontanari \\
Dipartimento di Matematica \\
Universit\`a degli Studi di Trento \\
Via Sommarive 14 \\
38123 Trento, Italy. \\
E-mail: fontanar@science.unitn.it
\end{tabular}
\end{center}
}
%\noindent
%Ciro Ciliberto \newline
%Dipartimento di Matematica \newline
%Universit\`a di Roma "Tor Vergata" \newline
%Via della Ricerca Scientifica 1\newline
%00133 Roma, Italy. \newline
%E-mail: cilibert@mat.uniroma2.it

%\hspace{0.3cm}

%\noindent
%Claudio Fontanari \newline
%Dipartimento di Matematica \newline 
%Universit\`a degli Studi di Trento \newline 
%Via Sommarive 14 \newline 
%38123 Trento, Italy. \newline
%E-mail: fontanar@science.unitn.it

 \maketitle
 %----------------------------------------pag.1-------------------------%
\chapter*{Introduzione \\
\LARGE{Funzioni di variabile complessa - Proprietà geometriche e proprietà proiettive delle figure.}}
\addcontentsline{toc}{chapter}{Introduzione - Funzioni di variabile complessa - Proprietà geometriche e proprietà proiettive delle figure}

\subsection*{Funzioni di variabile reale.}
Definizione di Dirichlet: $y=f(x)$ in un certo intervallo $AB$ quando ad ogni valore di $x$ compreso in $AB$ corrisponde in un modo qualsiasi un valore di $y$.

\begin{figure}[!htbp]
\label{fig:1}
\centering
\includegraphics[scale=0.3]{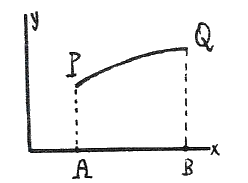}
\end{figure}

Se vogliamo che $y$ rappresenti geometricamente una curva, come si sa occorre fare delle restrizioni sulla continuità, ammissione di derivate prime. Finch\'{e} restiamo nel campo di funzioni di variabile reale nulla possiamo dire della prosecuzione della curva e proprietà fuori dei limiti $AB$ di definizione della curva. Analogamente per funzioni di $2$ variabili reali.

\subsection*{Funzione di variabile complessa.}
$y=f(x)$ risponda alle condizioni volute per essere sviluppabile in serie di Maclaurin:
\begin{gather}
\label{eq:varcomp}
f(x)=a_{0}+a_{1}x+a_{2}x^{2}+...
\end{gather}
Allora $f(x)$ è funzione analitica della variabile reale $x$ (Lagrange) e allora dando ad $x$ certi valori complessi la \eqref{eq:varcomp} è ancora convergente. Il campo di variabilità di $x$ perch\'{e} la \eqref{eq:varcomp} sia convergente si può rappresentare nel \emph{piano di Gauss} o della variabile complessa $x=\xi + i \eta$. Per ogni $x$ interno al \emph{cerchio di convergenza} $\sigma$ la \eqref{eq:varcomp} converge, per ogni $x$  esterno è divergente.

%----------------------------------------pag.2------------------------------%

\begin{figure}[!htbp]
\label{fig:1}
\centering
\includegraphics[scale=0.5]{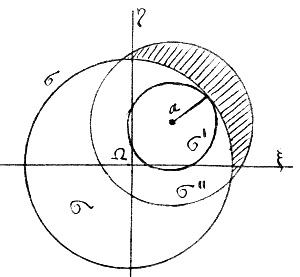}
\end{figure}
\noindent 
I punti reali $\eta=0$ per cui la \eqref{eq:varcomp} è convergente sono quelli interni al cerchio $\sigma$ posti sull'asse reale $\xi$.\\
Ammesso dunque che $f(x)$ sia una funzione analitica, da un campo reale essa resta anche definita in un campo complesso. Si può però andare avanti, e definire la funzione anche esternamente al detto cerchio di convergenza (Weierstrass): si prenda nell'interno di $\sigma$ un punto definito dal numero complesso $a$; si riesce a sviluppare la \eqref{eq:varcomp} in serie di potenze di $(x-a)$:
\begin{gather}
\label{eq:varcomp2}
f(x)=\alpha_{0}+\alpha_{1}(x-a)+\alpha_{2}(x-a)^{2}+ \dots
\end{gather}
Si dimostra che la \eqref{eq:varcomp2} è certamente convergente entro il cerchietto disegnato; e inoltre, in questo, \eqref{eq:varcomp} e \eqref{eq:varcomp2} danno lo stesso valore. Ma avviene in genere che la \eqref{eq:varcomp2} converge in un cerchio più ampio e concentrico di $\sigma'$. Abbiamo così esteso il dominio d'esistenza di $f(x)$ nel campo complesso. La \eqref{eq:varcomp2} è dedotta dalla \eqref{eq:varcomp} per \emph{prolungamento analitico} (Weierstrass).\\
E così ecc\dots si riesce a sviluppare la \eqref{eq:varcomp2} in una serie di potenze di $(x-b)$ ($b\equiv$ punto interno a $\sigma ''$) ecc \dots \\
Si ottengono così elementi successivi di funzione analitica, ciascuno dedotto dal precedente per prolungamento analitico. Proseguendo si ha una successione finita o infinita di sviluppi. 
%----------------------------------------pag.3------------------------------%
Le aree che si ottengono e che in parte si sovrappongono, costituiscono il dominio d'esistenza della funzione analitica o funzione di variabile complessa.\\
La funzione analitica è definita in tutto il suo dominio d'esistenza; anzi, se si conoscono i suoi valori nei punti d'un archetto del campo, allora la funzione è determinata completamente in tutto il suo dominio d'esistenza.

\subsection*{Curve analitiche e curve algebriche.}

Torniamo al piano $(x,y)$. Ora supponiamo $y$ funzione analitica di $x$, per modo che l'arco $PQ$ è un arco di \emph{curva analitica}. Allora intanto posso definire la curva per valori complessi dell'ascissa e dell'ordinata, e l'estensione dell'arco deve intendersi \emph{in un unico determinato modo} finch\'{e} $f(x)$ deve restare analitica. Si possono allora studiare non solo le proprietà differenziali, ma le proprietà della curva nella sua integrità (diametri, polarità, assi, asintoti, ecc \dots). Si dimostra che una funzione algebrica è una particolare funzione analitica. Così la Teoria delle curve algebriche si riattacca alla Teoria delle funzioni di variabile complessa. Se l'equazione della curva è di $n$- esimo grado, allora la funzione è ad $n$ rami, e si dimostra che si riattaccano l'uno all'altro per prolungamento analitico. La funzione $y$ è definita per valori reali 
%----------------------------------------pag.4----------------------------------%
e complessi di $x$; nella Teoria delle curve algebriche, non si fa distinzione tra punti reali e punti complessi della curva.

\subsection*{Proprietà geometriche delle figure.}
Il gruppo delle \emph{Similitudini} è il gruppo di trasformazioni per cui le proprietà geometriche delle figure (che con esso si definiscono) sono invarianti (Klein). \`{E} caratterizzato (Piani).  
%dubbio: incerta l'interpretazione di Piani, probabilmente intende quando ci troviamo nel piano
Trasforma: 
\begin{enumerate}
\item punti in punti
\item rette in rette
\item cerchi in cerchi.
\end{enumerate}
\noindent
I \emph{Rami di Geometria} studiano le proprietà delle figure invarianti di fronte a trasformazioni più generali delle similitudini. Per esempio il gruppo che trasforma punti in punti e cerchi in cerchi è il gruppo delle \emph{affinità circolari} (per esempio la proprietà di un cerchio di essere osculatore ad una curva in un punto, è invariante rispetto ad un'affinità circolare).

\subsection*{Proprietà Proiettive delle figure.}

Il gruppo che trasforma punti in punti e rette in rette definisce le \emph{Proprietà Proiettive}, ed è il gruppo delle \emph{collinearità} od \emph{omografie}.\\
Per esempio trasformazioni che cambiano una conica in una conica sono collineazioni; ma la proprietà per una conica di essere ellisse o iperbole o parabola, è proprietà geometrica.\\
Sono proprietà proiettive per un punto ed una retta l'essere polo e polare di una conica, ma  è proprietà geometrica quella di una punto di essere centro di una conica. Così 
%----------------------------------------pag.5-------------------------%
pure sono geometriche le proprietà focali, ecc \dots\\
Noi ci proponiamo lo studio delle \emph{Proprietà proiettive} delle curve algebriche piane e sghembe. Non ci può servire il sistema cartesiano di coordinate, ma il sistema proiettivo.

\chapter{Curve Algebriche piane}

%sezione 1%
%-------------------------------pag.5----------------------------------%

\section{Curve algebriche piane 
- Intersezioni con rette - Intersezioni di due curve - Teorema di B\'{e}zout.}
\label{sec:1}
\subsection{Curve algebriche piane d'ordine $n$.}
Sono rappresentate in coordinate cartesiane:
\begin{equation}
\label{eq:capo}
a+bx+cy+dx^{2}+exy+fy^{2}+\dots+px^{n}+\dots+ry^{n}=0.
\end{equation}
I coefficienti come pure le coordinate possono essere reali o complessi. Per avere la \eqref{eq:capo} in coordinate omogenee, basterebbe porre $x=\frac{X}{Z}$, $y=\frac{Y}{Z}$. \\
Nella \eqref{eq:capo} ci sono $\frac{(n+1)(n+2)}{2}$ termini; il numero dei \emph{parametri essenziali} è:

\begin{equation}
N=\frac{(n+1)(n+2)}{2}-1=\frac{n(n+3)}{2}
\end{equation} 
\noindent
Se la \eqref{eq:capo} deve passare per un punto ho un'equazione lineare non omogenea ad $N$ incognite (parametri); dati $N$ punti ho $N$ equazioni  non omogenee ad $N$ incognite (od omogenee ad $N+1$ incognite). Sicch\'{e} si prevede:
%---------------------------------pag.6---------------------------------%

\begin{center}
per $N$ punti di un piano passa sempre una $\mathcal{C}_{n}$ e generalmente una sola
\end{center}
\noindent
Ma può darsi il caso che qualcuna delle $N$ equazioni sia conseguenza delle rimanenti (può accadere  per particolari posizioni degli $N$ punti); allora per $N$ punti $\infty$ curve d'ordine $n$.\\
Ciò può accadere per posizioni qualunque degli $N$ punti? No. Prendiamo ad esempio una conica, e formiamo la matrice del sistema.
\[
\begin{pmatrix}
1 & x_{1} & y_{1} & x_{1}^{2} & x_{1}y_{1} & y_{1}^{2}\\
\dots & \dots & \dots & \dots & \dots & \dots\\
1 & x_{5} & y_{5} & x_{5}^{2} & x_{5}y_{5} & y_{5}^{2}\\
\end{pmatrix}
\]
Il caso di indeterminazione (pei $5$ punti infinite curve, cioè $\infty$ soluzioni del sistema delle 5 equazioni) si presenta quando si annullano \emph{tutti} i minori di quinto ordine della matrice. Si giungerebbe alla conclusione che dovrebbe essere $1=0$. Così se delle $N$ equazioni $K$ conseguissero dalle rimanenti $(N-K)$ si avrebbero, come si dice, $\infty^{K}$ curve per gli $N$ punti.\\

Esempio: per $n=3$, $N=9$; due cubiche s'intersecano in $9$ punti, cosicch\'{e} se per gli $N$ scegliamo proprio questi $9$, per essi passano $\infty$ cubiche; sappiamo anche che tutte le cubiche per $8$ di questi passano pel $9^{o}$; si sa che tutte le cubiche piane per $8$ punti, passano per un $9^{o}$ determinato dai primi.

\subsection{Intersezioni della \eqref{eq:capo} con la retta $y=\alpha x + \beta$.}
Intersezioni della \eqref{eq:capo} con la retta $y=\alpha x + \beta$; si hanno $n$ punti comuni, reali, immaginari, anche multipli. Fatta la sostituzione, la \eqref{eq:capo} potrebbe diventare di grado $(n-K)$, e
%-----------------------------pag.7----------------------------------------%
 si avrebbero $(n-k)$ intersezioni. Come spiegare ciò? Se il calcolo è fatto, come dovrebbe, in coordinate omogenee, la \eqref{eq:capo} risulta omogenea di grado $n$ in $x$ e $z$, e se non è soddisfatta da $z=0$ dà $n$ valori per $\frac{x}{z}$; ma può darsi anche che sia soddisfatta da $z=0$ e $x$ altri valori qualunque. Queste soluzioni non comparivano in coordinate non omogenee. Geometricamente: queste intersezioni stanno sulla retta all'infinito se adoperiamo coordinate cartesiane (retta $z=0$ in coordinate proiettive).

\subsection{Intersezioni della \eqref{eq:capo} con una curva algebrica d'ordine $m$.}

Vediamo prima la 
\subsubsection{Eliminazione di una incognita tra $2$ equazioni ad una sola incognita.}
\begin{equation}
\label{eq:cas}
\begin{cases}
a_{0}x^{m}+a_{1}x^{m-1}+ \dots +a_{m}=0 \qquad \text{ radici } x_{1}, x_{2}, \dots, x_{m}\\
b_{0}x^{n}+b_{1}x^{n-1}+\dots+b_{n}=0 \hspace{1.05cm} \text{ radici } y_{1}, y_{2}, \dots, y_{n}
\end{cases}
\end{equation}

Perch\'{e} le \eqref{eq:cas} abbiamo una radice comune, deve essere uguale a zero un certo polinomio in $a$ e $b$ (\emph{Risultante delle \eqref{eq:cas}}), la cui formazione 
costituisce la eliminazione dell'incognita. Si formino tutte le possibili differenze tra le prime radici e le seconde: dovrà essere zero il prodotto, che essendo una funzione omogenea (simmetrica) nelle radici, è esprimibile mediante i coefficienti. Cerchiamo dunque il Risultante delle \eqref{eq:cas} col metodo di Sylvester, pel caso di $m=2$, $n=3$; ma la dimostrazione è generale.\\
Dette \eqref{eq:cas} abbiano la radice comune $x_{0}$; si avrà:
%----------------------------------pag.8--------------------------------%
\[
\begin{cases}
\begin{matrix}
\qquad & \qquad & a_{0}x^{2}+ & a_{1}x+ &a_{2} &=0\\
\qquad & a_{0}x^{3} + & a_{1}x^{2}+ & a_{2}x \quad &  & =0\\
a_{0}x^{4}+ & a_{1}x^{3}+ & a_{2}x^{2} \quad & \qquad &  &=0\\
\qquad & b_{0}x^{3} + & b_{1}x^{2}+ &b_{2}x+ & b_{3} &=0\\
b_{0}x^{4} + & b_{1}x^{3}+ & b_{2}x^{2}+ & b_{3}x \quad &  &=0
\end{matrix}
\end{cases}
\]
Si hanno così $(m+n)$ equazioni, con la maggiore potenza $x^{m+n-1}$, che si possono riguardare lineari nelle incognite $x, x^{2}, \dots, x^{m+n}$. Questo sistema ha una soluzione $x=x_{0}, x^{2}=x_{0}^{2}, \dots$ dunque è zero il determinante dei coefficienti (\emph{Risultante}):
\[
R=
\begin{vmatrix}
a_{0} & a_{1} & a_{2} & 0 & 0\\
0 & a_{0} &  a_{1} & a_{2} & 0 \\
0 & 0 & a_{0} &  a_{1} & a_{2} \\
b_{0} & b_{1} & b_{2} & b_{3} & 0\\
0 & b_{0} & b_{1} & b_{2} & b_{3}\\
\end{vmatrix}=0
\]
Viceversa si dimostra che se $R=0$ le \eqref{eq:cas} hanno una radice comune.
Così, condizione necessaria e sufficiente perch\'{e} le \eqref{eq:cas} abbiano una radice comune, è $R=0$.
$R$ contiene omogeneamente i coefficienti di ciascuna di esse al grado dell'altra.
\subsubsection{Sistema di $2$ equazioni in $2$ incognite.} Posso ordinarle ad esempio rispetto ad $x$, ed ottengo due equazioni tipo \eqref{eq:cas}, in cui le $a$ e $b$ sono polinomi in y di grado uguale all'indice loro. Allora $R=0$ è una equazione in $y$, ogni cui radice dà luogo, sostituita nelle \eqref{eq:cas}, ad una radice $x$ comune. Di che grado è $R$? Se è $\nu$, $\nu$ sarà anche evidentemente il grado del nuovo determinante. 
\[
\begin{vmatrix}
1 & y & y^{2} & 0 &0\\
0 & 1  & y & y^{2} & 0\\
0 & 0 & 1  & y & y^{2}\\
 1  & y & y^{2} & y^{3} & 0\\
 0 &  1  & y & y^{2} & y^{3}\\
\end{vmatrix}
\]
\noindent
Se si modifica la seconda riga per $y$ si aumenta il grado di $1$, la $n$-esima per $y^{n-1}$ ancora di $(n-1)$; analogamente si opera pel secondo gruppo; poi dividendo la seconda colonna per 
%-----------------------------------pag.9---------------------------------%
$y$ si diminuisce il grado di $1$, la $(n+m)$-esima per $y^{n+m-1}$ ancora di $(n+m-1)$; il determinante diviene così numerico, cioè di grado zero in $y$, onde
\[
\nu +1+2+\dots+(n-1)+1+\dots+(m-1)-1-2-\dots-(m+n-1)=0
\]
cioè
\[
\nu=mn.
\]

\begin{teo}[Teorema di B\'{e}zout]
\label{}
Geometricamente: Due curve algebriche d'ordine $m$ ed $n$ hanno $mn$ intersezioni comuni.
\end{teo}
\noindent
Delle intersezioni potrebbero venire a coincidere, per contatti e punti multipli, e allora occorre verificare quante soluzioni vengono così assorbite.

\section{Comportamento d'una curva in un punto dato - Condizioni per un punto multiplo - Discriminante dell'equazione d'una curva. }
\label{sec:2}

\subsection{Comportamento nell'origine $O$.}
Curva per $O$ (manca il termine costante)
\begin{equation}
\label{eq:omo}
bx+cy+dx^{2}+exy+fy^{2}+\dots=0
\end{equation}
Almeno uno dei termini di $1^{o}$ grado non sia nullo in \eqref{eq:omo}.\\
 Retta generica per $O$:
\begin{equation}
\label{eq:gen}
y=\lambda x
\end{equation}
Per uno speciale valore di $\lambda$  
\[
y=-\frac{b}{c}x \qquad \text{(tangente in $O$),}
\]
si annullano in \eqref{eq:omo} i termini di $1^{o}$ grado, e la nuova \eqref{eq:omo} con la \eqref{eq:gen} hanno comune una radice doppia in $O$\footnote{Ch\'{e} si può raccogliere $x^{2}$ a fattor comune}. Se la curva tocca in $O$ l'asse $x$ la sua equazione è del tipo:
%-----------------------------------pag.10--------------------------------%
\begin{equation}
\label{eq:omo1}
cy+dx^{2}+exy+fy^{2}+gx^{3}+\dots=0
\end{equation}
Tangente: 
\begin{equation}
\label{eq:tan}
y=0
\end{equation}
Intersezioni: $$dx^{2}+gx^{3}+\dots =0;$$ se ne hanno due sole in $O$ se $d \neq 0$; ma se $d =0$, almeno $3$ (\emph{flesso}). Allora la \eqref{eq:omo1} diviene: equazione d'una curva che ha in $O$ con $y=0$ almeno un contributo di $2^{o}$ ordine
\[
y(c+cx+fy)+gx^{3}+\dots=0
\]
Eseguendo una rotazione d'assi coordinate si vede subito: per verificare se nella \eqref{eq:omo} la tangente in $O$ ($bx+cy=0$) ha ivi un contatto d'ordine superiore al $1^{o}$, basta esaminare se il gruppo dei termini di $1^{o}$ grado divide il gruppo dei termini di $1^{o}$ e $2^{o}$ grado. Si noti che allora, la conica ottenuta uguagliando a zero il gruppo di $1^{o}$ e $2^{o}$
 grado, che si vede ha in $O$ anche un contatto d'ordine superiore al $1^{o}$ con \eqref{eq:omo}, si spezza.\\

Se la equazione d'una curva comincia con termini di grado $r$ la curva ha un punto $r$-uplo in $O$, caratterizzato dalla proprietà che la retta generica per esso ha ivi $r$ punti comuni colla curva (poich\'{e} colla sostituzione ci si riduce ad una equazione nella sola $x$ con $r$ radici uguali a zero). Esistono però $r$ rette, \emph{Tangenti principali}, ciascuna delle quali ha almeno $r+1$ intersezioni colla curva (in $O$). L'equazione complessiva delle $r$ tangenti principali si ottiene uguagliando
%----------------------------------pag.11--------------------------------%
  a zero il gruppo dei termini di grado $r$ (equazione omogenea in $x$, $y$, che risolta rispetto ad $y$ dà $r$ radici di $1^{o}$ grado).
Se le tangenti sono distinte si ha un punto $r$-uplo ordinario, se no singolare.
Si vede chiaramente che:\\
perch\'{e} una curva abbia in un punto assegnato una molteplicità $r$ occorrono $$1+2+\dots+r=\frac{r(r+1)}{2}$$ condizioni lineari (coeff.=0). Si deduce come caso particolare: volendo una $\mathcal{C}_{n}$ passare con molteplicità $n$ per un punto, disponiamo solo di $n$ parametri.

\subsubsection{Intersezioni in $O$ di due curve, $\mathcal{C}_{m}$ e $\mathcal{C}_{n}$ - passanti per $O$ con molteplicità $r$ ed $s$.}
Col metodo di Sylvester 
%dubbio: nel manoscritto è scritto incorrettamente Sylvester (senza y)%
si può trovare il risultante in $y$. In esso il grado più basso in cui comparisce la $y$ è generalmente (come si vede con lunga dimostrazione) $rs$; dunque $rs$ radici nulle, cui corrispondono $rs$ radici  nulle per $x$: Dette curve hanno $rs$ intersezioni in $O$. Questo quando le $r$ tangenti principali a $C_{m}$, differiscono dalle $s$ a $C_{n}$. Se si ha qualche coincidenza si verificano contatti.

\subsection{Comportamento in un $P$ qualsiasi 
della curva $f(x,y)=0$, ossia $f(x_{0}+X,y_{0}+Y)$.}

Sviluppando secondo le potenze crescenti di $X$ e $Y$, l'equazione diviene\footnote{$f(x_{0}, y_{0})=0$, poich\'{e} la curva passa per $P$.}:
\begin{equation}
0=(x-x_{0})\frac{\partial f}{\partial x_{0}}+(y-y_{0})\frac{\partial f}{\partial y_{0}}+\frac{1}{2} \Biggl\{ (x-x_{0})^{2} \frac{\partial ^{2}f}{\partial x_{0}^{2}} + 2(x-x_{0})(y-y_{0})\frac{\partial ^{2}f}{\partial x_{0} \partial y_{0}} + \dots \Biggr\} + \dots
\end{equation}

\begin{figure}[!h]
\centering
\includegraphics[scale=0.3]{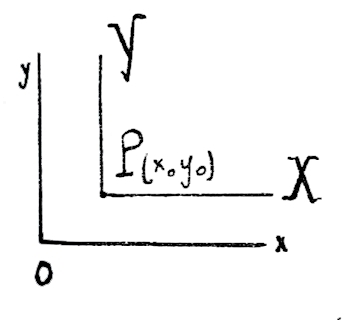}
%\caption{}
\end{figure}
\noindent
Ed ora si può ragionare come precedentemente. Per esempio, 
%----------------------------------pag.12---------------------------------%
se esiste uno dei termini a $1^{o}$ grado $P$ è punto semplice di $f$; equazione della tangente:
\[
(x-x_{0})\frac{\partial f}{\partial x_{0}}+ (y-y_{0}) \frac{\partial f}{\partial y_{0}}=0\\
\]
che in coordinate omogenee si riduce a
\[
x \frac{\partial f}{\partial x_{0}}+ y \frac{\partial f}{\partial y_{0}}+ z \frac{\partial f}{\partial z_{0}}=0
\]
(basta porre, per vederlo $z=z_{0}=1$, e tenere conto del Teorema d'Eulero:
\[
x_{0} \frac{\partial f}{\partial x_{0}}+\dots =n f(x_{0},y_{0},z_{0})=0;
\]
sottraendo questa dalla precedente si ottiene proprio l'antiprecedente).
L'ultima equazione permette anche di scrivere la tangente nel punto all'infinito.
\noindent
Condizioni in coordinate omogenee perch\'{e} $P$ sia doppio:
\[
\frac{\partial f}{\partial x_{0}}=0 \qquad \frac{\partial f}{\partial y_{0}}=0 \qquad \frac{\partial f}{\partial z_{0}}=0
\]
Così ecc. per punto triplo, ecc; si tiene conto delle identità d'Eulero, che non sono più una.

\subsection{Condizioni per un punti multiplo - Invarianti d'una funzione - Discriminante.}
Data la curva $f(x,y)=0$ come si verifica che essa ha un punto doppio in qualche posto? Deve essere 
$$\frac{\partial f}{\partial x_{0}}=0 \quad \frac{\partial f}{\partial y_{0}}=0 \quad \frac{\partial f}{\partial z_{0}}=0,$$ 
cioè caso generale: 
$$\varphi_{m}(x,y)=0 \, , \, \psi_{n}=0 \, ,  \, \chi_{p}=0$$ 
le cui radici danno il punto almeno doppio. 
Analogamente alla svolta Teoria delle eliminazioni si forma il Risultante. Le soluzioni delle ultime due siano $(x_{1},y_{1}), (x_{2}, y_{2}), \dots ,$ $ (x_{np}, y_{np})$. 
Dovrà essere 
\[
\varphi(x_{1},y_{1}) \cdot \varphi(x_{2},y_{2}) \dots  \varphi(x_{np},y_{np})=0
\]
Viceversa se questo prodotto è zero le equazioni precedenti hanno una soluzione comune. Questo prodotto (Risultante) contiene i coefficienti della $\varphi$ al grado $np$, ed è una funzione simmetrica delle coppie di soluzioni di $\psi$ e $\chi$. Si vede che: il Risultante di tre equazioni 
%-------------------------------pag.13---------------------------------%
con due incognite è una funzione razionale intera dei coefficienti delle tre equazioni, e contiene i coefficienti di ciascuna omogenei e ad un grado uguale al prodotto dei gradi delle due rimanenti.\\
Le equazioni in coordinate omogenee $\left(\frac{\partial f}{\partial x}=0, \frac{\partial f}{\partial y}=0, \frac{\partial f}{\partial z}=0\right)$ sono di grado $(n-1)$, quindi il Risultante (\emph{Discriminante di $f$}) contiene i coefficienti della prima al grado $(n-1)^{2}$, come della seconda e della terza. Condizione necessaria e sufficiente perch\'{e} $\mathcal{C}_{n}$ abbia un punto multiplo, è che sia nullo il \emph{Discriminante} della equazione (funzione razionale intera dei coefficienti dell'equazione, contenente questi al grado $3(n-1)^{2}$).\\
\noindent
Con operazioni razionali si riesce a costruire il discriminante. Per esempio per una conica
\[
a_{11}x^{2}+a_{22}y^{2}+a_{33}z^{2}+a_{12}xy+\dots=0
\]
le condizioni equazioni sono 
%dubbio: nel manoscritto è riportata la doppia dicitura equazioni condizioni

\begin{gather*}
a_{11}x+a_{12}y+a_{13}z=0,\\
a_{21}x+a_{22}y+a_{23}z=0,\\
a_{31}x+\dots=0
\end{gather*}
ed il discriminante sarà il determinante\footnote{Da esso si vede che tre coniche di un fascio si spezzano in rette.}:
\[
\begin{vmatrix}
a_{11} & a_{12} & a_{13}\\
a_{21} & a_{22} & a_{23}\\
a_{31} & a_{32} & a_{33}\\
\end{vmatrix}
\]

Se ho due curve del $3^{o}$ ordine $f=0$, $f_{1}=0$, quante curve del fascio $f+\lambda f_{1}=0$ hanno un punto doppio? Si vede subito: $\lambda$ compare nel discriminante al $12^{o}$ grado: $12$ curve hanno un punto doppio.

Ritorniamo alla curva $f=0$; sia il Discriminante $\Delta(a,b,c \dots)=0$, e sottoponiamo la $f$ ad una collineazione che è rappresentata da una trasformazione lineare:
%-----------------------------------pag.14--------------------------------%
\[
x= \alpha x'+ \beta y' +\gamma z' ; \quad y= \alpha_{1} x'+ \beta_{1} y' +\gamma_{1} z'
\]
Questa manda la $f$ in una nuova $f_{1}(x',y',z')$ il cui nuovo discriminante sarà 
 $\Delta (a',b',c', \dots)$. Se la $f$ ha un punto doppio, naturalmente anche la $f_{1}$ ne avrà uno:
\[
\Delta (a,b,c, \dots)=0\, \Longleftrightarrow \,\Delta (a',b',c', \dots)=0
\]
\noindent
Si dimostra che il rapporto dei $\Delta$ non può dipendere dai coefficienti di $f$, ma solo da $\alpha, \beta, \gamma, \dots$ della sostituzione lineare.\\
Invariante d'una funzione omogenea di più variabili è una espressione razionale intera dei coefficienti che non cambia valore salvo un fattore uguale ad una potenza (qualsiasi) del determinante della sostituzione, quando la espressione stessa venga calcolata coi coefficienti della funzione trasformata collinearmente.\\
Il Discriminante è un invariante della $f$. Una conica possiede un solo invariante in questo senso, ed è il discriminante. In una curva, $\Delta=0$ esprime una particolarità proiettiva.

\section{Tangenti da un punto; prima polare del punto - Classe di una curva con punti doppi e cuspidi - Curva reciproca - Formole di Poncelet.}
\label{sec:3}

\subsection{Tangenti da un punto.}
Data
\begin{equation}
\label{eq:func}
f(x,y,z)=0
\end{equation}
Equazione della tangente in $Q$:
\begin{equation}
\label{eq:tang}
X\frac{\partial f}{\partial x}+Y\frac{\partial f}{\partial y }+ Z\frac{\partial f}{\partial z}=0
\end{equation}
%---------------------------------pag.15----------------------------------%

\begin{figure}[!htbp]
\centering
\includegraphics[scale=0.5]{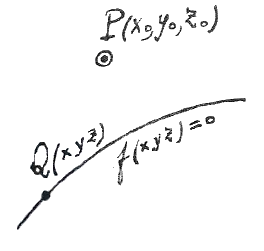}
%\caption{}
\end{figure}

Condizione perch\'{e} questa tangente passi per $P$:
\begin{equation}
\label{eq:cond}
x_{0} \frac{\partial f}{\partial x} + y_{0} \frac{\partial f}{\partial y} + z_{0}\frac{\partial f}{\partial z }=0
\end{equation}
Questa \eqref{eq:cond} ove si considerino variabili $x,y,z$, è una curva, la \emph{Prima Polare} del punto  rispetto alla \eqref{eq:func} (nome preso dalla Teoria delle coniche). Le $n(n-1)$ intersezioni di \eqref{eq:func} e \eqref{eq:cond} sono i punti di tangenza.
$n(n-1)$ è \emph{la classe di $f$}. Ma se per alcune intersezioni $\frac{\partial f}{\partial x}=\frac{\partial f}{\partial y}=\frac{\partial f}{\partial z}=0$ (Punti multipli di $f$ allora esistono), esse rendono indeterminata la equazione della tangente, \eqref{eq:tang}; esse soddisfano \eqref{eq:func} e \eqref{eq:cond} senza essere punti di tangenza. Così le intersezioni di \eqref{eq:func} e \eqref{eq:cond} cadono nei punti di tangenza e nei punti multipli. Si dimostra che la tangente in flesso conta per due; come è evidente per la tangente doppia; e come se $P$ è sulla curva.

\subsection{Classe di $f$ con punti doppi e cuspidi.}
$f$ abbia un punto doppio, che assumiamo come origine
\[
f=ax^{2}+2bxy+cy^{2}+dx^{3}+\dots=0
\]

\begin{figure}[!htbp]
\centering
\includegraphics[scale=0.6]{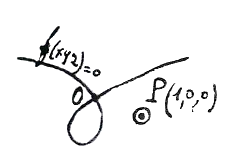}
%\caption{}
\end{figure}
\noindent
Introduciamo le coordinate omogenee, o prendendo $P(1,0,0)$ come vertice del triangolo fondamentale (coordinate proiettive), o proiettando la curva. 
%dubbio: nella scannerizzazione del manoscritto risulta tagliata la frase, forse ``in modo che P vada all'infinito.'' %
\\ 

Prima Polare di $P$\footnote{Questa è ottenuta con derivazione rapportata
%dubbio: l'interpretazione dell'abbreviazione rapp. è dubbia, forse ``rapportata''%
ad $x$ della coppia di tangenti in $O$ alla curva. Quindi è la polare di detta coppia (conica), ossia la loro coniugazione armonica rispetto a $P$. Essendo distinte, ne è distinta.}: 
$0=2(ax+by)+3dx^{2}+\dots$\\
\noindent
Se le tangenti Principali in $O$ ad $f$ sono distinte, allora la tangente in $O$ alla prima polare è distinta da esse. 
%----------------------------------pag.16--------------------------------%
Come si vede $O$ è multiplo secondo $2$ di $f$, e semplice per la prima polare. Un punto doppio assorbe due unità sulla classe della curva. Ora $f$ abbia una cuspide in $O$, asse $y$ tangente cuspidale:
\[
f=x^{2}+(ax^{3}+3b x^{2}y+3cxy^{2}+dy^{3})+ex^{4}+\dots=0
\]
\begin{figure}[!htbp]
\centering
\includegraphics[scale=0.6]{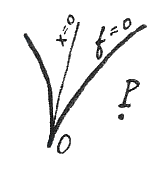}
%\caption{}
\end{figure}
\\Consideriamo anche qui $P(1,0,0)$; equazione della prima polare:
\[
\frac{\partial f}{\partial x}=2x+3(ax^{2}+2bxy+cy^{2})+4ex^{3}+\dots=0
\]
Ma qui la prima polare ha anch'essa l'asse $y$ per tangente in $O$, per cui essa passa semplicemente. Per liberarci da questa difficoltà, invece di considerare le radici in $O$ di $f=0$ e $\frac{\partial f}{\partial x}=0$, consideriamo le radici in $O$ di 
\[
2f-x\frac{\partial f}{\partial x }=0 \qquad \text{e} \qquad \frac{\partial f}{\partial x }=0,
\]
sistema equivalente. 
La prima $0=-ax^{3}+3cxy^{2}+2dy^{3}+2ex^{4}+\dots$ è una curva che ha nell'origine un punto triplo, in cui l'asse $y$ non è più tangente principale perch\'{e} $x$ non è fattore comune ai termini di $3^{o}$ grado. Si vede: una cuspide assorbe tre unità sulla classe della curva: quindi classe di $\mathcal{C}_{n}$ con $d$ punti doppi ordinari e $r$ cuspidi:
\begin{equation}
\nu=n(n-1)-2d-3r
\end{equation}
Questa è la $1^{a}$ formola di Poncelet.

\subsection{Curva reciproca di $f$ - Formole di Poncelet.}
La Correlazione è una corrispondenza che muta punti in rette, punti allineati in rette per un punto. Se ad ogni punto si fa corrispondere la polare si ha la Polarità. Così se $P$ descrive  $f$, la polare 
%dubbio: le parole conica e retta vengono aggiunte forse in matita sugli appunti%
descriverà
%---------------------------------pag.17--------------------------------%
 la \emph{Reciproca di $f$}, $F$. Viceversa ad ogni punto di $F$ corrisponde una tangente di $f$. Analiticamente, la più generale correlazione si scrive: $u, v, w$ (coordinate di rette) funzioni lineari di $x,y,z$.
%dubbio: qui probabilmente manca z a causa della fotocopiatura%
Per esempio se nell'un piano 
%dubbio: incerta l'interpretazione 
abbiamo la curva 
\begin{equation}
\label{eq:a}
f(x,y,z)=0
\end{equation}
ed eseguiamo la trasformazione $x=u, y=v, z=w$, veniamo a scrivere 
\begin{equation}
\label{eq:b}
f(u,v,w)=0
\end{equation}
che rappresenta la reciproca di \eqref{eq:a} in coordinate di rette: 
una stessa equazione rappresenta una curva o la sua reciproca secondoch\'{e} le variabili si riguardano coordinate di punti o di rette.\\
Data la \eqref{eq:a}, se si vuole la reciproca \eqref{eq:b} in coordinate di punti: costruisco la tangente $t$ in un punto della \eqref{eq:a}:  $$X\frac{\partial f}{\partial x}+Y\frac{\partial f}{\partial y }+ Z\frac{\partial f}{\partial z}=0;$$ le coordinate di rette di $t$ sono $$u=\frac{\partial f}{\partial x}, ~ v= \frac{\partial f}{\partial y}, ~ w=\frac{\partial f}{\partial z }.$$
 Questa $t$ corrisponde ad un punto
 \begin{equation}
 \label{eq:c}
 \xi=\frac{\partial f}{\partial x},~ \eta = \frac{\partial f}{\partial y}, ~\zeta =\frac{\partial f}{\partial z }
 \end{equation}
Al variare di $(x,y,z)$ sulla \eqref{eq:a}, varia $P(\xi,\eta,\zeta)$ e descrive la curva reciproca $F$. Sicch\'{e} $F$ è rappresentata parametricamente $(\xi,\eta,\zeta)$ dalle \eqref{eq:c} legate dalla \eqref{eq:a}. Sicch\'{e} per l'equazione di $F$ in coordinate $(x,y,z)$  basta eliminare $\xi,\eta,\zeta$ fra le \eqref{eq:c} e la \eqref{eq:a}. $F$ è d'ordine $n(n-1)$ quando $f$ non ha punti multipli, altrimenti l'ordine si abbassa. Ricordando la definizione di ordine e di classe si 
%----------------------------------------pag.18--------------------------%
 vede subito: date due curve reciproche $f$ ed $F$, l'ordine dell'una è uguale alla classe dell'altra.\\
\noindent
Ad un punto doppio ordinario di una curva, corrisponde nella reciproca una tangente doppia. Ad una cuspide una tangente d'inflessione. Così $f_{n}$ abbia classe $\nu$, $d$ punti doppi ordinari, $r$ cuspidi, $\delta$ tangenti doppie, $\rho$ flessi. Sono evidenti le due formole di Poncelet
\begin{gather}
\label{eq:P1}\nu = n(n-1)-2d-3r\\
\label{eq:P2}n=\nu(\nu-1)-2\delta-3\rho
\end{gather}
(Ma potrebbero aversi curve con particolarità più complicate). Se $f$ manca di punti doppi e cuspidi, $F$ deve averne. E qui si spiega il paradosso di Poncelet: esempio cubica $d=r=0$: $\nu=6$, quindi sembra $n=30$. Questa difficoltà aveva fatto nascere il sospetto che fosse in una curva qualunque $n=\nu$. Ma in realtà, se $d$ ed $r$ sono zero, non lo saranno $\delta$ e $\rho$, e così si riabbassa a tre. Possono darsi solo tre casi;
\[
\begin{matrix}
n & d & r & & \nu & \delta & \rho\\
3 & 0 & 0 &  & 6 & 0 & 9\\
3 & 1 & 0 & \Rightarrow & 4 &0 &3\\
3 & 0 & 1 & & 3 & 0 & 1
\end{matrix}
\]
non può essere $d=2$; la retta avrebbe $4$ intersezioni, ecc. Per la stessa ragione è sempre $\delta=0$.\\
Per le quartiche piane, colle \eqref{eq:P1} e \eqref{eq:P2} troveremmo un numero esuberante di casi.\\

%--------------------------------------pag.19------------------------------%

\section{Flessi di una curva - Curva Hessiana - Formole di Pl\"{u}cker.}
\label{sec:4}

\subsection{Flessi d'una curva piana Hessiana.}
Sappiamo che condizione perch\'{e} $P_{0}(x_{0}, y_{0},z_{0})$ sia flesso della curva per $P_{0}$:
\begin{equation}
\label{eq:uno}
0= X\frac{\partial f}{\partial x_{0}}+Y\frac{\partial f}{\partial y_{0} }+ Z\frac{\partial f}{\partial z_{0}}+\frac{1}{2} \left( X^{2}\frac{\partial^{2} f}{\partial x_{0}^{2}}+2XY\frac{\partial^{2} f}{\partial x \partial y }+\dots   \right)
\end{equation}
è che la conica dei termini di $1^{o}$ e $2^{o}$ grado si spezzi; come si vede, se la conica (per $P_{0}$) si spezza, tangente in $P_{0}$ ad essa sarà $$X\frac{\partial f}{\partial x_{0}}+Y\frac{\partial f}{\partial y_{0} }=0,$$ e lo spezzamento consta di questa retta ed un'altra (non per $P_{0}$; sarebbe punto doppio) cioè il gruppo dei termini a $1^{o}$ grado divide il gruppo dei termini a $2^{o}$ grado, e allora la curva primitiva ha un flesso in $P_{0}$. Condizione necessaria e sufficiente perch\'{e} la curva primitiva abbia un flesso in $P_{0}$ è che la conica detta, semplice in $P_{0}$, sia degenere, cioè dev'essere zero il discriminante della conica:
\begin{equation}
\label{eq:due}
\begin{vmatrix}
\frac{\partial^{2} f}{\partial x_{0}^{2}} & \frac{\partial^{2} f}{\partial x_{0} \partial y_{0}} & \frac{\partial f}{\partial x_{0}}\\
\frac{\partial^{2} f}{\partial x_{0} \partial y_{0}} & \frac{\partial^{2} f}{\partial y_{0}^{2}} & \frac{\partial f}{\partial y_{0}}\\
\frac{\partial f}{\partial x_{0}} & \frac{\partial f}{\partial y_{0}} & 0
\end{vmatrix}=0
\end{equation}
In coordinate omogenee:

\begin{equation}
\label{eq:due1}
\begin{vmatrix}
\frac{\partial^{2} f}{\partial x_{0}^{2}} & \frac{\partial^{2} f}{\partial x_{0} \partial y_{0}} & \frac{\partial^{2} f}{\partial y_{0} \partial z_{0}}\\
\frac{\partial^{2} f}{\partial y_{0} \partial x_{0}} & \frac{\partial^{2} f}{\partial y_{0}^{2}} & \frac{\partial^{2} f}{\partial y_{0} \partial z_{0}}\\
\frac{\partial^{2} f}{\partial z_{0}\partial x_{0}} & \frac{\partial^{2} f}{\partial z_{0} \partial y_{0}} & \frac{\partial^{2} f}{\partial z_{0}^{2} }
\end{vmatrix}=0
\end{equation}
\noindent
Possiamo partire dalla \eqref{eq:due1} per ritrovare la \eqref{eq:due}, tenendo conto del Teorema d'Eulero sulle funzioni omogenee\footnote{La \eqref{eq:uno} in coordinate omogenee è $X^{2} \frac{\partial^{2} f}{\partial x_{0}^{2} }+Y^{2}\frac{\partial^{2}y}{\partial y_{0}^{2}}+Z^{2}\frac{\partial^{2} f}{\partial z_{0}^{2}}+ 2XY \frac{\partial^{2} f}{\partial x_{0} \partial y_{0}}+2XZ \frac{\partial^{2} f}{\partial x_{0} \partial z_{0} }+2YZ \frac{\partial f}{\partial y \partial z}$. Infatti si ricordi il Teorema d'Eulero, si sottragga, e si ottiene la \eqref{eq:uno}.}.
Moltiplichiamo l'ultima verticale di \eqref{eq:due1} per $z$, ed aggiungiamo agli elementi quelli rispettivi delle prime due colonne moltiplicate per $x$ ed $y$. Eseguiamo 
%----------------------------------------pag.20-----------------------------%
la stessa operazione sull'ultima orizzontale di \eqref{eq:due1}; tenendo conto che \eqref{eq:due1}=0, si ottiene precisamente \eqref{eq:due}. Però \eqref{eq:due1} è più comoda, simmetrica, e tiene conto di eventuali punti all'infinito.\\ 
Il Determinante \eqref{eq:due1}, Jacobiano delle tre derivate prime di $f$, è l'\emph{Hessiano} di $f$, e la \eqref{eq:due1} è la \emph{Curva Hessiana} di $f$, d'ordine $3(n-2)$. Come si vede ogni flesso di $f$ sta sulla Hessiana; e viceversa ogni punto semplice di $f$ che stia sulla Hessiana è flesso su $f$. Le intersezioni di $f$ colla Hessiana cadono nei flessi e nei punti multipli di $f$ (che annullano le derivate prime e seconde). Ogni punto del piano ha la prima polare; $\infty^{2}$ 
prime polari d'ordine $(n-1)$; una generica prima polare non ha punti multipli all'infuori dei punti sulla $f$; ma ve ne sono $\infty^{1}$ ciascuna delle quali possiede un punto doppio; il luogo di questi è la Hessiana:  la Hessiana è il luogo dei punti doppi delle $\infty^{1}$ prime polari che hanno un punto doppio fuori di quelli della curva. Ogni intersezione di $f$ con $H$ è un flesso. Anzi si vede che un flesso conta per una sola intersezione di $f$ con $H$. Dunque il numero di flessi di $f_{n}$ è $3n(n-2)$. Vediamo come si abbassa il numero dei flessi, $\rho$, quando $f$ ha punti doppi ordinari e cuspidi, per cui passa $H$, e che assorbono delle intersezioni di $f$ con $H$, numero che vogliamo determinare. Si vedrebbe che un punto doppio di $f$ è doppio per $H$, una
%----------------------------------------pag.21---------------------------%
 cuspide di $f$ è tripla per $H$: evitiamo la discussione. In ogni punto doppio di $f$ cadono $\alpha$ intersezioni con $H$, cuspide $\beta$.
 %dubbio: la lettera beta è tagliata nella scannerizzazione del manoscritto%

\subsection{Formole di Pl\"{u}cker.}
Si avrebbe 
\begin{equation}
\label{eq:tre}
\rho=3n(n-2)-\alpha d-\beta r
\end{equation}
Moltiplichiamo la \eqref{eq:P1} per $3$ e sottragghiamovi la \eqref{eq:tre}: si avrà 
\[
(6-\alpha)d+(9-\beta)r-\rho=3(n-\nu),
\]
e dualmente: 
\[
( 6-\alpha)\delta+(9-\beta)\rho-r=3(\nu-n).
\]
Sommando membro a membro: 
 \begin{equation}
 \label{eq:qua}
 (6-\alpha)(d+\delta)+(8-\beta)(r+\rho)=0
 \end{equation}
Questa sarebbe soddisfatta per $\alpha=6$ e $\beta=8$ da qualunque curva. Se le differenze non sono nulle, la \eqref{eq:qua} dà $d+\delta=K(r+\rho)$ ($K=cost$).\\
\`{E} possibile che per ogni curva algebrica passi una tal relazione?\\
Prendiamo una generale\footnote{Per curva \emph{generale} s'intende una curva senza punti multipli. \\ Per curva \emph{generica} s'intende una curva qualunque.\\
Questa denominazione si giustifica: una curva generica non può avere punti multipli che non in casi eccezionali, come si può facilmente vedere tenendo conto delle condizioni.}
%dubbio: nel manoscritto riportato in cima alla pagina e con un asterisco%
 $d=r=0$.
%dubbio: nel manoscritto poco leggibile%
Dovrebbe essere $\delta=K\rho$ per qualunque $n$. La seconda formola di Poncelet diverrebbe $n=(\nu-1)\nu-(2k+3)\rho$, cioè $$n=n(n-1)\left[n(n-1)-1\right]-(2K+3)\cdot3n\cdot(n-2).$$
Una uguaglianza di questo genere dovrebbe sussistere, qualunque $n$, per $K$ fisso; ma questa non può essere un'identità, è assurda in generale, è assurdo dunque $\alpha\neq 6$ e $\beta\neq 8$. Abbiamo le formole, che possiamo chiamare anche tutte di Pl\"{u}cker:
\begin{enumerate}[label=\emph{\roman*}), leftmargin=1.5cm]
\item $\nu=n(n-1)-2d-3r$ \label{eq:FP1}
\item $n=\nu(\nu-1)-2\delta-3\rho$ \label{eq:FP2}
\item $\rho=3n(n-2)-6d-8r$ \label{eq:FP3}
\item $r=3\nu(\nu-2)-6\delta-8\rho$ \label{eq:FP4}
\end{enumerate}
Non sono indipendenti, ed è perciò che abbiamo potuto calcolare $\alpha$ e $\beta$; così da \ref{eq:FP1} e \ref{eq:FP3} che da \ref{eq:FP2} e \ref{eq:FP4} possiamo avere $3 \nu -\rho=3n-r$.\\

La \ref{eq:FP1}, \ref{eq:FP2}, \ref{eq:FP3}, \ref{eq:FP4} fissano legami necessari tra i $6$ caratteri d'una curva,
%dubbio: nel manoscritto d. curva%
ma non è detto che tutte le loro soluzioni aritmetiche, anche
%------------------------------------pag.22------------------------------%
positive diano luogo a curve esistenti. Alcune soluzioni vanno scartate (es. se il genere venisse negativo, la curva non esiste,
 %dubbio: esisterebbe %
  oppure potrebbe trattarsi di una curva spezzata), e questo scartamento è facile per le classi basse.\\
 Proiettando una $\mathcal{C}_{n}$ e la sua $H$ su un altro piano, vengono un'altra $\mathcal{C}_{n}'$ colla sua $H'$. Ciò condurrebbe a dire che il determinante Hessiano ha colla equazione della curva $\mathcal{C}$ una relazione in qualche senso analoga a quella vista negli invarianti. L'Hessiano ha i coefficienti funzioni di $a,b,c$. Si trasformi dunque collinearmente $\mathcal{C}$ in $\mathcal{C}'$; $a,b,c,\dots$ si trasformano in $a',b',c', \dots$. $H'$ non differisce da $H$  che per un fattore, potenza del determinante della sostituzione, indipendente dai coefficienti e dalle variabili della $\mathcal{C}$. Cioè l'Hessiano è covariante dell'equazione della curva. Esempio: conica; l'Hessiano diviene il discriminante, e l'Hessiano invariante diviene covariante. 

\section{Fasci di curve - Relazioni tra i punti base.}
\label{sec:5}

\subsection{Fasci di curve piane.}
Date $f_{n}$ e $\varphi_{n}$, la curva $f+\lambda \varphi=0$ ($\lambda$ parametro) al variare di $\lambda$ descrive un ``fascio'' di curve di grado $n$; $f$ e $\varphi$  non hanno importanza particolare nel fascio. \`{E} chiaro che ogni punto comune  
%----------------------------------------pag.23----------------------------%
ad $f$ e $\varphi$ è comune a tutte le curve del fascio (\emph{Punto base}). \\
Se due curve d'un fascio passano in un punto con molteplicità $r$, tutte le altre del fascio hanno ivi un punto $r$-uplo. Se però le due prime hanno anche comune le $r$ tangenti principali, fra le curve del fascio ve ne sarà una ed una sola avente in quel punto una molteplicità superiore ad $r$.\\
Ciò si vede immediatamente; la seconda parte si verifica precisamente quando nelle due curve i coefficienti dei termini di grado $r$ sono proporzionali: $$a:b:c: \dots= a':b'c':\dots$$\\
Viceversa, se $f$ ha molteplicità $r$ in $0$, e $\psi$
%dubbio: nel manoscritto qui viene riportato Psi, ma secondo me è phi?%
$(r+1)$, tutte le altre curve del fascio hanno ivi molteplicità $r$, e con le stesse tangenti principali.\\

\begin{teo}
\label{}
Per ogni punto del piano che non sia base, passa una ed una sola curva del fascio determinato da $f_{n}$ e $\varphi_{n}$.
\end{teo}
Infatti:
%dubbio: nel manoscritto riportato il primo pezzo con un asterisco, io ho deciso di inserirlo, altrimenti si può riportare come nota a piè di pagina, ma sembra parte integrante della dimostrazione%
Sia dato il punto $P_{1}(x',y')$; sostituendo queste coordinate in $f+\lambda \varphi=0$, ricaviamo $\lambda$. Dunque intanto per $P_{1}$ passa certamente almeno una curva del fascio. Si consideri ora una irriducibile $\psi_{n}$ che passi per le $n^{2}$ intersezioni di $f$ e $\varphi$, e precisamente in modo d'avere con ciascuna curva del fascio in ciascun punto base tante intersezioni quante due curve del fascio. Allora $\psi$ fa parte del fascio. Infatti: consideriamo la curva del fascio che passa per un generico punto di $\psi$.
\begin{figure}[!htbp]
\centering
\includegraphics[scale=0.4]{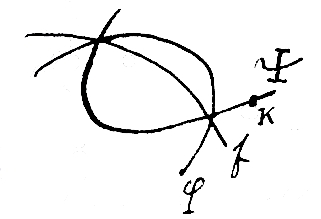}
%\caption{}
\end{figure}
%----------------------------------------pag.24-------------------------%
\\
Allora questa curva del fascio ha con la $\psi$ $(n^{2}+1)$ intersezioni, ma fuori dei punti base; essa deve coincidere allora con la $\psi$ in tutto o in parte. Quindi $\psi$ essendo irriducibile fa parte del fascio (\emph{Principio di Lam\'{e}}). Enunciamo alcune conseguenze. La proprietà che per ogni punto passa una sola curva del fascio serve a caratterizzare completamente il fascio, tra i sistemi algebrici di curve algebriche. Spieghiamo: sia un'equazione algebrica che contenga razionalmente un solo parametro essenziale, $k$: $$f(x,y,k)=0.$$ Supponiamo che per un punto generico $(x,y)$ passi una sola curva di questo sistema semplicemente infinito.
%dubbio: forse da rendere con semplicemente infinito%
Allora si dimostra che il sistema è un fascio, cioè la curva è $f(x,y,\lambda)=0$, lineare in $\lambda$.\\
Noi supporremo sempre si tratti di punti base tutti semplici. Due curve d'ordine $n$  si seghino in $n^{2}$ punti base semplici per ambedue. Allora, per $n^{2}$ punti base passano $\infty^{1}$ curve d'ordine $n$, e questo sebbene $n^{2}\geq \frac{n(n+3)}{2}$ appena $n > 3$.
%dubbio: nel manoscritto n>1 ma sembra corretto a matita in n>3 che è il risultato della disequazione...%
Ciò non è una contraddizione, ma vuol dire che degli $n^{2}$ punti base solo $\frac{n(n+3)}{2}-1$ sono indipendenti; ed ogni curva per essi $N-1$ passerà pei rimanenti $\frac{(n-1)(n-2)}{2}$.\\
%dubbio: forse l'errore n>1 abbia delle conseguenze sui calcoli successivi da ricontrollare%
\noindent
Presentiamo il risultato altrimenti: le curve d'ordine $n$ per $N-1$ punti generici formano un fascio, e si vede così: considerando una $\mathcal{C}_{n}$ per essi ho $N-1$ equazioni lineari non omogenee ad $N$ incognite
%----------------------------pag.25-----------------------------------------%
(parametri); ad una incognita ($\lambda$) dà un valore arbitrario. Resta un sistema di $N-1$ equazioni lineari con $N-1$ incognite che risolvo. Se applico la Regola di Cramer mi accorgo che ciascuno dei coefficienti da determinare dipende linearmente da $\lambda$, in modo che l'equazione della curva generica per gli $N-1$ punti sarà della forma:
\[
(a+a'\lambda)x^{n}+(b+b'\lambda)x^{n-1}y+\dots=0
\]
Se separo qui i termini contenenti $\lambda$ da quelli non contenentilo si vede che ho la equazione d'una curva che al variare di $\lambda$ dà un fascio. Tutte le curve d'ordine $n$ per $N-1$ punti generici del piano, passano in conseguenza per altri $\frac{(n-1)(n-2)}{2}$ determinati dai primi.\\
Per esempio tutte le cubiche per $8$ punti passano per un $9^{o}$ determinato dai primi $8$. Le coordinate del $9^{o}$ punto devono dipendere linearmente da quelli dei primi $8$. Cioè con la sola riga dati $8$ punti si può determinare il $9^{o}$. Per $n=4$ il problema è di $3^{o}$ grado, non basta n\'{e} la riga n\'{e} il compasso.

\section{Applicazioni alle cubiche: configurazione dei flessi - Fasci di cubiche - Fascio determinante da cubica ad Hessiana, e flessi.}
\label{sec:6}

\begin{teo} 
Se degli $n^{2}$ punti base di un fascio, $nm$ appartengono ad una curva irriducibile $\psi$ d'ordine $m$, i rimanenti $n(n-m)$ stanno su una curva d'ordine $(n-m)$.
\end{teo}
%----------------------------pag.26-----------------------------------------%
\begin{figure}[!htbp]
\centering
\includegraphics[scale=0.5]{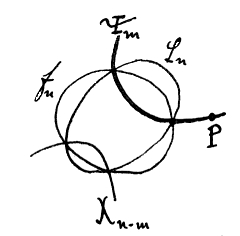}
%\caption{}
\end{figure}
\noindent
Tra le curve del fascio $f+\lambda \phi=0$ posso considerarne una per $P$; questa ha colla $\psi$ $mn+1$ intersezioni e coincide con essa 
%dubbio: nel manoscritto la frase successiva è inserita in parentesi quadre
[un'altra non la può segare che in $nm$ punti]. La nostra curva del fascio si spezza nella $\psi$ ed in una $\chi_{(n-m)}$ che deve contenere gli altri punti base.\\
Per un valore di $\lambda$ che conduca allo spezzamento
\[
f+\lambda \varphi\equiv \psi \chi \qquad f=\psi \chi -\lambda \varphi \qquad \text{ e viceversa.}
\]
\noindent
Applichiamo il teorema alle cubiche:\\
se delle $9$ intersezioni di $2$ cubiche tre sono allineate, le rimanenti $6$ appartengono ad una conica, e viceversa. Segue immediatamente il Teorema di Pascal sulle coniche. I tre lati costituiscono un trilatero che è una particolare curva del $3^o$ ordine. Si scrive scrivendo le equazioni dei tre lati e moltiplicandole. Una seconda curva del $3^o$ ordine si ha alternando i vertici.

\begin{figure}[!htbp]
\centering
\includegraphics[scale=0.8]{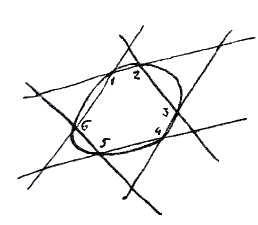}
%\caption{}
\end{figure}
\noindent
Queste due curve determinano un fascio del $3^o$ ordine. I $9$ punti base sono: $6$ intersezioni su una conica, e gli altri tre devono stare su una retta. Appartengono a $1-2$ con $4-5$, $2-3$ con $5-6$ e $3-4$ con $6-1$ (punti d'incontro dei tre lati opposti dell'esagono).\\
\begin{figure}[!htbp]
\centering
\includegraphics[scale=0.6]{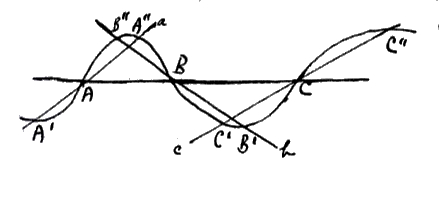}
%\caption{}
\end{figure}

Prendiamo una cubica piana schematica; prendiamo una retta secante in $3$ punti.
Rette che la 
%----------------------------pag.27-----------------------------------------%
incontrano ciascuna in altri due punti, $a, b, c$ formano una curva del $3^o$ ordine, e i $9$ punti disegnati sono base del fascio; $A, B, C$ sono allineati, gli altri $6$ staranno su una conica; se altri  tre sono allineati, la conica si spezza in $A'B'C'$ e $A''B''C''$. Considerando che $A',B',C'$ si vadano avvicinando indefinitamente ad $A, B, C$: le tangenti ad una cubica in $3$ punti allineati, la incontrano di nuovo in $3$ punti allineati (i \emph{tangenziali} di $3$ punti allineati sono anch'essi allineati). Segue:

\subsection{Flessi di una cubica piana.}
La retta che sega la cubica in $2$ flessi, la segherà in un terzo flesso. La configurazione dei flessi di una cubica piana si compone di $9$ punti (flessi) e $12$ rette d'inflessione, ciascuna delle quali contiene \emph{tre}flessi dei $9$, mentre per ciascuno dei punti passano $4$ di queste rette (ma dei $9$ flessi $3$ sono reali, $6$ immaginari).\\
Si possono con le $12$ rette d'inflessione costruire $4$ triangoli in modo che sui lati di ciascuno si trovino situati complessivamente i $9$ flessi (\emph{Triangoli d'inflessione}).
\begin{figure}[!htbp]
\centering
\includegraphics[scale=0.5]{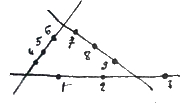}
%\caption{}
\end{figure}
Evidentemente per ogni retta d'inflessione resta determinato un triangolo d'inflessione. Ma così detti triangoli vengono contati tre volte. Quindi si hanno $\frac{12}{3}=4$ triangoli d'inflessione. Se consideriamo due dei tre triangoli d'inflessione, i $6$ lati s'incontrano
%----------------------------pag.28-----------------------------------------%
complessivamente in $9$ flessi. Ogni triangolo è una curva degenere del $9^o$ ordine; sicch\'{e} i $9$ flessi sono punti base d'un fascio di cubiche; ciò che sapevamo. In questo particolare fascio di cubiche (determinato dalla cubica e sua Hessiana) vi sono $4$ cubiche degeneri in $3$ rette, ciò che non succede per un fascio di curve generali. I tre vertici d'ogni triangolo sono punti doppi sicch\'{e} i $12$ punti doppi sono i $12$ vertici dei triangoli d'inflessione. Altra proprietà del fascio: \\
i $9$ punti base, sono flessi per tutte le cubiche del fascio determinato dalla $\mathcal{C}$ ed $H$, anche per la Hessiana.\\
\begin{figure}[!htbp]
\centering
\includegraphics[scale=0.6]{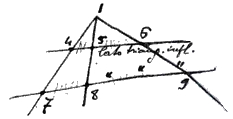}
%\caption{}
\end{figure}
\noindent
Da $\underline{1}$ proiettiamo $4,5,6$; otteniamo $9$ rette, ciascuna deve contenere un ulteriore flesso (nel caso reale non si può costruire). Consideriamo una cubica del fascio, passerà per tutti questi punti (come per $2$ e $3$); segata da $\underline{456}$ in $3$ punti, da $\underline{789}$ in altri tre. Unisco $4-7, 5-8, 6-9$. Queste tre rette vanno a incontrare la cubica in $3$ punti allineati che vengono a riunirsi in $\underline{1}$. Cioè $\underline{1}$ è flesso per ognuna delle cubiche che passa per quei $7$ punti. Così formando il fascio della cubica colla Hessiana, tutte le curve del fascio hanno i $9$ punti base come flessi.Segue: ciascuna di queste cubiche ha la sua Hessiana nel fascio stesso. Segue ancora: vogliamo determinare algebricamente le coordinate dei $9$ flessi di una cubica. Dobbiamo intersecare la $\mathcal{C}$ con la $H$, che porta ad un'equazione di $9^o$
%----------------------------pag.29-----------------------------------------%
grado in $x$. Ma questa risoluzione si può abbassare, come vedremo ora, ed essere risolta con equazioni di $3^o$ e $4^o$ grado, con metodi algebrici (estrazione di radice da funzioni razionali dei coefficienti; equazioni di grado  $>4$
%dubbio: ho sciolto >4 in maggiore di 4%
non si possono risolvere con radicali).\\
Formo il fascio $\mathcal{C}\, H$, e so che $4$ curve degenerano in triangoli; vuol dire che devonsi avere $4$ valori di $\lambda$ che spezzano $\mathcal{C}+\lambda H$ in $3$ rette. Per averle: uguaglio a $0$ il discriminante (di $\mathcal{C}+\lambda H$), equazione di $12^o$ grado; deve allora dare $4$ radici triple, si riduce quindi ad un'equazione di $4^o$ grado. Fin qui tutto razionalmente. Trovo così uno dei $4$ triangoli d'inflessione ($\mathcal{C}+\lambda_{1} H=0$),
%dubbio: nel manoscritto sembra che ci sia scritto \lambda_1  ha senso? %  
equazione di $3^o$ grado; ciascuno dei lati incontra la $\mathcal{C}$ in $3$ flessi. Sicch\'{e} intersecando $\mathcal{C}$ coi tre lati detti, trovo i $9$ flessi.

\section{Sistemi lineari di curve di dimensione $r$ - Teoremi - Punti multipli e punti base - Sistemi di curve degeneri - Teorema di Bertini.}
\label{sec:7}

\subsection{Sistemi lineari di dimensione $r$.}
Siano $(r+1)$ curve piane
\begin{gather}
\label{eq:sis1}
f_0(x,y)=0 \qquad f_1(x,y)=0 \qquad \dots \qquad f_r(x,y)=0
\end{gather}
Si consideri la combinazione lineare:
\begin{gather}
\label{eq:sis2}
\lambda_0f_0+\lambda f_1+\dots+\lambda_r f_r=0
\end{gather}
%---------------------------pag.30-----------------------------------------%
Potrebbe accadere che per valori particolari delle $\lambda$ la \eqref{eq:sis2} divenisse una identità; allora alcune delle $f$ sono combinazioni lineari delle altre; cioè allora la \eqref{eq:sis2} si potrebbe esprimere come combinazioni lineare di meno di $(r+1)$ polinomi. Se invece le \eqref{eq:sis1} sono tutte linearmente indipendenti, al variare dei rapporti di $r$ parametri al rimanente (parametri essenziali), la \eqref{eq:sis2} descrive un sistema $\infty^r$ di curve. Le \eqref{eq:sis1} non hanno nessuna particolarità nel sistema, ma possono venire sostituite da altre $(r+1)$ curve del sistema linearmente indipendenti. Per questo: si prendano $(r+1)$ curve \eqref{eq:sis2}; se qualcuna è combinazione lineare delle altre, il determinante dei parametri è zero, e viceversa. Per avere curve linearmente indipendenti il determinante non si deve annullare. Allora possiamo ritrovare tutte le curve del sistema $\infty^r$: un sistema $\infty^r$ è determinato da $(r+1)$ curve arbitrarie linearmente indipendenti.\\
Se i parametri sono legati da $K$ relazioni lineari indipendenti, la \eqref{eq:sis2} rappresenta un sistema $\infty^{r-K}$ di curve. Se $K=r$, si arriva ad una sola curva.\\
Tutte queste proprietà si ricordano meglio se si fa corrispondere alla \eqref{eq:sis2} il punto dello spazio ad $r$ dimensioni ($S_r$) di coordinate omogenee $\lambda_0 \lambda_1 \dots \lambda_r$; quando queste variano
%----------------------------pag.31-----------------------------------------%
il punto varia entro $S_r$. Se i parametri ($\lambda$) sono legati da una relazione (lineare omogenea), il punto varia entro un iperpiano ($S_{r-1}$); un numero conveniente di iperpiani si intersecano in un piano, con uno di più in una retta, con uno di più in un punto ($r$ iperpiani). Arrivati al punto, per l'interpretazione detta possiamo riferirci a curve (corrispondenti ai punti).\\
Nella Geometria Proiettiva degli iperspazi si hanno forme generatrici di $r^{a}$ specie ($\infty^r$ elementi); se l'elemento generatore è il punto, si ha $S_r$ punteggiato, il cerchio il sistema $\infty^r$ di cerchi (in $S_r$), la curva d'ordine $n$ il sistema lineare $\infty^r$ di curve d'ordine $n$ (in $S_r$). La Geometria proiettiva non ha bisogno di precisare la natura dell'elemento generatore.
Tutte le proprietà proiettive delle forme di $r^{ma}$ specie si trasportano senz'altro ai sistemi lineari $\infty^r$ di curve piane.\\
Si viene dunque a costruire la Geometria Proiettiva degli iperspazi senza specificare l'elemento generatore. Così si vengono a stabilire solo le proprietà che dipendono esclusivamente dal modo come i parametri entrano nella combinazione lineare, ma non dipendono invece dalla natura delle curve considerate. \\
Se un punto è comune a due curve del sistema, con 
%----------------------------pag.32-----------------------------------------%
molteplicità $i$, sarà comune a tutte le curve del sistema (\emph{Punto base}) con molteplicità $i$; si può avere molteplicità superiore per curve speciali. Per esempio scegliamo una curva del sistema che passi doppiamente per $O$.
Devono mancare i due termini di $1^o$ grado (coefficienti zero), quindi si vengono ad avere due equazioni lineari omogenee nei $\lambda$. 
Prendiamo un punto non base. Considero le curve del sistema $\infty^r$ per esso: una condizione lineare omogenea tra le $\lambda$, quindi dette curve formano un sistema lineare $\infty^{r-1}$.
E così: per $r$ punti in posizione generica, una ed una sola curva del sistema $\infty^r$.
Questa proprietà può servire a definire un sistema lineare di dimensione $r$: supponiamo avere una equazione contenente algebricamente un certo numero di parametri; al variare di questi si descrive un sistema algebrico di curve. Supponiamo 
sapere che per $r$ punti generici passa una ed una sola curva del sistema. Allora il sistema ``algebrico'' è naturalmente un sistema lineare; cioè i parametri o entrano linearmente nella equazione, o possono sostituirsi con nuovi che vi entrino linearmente. Cioè:
\begin{teo} 
Un sistema algebrico di curve $\mathcal{C}_n$ tale che per $r$ punti passa una ed una sola curva è lineare.
\end{teo}
\noindent
Non lo dimostriamo.\\
Per esempio: tutte le $\mathcal{C}_n$ del piano formano un sistema lineare di dimensione $\frac{n(n+3)}{2}$. 
%----------------------------pag.33-----------------------------------------%
\begin{teo} 
Una curva generica d'un sistema lineare non può avere punti multipli fuori dei punti base del sistema.
\end{teo}
\noindent
Per un sistema \emph{continuo}
 qualsiasi di curve algebriche o trascendenti, semplicemente $\infty$, $\lambda$ parametro
\begin{equation}
\label{eq:1sis3}
f(x,y,\lambda)=0
\end{equation}
Al variare di $\lambda$ la \eqref{eq:1sis3} descrive un sistema continuo di curve.
Allora resta determinato l'inviluppo del sistema che è il luogo dei punti d'incontro di curve infinitamente vicine, ossia la curva toccata da tutte le curve del sistema. Se una generica possiede un punto doppio variabile con $\lambda$, questo descrive una certa linea. \underline{Questa linea fa parte della curva inviluppo}. Così pei punti multipli. Così, un punto comune a due curve è punto base, e quindi un punto doppio d'una curva è punto base. E così è dimostrato il Teorema. Vediamo come si stabilisce il \underline{lemma} sottolineato. Per ogni $\lambda$, coordinate del punto multiplo:
\begin{equation}
\label{eq:coor}
x=\varphi(\lambda) \qquad y=\psi(\lambda)
\end{equation}
Sostituendo nella \eqref{eq:1sis3}, la \eqref{eq:1sis3} è soddisfatta per qualunque $\lambda$. Allora possiamo derivarla rispetto a $\lambda$:
\begin{equation}
\label{eq:deriv}
\frac{\partial f}{\partial \lambda}=f'_\lambda(x,y,\lambda)=0
\end{equation}
%----------------------------pag.34-----------------------------------------%
La \eqref{eq:deriv} è verificata sostituendo nella \eqref{eq:1sis3} le \eqref{eq:coor}. Per determinare l'inviluppo della \eqref{eq:1sis3} basta eliminare $\lambda$ tra \eqref{eq:1sis3} e \eqref{eq:deriv}.\\
Essendo queste verificate dalle \eqref{eq:coor}, vuol dire che le coordinate \eqref{eq:coor} verificano anche l'equazione dell'inviluppo. Il punto \eqref{eq:coor} dunque per qualunque $\lambda$ appartiene alla curva inviluppo.
Consideriamo una curva generica
\begin{equation}
\label{eq:sis3b}
f(x,y)+\lambda \varphi (x,y)=0
\end{equation}
\noindent
La generica \eqref{eq:sis3b} possieda un punto doppio variabile con $\lambda$. Allora le coordinate 
%dubbio: nel manoscritto vengono indicate con (4')...come riportarlo con la numerazione di qui?%
di questo punto multiplo, soddisfano la \eqref{eq:sis3b}, e anche la derivata
\begin{equation}
\label{eq:derivb}
\varphi(x,y)=0
\end{equation}
Quel punto multiplo appartiene dunque a due curve, \eqref{eq:sis3b}, \eqref{eq:derivb}, quindi è un punto base. Solo in questi possono cadere punti multipli della curva generica del fascio; non ci sono punti multipli variabili. Punti base variabili si presentano in curve degeneri. Esempio: coniche degeneri con una retta fissa.
I punti della curva fissa sono comuni a tutte le curve del fascio.

\begin{figure}[!htbp]
\centering
\includegraphics[scale=0.6]{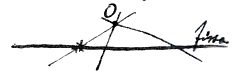}
%\caption{}
\end{figure}

\subsection{Sistemi lineari di curve tutte degeneri o riducibili.}
Ve ne sono due tipi ed uno misto.
%dubbio: incerta la resa con itemize, forse meglio enumerate avendo come label[1^o Tipo:]%
\begin{itemize}
\item \emph{$1^o$ Tipo}. La curva generica si compone di una parte fissa per
%----------------------------pag.35-----------------------------------------%
tutte, $\varphi=0$, ed una variabile in un sistema lineare:
\[
\psi=\lambda_0 f_0+\dots + \lambda_r f_r=0
\]
La curva $\varphi \cdot \psi=0$ costituisce un sistema lineare di curve tutte riducibili.
\item \emph{$2^o$ Tipo}. Si consideri un polinomio omogeneo di grado $n$ in $x,y$: $F(x,y)=0$. Si vede ch'esso si riduce ad $F\left(\frac{x}{y}\right)=0$ cioè ad $n$ rette 
\[
\biggl\{\frac{x}{y}=k_1, \dots, \frac{x}{y}=k_n \biggr\}
\]
per $O$.
Si consideri l'involuzione d'ordine $n$ e dimensione $r$: $$\lambda_0 F_0 +\dots + \lambda_r F_r=0;$$ si hanno lo stesso $n$ rette, ma variabili coi parametri: $\infty^r$ gruppi.\\

\emph{Genericamente}\\
Essendo $\varphi$ e $\psi$ due polinomi in $(x,y)$, ed $F$ omogenea in $\varphi, \psi$, di grado $n$; consideriamo $F(\frac{\varphi}{\psi})=0$ (a cui si riduce $F(\varphi, \psi)=0$); si hanno \underline{$n$} curve 
\[
\biggl\{\frac{\varphi}{\psi}=K_1, \dots, \frac{\varphi}{\psi}=K_n \biggr\}
\]
del fascio $\varphi=0$, $\psi=0$. Si consideri l'$\text{involuzione}^r_n$:
\[
\lambda_0 F_0+\dots \lambda_r F_r=0;
\]
si ha lo stesso una curva composta di $n$ curve, una variabile coi parametri restando in uno stesso fascio: $\infty^r$ gruppi di curve, che costituiscono un sistema lineare di curve tutte spezzate.
\item  \emph{$3^o$ Tipo}. Misto.
\end{itemize}

\begin{teo}[Teorema di Bertini.]
Oltre i due tipi visti di sistemi lineari di curve riducibili, e al tipo misto, non ve n'è nessun altro tipo.
\end{teo}
\noindent
Se la curva generica d'un sistema lineare è riducibile, o si compone di una parte fissa per tutte e una variabile irriducibile, o di più curve variabili nello stesso fascio, o i due fatti si presentano insieme.
%dubbio:sembra presente un rimando a matita al paragrafo seguente ma risulta illeggibile da fotocopia%
%----------------------------pag.36-----------------------------------------%
\section{Teorema fondamentale di  N\"{o}ther sulle curve passanti per le intersezioni di due curve date.} 
\label{sec:8}

\subsection{Sistema completo di curve o determinato dai punti base.}
\`{E} l'insieme di \emph{tutte} le curve d'ordine $n$ che passano nel punto $a_{1}$  con molteplicità $\alpha_1$, \dots $a_m$ con molteplicità $\alpha_m$. Bisogna ben assegnare i punti con le molteplicità, perch\'{e} \emph{oltre questi} potrebbero comparire altri punti base (Es. cubiche per $8$ punti).\\
Le curve d'ordine $n$ dipendono da $\frac{n(n+3)}{2}$ parametri; la molteplicità $\alpha$ porta $\frac{\alpha(\alpha+1)}{2}$ condizioni lineari; restano disponibili parametri
\[
r=\frac{n(n+3)}{2}-\sum\frac{\alpha(\alpha+1)}{2}
\]
%dubbio: sum o sigma? sembra \sum_i=1^m \frac{\alpha_i(\alpha_i+1)}{2}%
purch\'{e} le condizioni imposte dai punti base siano indipendenti, altrimenti restano disponibili parametri $r'=r+s$; $r\equiv$ dimensione \emph{virtuale}, $r'\equiv$ dimensione \emph{effettiva} del sistema. Se $s=0$ (condizioni indipendenti) il sistema è \emph{Regolare}; se no \emph{Sovrabbondante}, $s$ è la sovrabbondanza (Es. fascio di cubiche per $8$ punti sistema regolare, per $9$ sovrabbondante, $s=1$).\\
Ora supponiamo dati i punti e le molteplicità; posso considerare per essi curve $\mathcal{C}_n$
d'ordini crescenti; se $n$ è troppo basso non si avranno curve per essi. 
Può darsi che pei primi valori di $n$ a partire dai quali il sistema esiste, il sistema riesca sovrabbondante. Vedremo che appena $n$ supera un certo limite, il sistema diviene regolare.\\
%----------------------------pag.37-----------------------------------------%
Quante condizioni un gruppo di punti presenta alle curve di dato ordine che passano con date molteplicità per essi (s'intende per valori relativamente bassi di $n$). 
La risposta si può dare in pochi casi. Per esempio, se il gruppo è dato dall'intersezione di due curve.

\begin{teo} 
La equazione d'ogni $\mathcal{C}_n$ che passi semplicemente per le intersezioni supposte semplici di $\varphi_p$ e $\psi_q$ si può porre sotto la forma
\[
f=A\varphi+B\psi=0
\]
dove $A$ e $B$ sono due polinomi che introdotte le coordinate omogenee, sono di grado $A_{n-p}$ e $B_{n-q}$. (Il Reciproco è evidente). 
\end{teo}

\begin{enumerate}[label=\emph{\Roman* $^a$ Parte}, leftmargin=1.5cm]
\item Riguarda il caso se, ad esempio, $\psi$ è una retta: $\psi\equiv y=0$.\\
\begin{figure}[!htbp]
\centering
\includegraphics[scale=0.6]{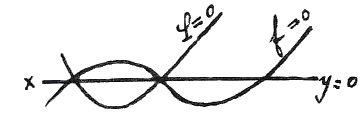}
%\caption{}
\end{figure}

$f=0$ incontra $y=0$ nel gruppo $f(x,0)=0$\\
$\varphi=0$ incontra $y=0$ nel gruppo $\varphi(x,0)=0$\\
Deve essere per ipotesi $f(x,0)=A(x)\varphi(x,0)$. Formiamoci $$F(x,y)=f(x,y)-A(x)\varphi(x,y).$$ Se pongo $y=0$, il polinomio $F$ si annulla, dunque è divisibile per $y$\\ $F=y\cdot B(x,y)$. E allora si deduce subito il nostro Teorema
\[
f(x,y)=A(x)\varphi(x,y)+B(x,y)y
\]
\item Se il Teorema è vero per un certo valore di $n$, è vero per valori inferiori. Dimostriamo che se è vero per curve d'ordine $n+1$ è vero per curve d'ordine $n$. Traccio una retta che non passi pei punti, $y=0$. Allora $f_n=0$ e $y=0$ danno una $C_{n+1}$ che passa per le
%----------------------------pag.38-----------------------------------------%
intersezioni di $\varphi $ e $\psi$. Allora sarà ($A$ e $B$ polinomi):
\begin{equation}
\label{eq:2parte1}
yf\equiv A\varphi +B \psi
\end{equation}
Pei punti intersezioni di $y=0$ e $\psi=0$, essi devono stare evidentemente anche sulla $A\varphi=0$; ma non stanno sulla $\varphi=0$, dunque staranno sulla $A=0$. Allora per la $I^a$ parte $A\equiv A'y+B'\psi$, e quindi 
\begin{equation}
\label{eq:2parte2}
yf\equiv (A'y+B'\psi)\varphi+B\psi\equiv A'y\varphi+(B'\varphi+B)\psi
\end{equation}
$y$ dovrà dividere l'ultimo termine, ma non divide $\psi$, dovrà dividere l'altro fattore: $B'\varphi+B=B''y$; allora
\[
f\equiv A'\varphi+B''\psi
\]
\item Il Teorema sussiste per le curve d'ordine $n$ abbastanza alto.
Prendo $n$ tanto alto che le $pq$ intersezioni di $\varphi$ e $\psi$ presentino alle curve d'ordine $n$ condizioni tutte indipendenti (ciò che è possibile pel \underline{lemma}). %enunciato alla fine della seconda pagina precedente
%dubbio: qui vale la pena inserire un ref con la citazione del lemma%
Vedremo ora che per $n$ così alto il Teorema è vero, cioè si avrà 
\begin{equation}
\label{eq:3parte3}
f_n\equiv A_{n-p}\varphi_p+B_{n-q}\psi_q=0
\end{equation}
Dimensione delle curve $\mathcal{C}_n$
 per le $pq$ intersezioni di $\varphi$ e $\psi$ $\left[\left[n\right]=\frac{n(n+3)}{2} \right]$:
\[
r=\left[ n\right] -1 -pq
\]
Dimensione delle curve $\mathcal{C}_n$ che si possono porre sotto la forma \eqref{eq:3parte3}, che pel reciproco del Teorema passano per dette $pq$ intersezioni:
\begin{equation}
\label{eq:3parte4}
r'=\left[n-p\right]+\left[n-q\right]-1-\nu=d-\nu
\end{equation}
%----------------------------pag.39-----------------------------------------%
Sia $\nu$ il numero dei parametri arbitrari di cui posso disporre per mettere $f$ sotto la forma $A\varphi +B\psi$. Vedremo $r=r'$, e così il Teorema sarà dimostrato. Data una curva $f$ rappresentabile sotto la forma \eqref{eq:3parte3}, sono dunque dati\\ $d-\nu$ parametri, ne restano $\nu$ arbitrari, cioè \emph{una stessa} curva $f$ rappresentabile sotto la forma \eqref{eq:3parte3}, si può rappresentare sotto $\infty^\nu$ forme. Ora vediamo di trovare $\nu$. \emph{Una} $f$ sia rappresentabile:
\begin{equation}
\label{eq:3parte5} 
A\varphi+B\psi\equiv A'\varphi+B'\psi
\end{equation}
Vediamo con che arbitrarietà variano $A'$ e $B'$; si ottiene
\[ 
\left. 
\begin{aligned} 
A'\equiv & A-\psi X_{n-p-q}\\
B'\equiv & B+\varphi X_{n-p-q}
\end{aligned} 
\right\} 
\quad
A' \varphi+B'\psi\equiv(A-\psi X)\varphi+(B+\varphi X)\psi
\]
Come si vede, $A'$ e $B'$ dipendono \emph{dai coefficienti} d'un polinomio di grado $n-p-q$, cioè da $\left[n-p-q\right]$ parametri. Dunque ci occupiamo di $n\geq p+q$. Allora la $f_n$ \emph{data} si può rappresentare in $\infty^{\left[n-p-q\right]}$ modi sotto la forma \eqref{eq:3parte3}, $\nu=\left[n-p-q\right]$.\\ 
Sostituendo si ha $r=r'$.
\end{enumerate}

\subsection{Sistema completo di curve $\mathcal{C}_n$ per le intersezioni di $\varphi$ e $\psi$.}
La \eqref{eq:3parte4}, per $n\geq p+q$, (ove sia $\nu=\left[n-p-q\right]$) ci dà la dimensione effettiva delle curve $\mathcal{C}_n$ per le $pq$ intersezioni di $\varphi$ e $\psi$. Se $n>p+q$\footnote{Allora la rappresentazione è unica. Infatti, dalla \eqref{eq:3parte5} $(A-A')
%dubbio: nel manoscritto riportato A con il pedice, ma sembra un errore%
\varphi=(B-B')\psi$; così $\varphi$ deve dividere $B'-B$; ma se $n<p+q$ è di grado minore, e ciò è impossibile; cioè $A'=A, B'=B$.},
%dubbio: nella nota non è molto chiara la costruzione della frase...%
da $r'$ si deve sopprimere $\nu$, cioè per avere la dimensione effettiva basta contare i parametri della \eqref{eq:3parte3}.
%dubbio: nel manoscritto dubbio il riferimento alla 3 o alla csi%
Così per ogni valore di $n$ sappiamo trovare la dimensione effettiva del sistema. Così se $n\geq p+q$ la dimensione effettiva è la $r'$ al completo, ed $r'=r$; cioè le
%----------------------------pag.40-----------------------------------------%
condizioni imposte dalle intersezioni di $\varphi$ e $\psi$ sono indipendenti appena sia $n\geq p+q$. Anzi vedremo anche indipendenza per $n=p+q-1$ ed $n=p+q-2$, ma non per valori inferiori.\\
Dunque, se $n<p+q$ la rappresentazione di \emph{una} $f_n$ è unica a meno di una costante moltiplicativa, ed i parametri sono in numero:
\begin{equation}
\label{eq:alpha}
r=\left[n-p\right]+\left[n-q\right]-1
\end{equation}
Se $n\geq p+q$ (arbitrarietà di $X$)
\begin{equation}
\label{eq:beta}
r=\left[n-p\right]+\left[n-q\right]-1- \left[n-p-q\right]=\frac{n(n+3)}{2}-pq
\end{equation}
Quindi la \eqref{eq:alpha} la posso scrivere 
\[
r=\frac{n(n+3)}{2}-pq+\left[n-p-q\right]
\]
Il sistema sarebbe regolare senza l'ultimo termine. La sovrabbondanza è:
\[
\left[n-p-q\right]=\frac{(n-p-q+1)(n-p-q+2)}{2}=\frac{(p+q-n-1)(p+q-n-2)}{2}
\]
Risulta: le curve $\mathcal{C}_n$ che passano semplicemente per le $pq$ intersezioni semplici di $\varphi_p$ e $\psi_q$ formano un sistema regolare (sovrabbondanza $=0$) se $n\geq p+q-2$, sovrabbondante per valori inferiori di $n$. 

\begin{teo}[Principio di Lam\'{e}]
Se $\varphi, \psi, f$ hanno lo stesso ordine $n$, le curve $\varphi, \psi, f$ per le stesse $n^2$ intersezioni formano fascio; sicch\'{e} la equazione d'una di esse si può scrivere come combinazione lineare delle altre due; allora $A$ e $B$ sono costanti.
\end{teo}

\subsection{Teorema di N\"{o}ther.}
\begin{teo} [Teorema di N\"{o}ther]
Se $f=0$ è l'equazione d'una curva per le $pq$ intersezioni di $\varphi $ e $\psi$, con molteplicità $\alpha+\beta-1$ per un punto comune a $\varphi_{p, (\alpha)}, \psi_{q, (\beta)}$ (che vi passano rispettivamente con molteplicità $(\alpha)$ e $(\beta))$,
%----------------------------pag.41-----------------------------------------%
allora $f$ può porsi $$f=A\varphi+B\psi$$ 
dove $A=0$ è curva pel detto punto con molteplicità $\beta -1$\\
dove $B=0$ è curva pel detto punto con molteplicità $\alpha -1$. 
\end{teo}
\noindent
Ci basta esaminare il caso semplice che non avvengano contatti pei punti multipli detti (Tangenti principali delle due curve distinte). Anche qui tre parti:

\begin{enumerate}[label=\emph{\Roman* $^a$ Parte}, leftmargin=1.5cm]
\item Se, ad esempio $\psi=0$ è retta, la dimostrazione è la stessa del Teorema precedente \footnote{In cui si tenga conto della molteplicità. Dunque per ipotesi $f$ passa con molteplicità \emph{$\alpha$} 
%dubbio: sottolineare alpha?%
per ogni punto in cui $\varphi$ passa con molteplicità $\alpha$, ed $y=0$ naturalmente con molteplicità $1$. Allora $F$ che appartiene al fascio ($A$ molteplicità $0$) passa anche ivi $\alpha$ volte, ed $f_{(\alpha)}\equiv A_{(0)}\varphi_{(\alpha)}+B_{(\alpha-1)}y_{(1)}$.}.
\item Se il Teorema è vero per le curve d'ordine $n+1$, sarà vero per l'ordine $n$.\\
Si ripete la dimostrazione precedente $yf\equiv A_{(\beta-1)}\varphi+B_{(\alpha-1)}\psi$. Si giunge: $A\equiv A'y+B'\psi$. Consideriamo la curva $A'y\equiv A-B'\psi$; appartiene al fascio $A=0$~ $B'\psi=0$. Ora $A=0$ ha molteplicità $\beta-1$, $\psi B'=0$ almeno molteplicità $\beta$ nelle intersezioni di $\varphi$ e $\psi$; la curva generica del fascio allora, molteplicità almeno $\beta-1$; sicch\'{e} $A'y=0$ ha almeno molteplicità $\beta-1$. Verrà come al Teorema precedente
\[
f\equiv A'\varphi+B''\psi
\]
$B''\psi=0$ appartiene al fascio $f=0$ $A'\varphi=0$. Nel punto che ci interessa, $f$ ha molteplicità $\alpha+\beta-1$, $A'\varphi$ ~ $\alpha+\beta-1$, sicch\'{e} ogni curva del fascio, come $B''\psi$, almeno molteplicità $\alpha +\beta -1$. Ma $\psi$ vi passa con molteplicità $\beta$, dunque $B''$ molteplicità $\alpha -1$.
%----------------------------pag.42-----------------------------------------%
\item Il Teorema sussiste per $n$ abbastanza alto che le curve d'ordine $n$ per le $pq$ intersezioni di $\varphi$ e $\psi$ formino un sistema regolare.
%dubbio: costruzione della frase?%
 Allora possiamo determinare la dimensione $\Sigma$ delle curve d'ordine $n$ con molteplicità voluta pei punti d'incontro di $\varphi$ e $\psi$, visto che formano un sistema regolare; vedremo la dimensione $\Sigma'$ delle curve $\mathcal{C}^n$ che si possono scrivere $A\varphi+B\psi$, e vedremo $\Sigma=\Sigma'$.
\[
\Sigma= \left[n\right]-1-\sum\frac{(\alpha+\beta)(\alpha+\beta-1)}{2}
\]
%dubbio: la seconda sigma resa come somma%
Ora vediamo la $\Sigma'$ delle $f$ rappresentabili:
\[
f\equiv A^{[\beta-1]}\varphi+B^{[\alpha-1]}\psi\equiv (A-\psi X) \varphi +(B-\varphi X)\psi
\]
dove $X$ è un polinomio arbitrario. $A-\psi X=0$ passa con molteplicità 
%dubbio: almeno sembra aggiunto a matita%
$\beta-1$, qualunque $X$, purch\'{e} di grado
%dubbio: a matita aggiunto ``grandezza?''%
conveniente. Allora:
\begin{align}
\Sigma'&=\left[n-p\right]-\sum\frac{\beta(\beta-1)}{2}+\left[n-q\right]-\sum\frac{\alpha(\alpha-1)}{2}-1-\left[n-p-q\right]\\
&= \frac{n(n+3)}{2}-pq-\sum\frac{\alpha^2+\beta^2-\alpha-\beta}{2} \qquad (pq=\sum \alpha\beta)\\
%dubbio: sono tutte sigma o ci sono anche somme? v.pag36 per la prima volta in cui compare sigma ma lì l'abbiamo interpretata come somme%
\end{align}
Quindi $\Sigma'=\Sigma$, e il Teorema è dimostrato. Il passaggio al caso più generale porterebbe a introdurre intersezioni che cadono in punti infinitamente vicini; per \emph{contatti} si ricade nella dimostrazione precedente.
\end{enumerate}

%----------------------------pag.43-----------------------------------------%
\begin{teo}[Teorema di Cayley]
Se una curva d'ordine $n$ passa pei punti base d'un fascio di curve d'ordine $p<n$, allora la curva può generarsi mediante due fasci proiettivi: l'uno il nominato d'ordine $p$, l'altro un opportuno d'ordine $n-p$.
\end{teo}
Sia il fascio
\begin{equation}
\label{eq:cay1}
\varphi_p-\lambda \psi_p=0
\end{equation}

I $p^2$ punti base \emph{siano semplici}. Allora una $F_{p+p'}$ passante per essi può scriversi:
\[
F_{p+p'}=A_p'\varphi_p+B_p'\psi_p=0.
\]
\noindent
Consideriamo il fascio
\begin{equation}
\label{eq:cay2}
B_{p'}+\lambda A_{p'}=0
\end{equation}
\noindent
Se nei fasci \eqref{eq:cay1} e \eqref{eq:cay2}  chiamo corrispondenti due curve che provengono dallo stesso valore di $\lambda$, i due fasci vengono riferiti biunivocamente. La corrispondenza tra i due fasci si dice \emph{Proiettiva}. Effettivamente se si trattasse di due fasci di rette si avrebbe la ordinaria corrispondenza proiettiva.\\

\begin{figure}[!htbp]
\centering
\includegraphics[scale=0.5]{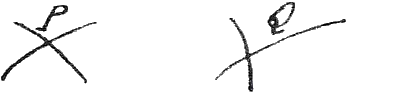}
%\caption{}
\end{figure}
\noindent
Anche se non si trattasse di fasci di rette associando alle curve le tangenti in uno dei punti base, mentre una curva varia descrivendo il primo fascio e l'altra il secondo, le tangenti a $P$ e $Q$ variano descrivendo due fasci proiettivi di rette.\\
Le $pp'$
%dubbio: nel manoscritto p_1 ma in realtà si riferisce a p' in eq:cay2% 
intersezioni di due curve corrispondenti \eqref{eq:cay1} e \eqref{eq:cay2}, al variare di $\lambda$ descrivono un luogo che si ottiene eliminando $\lambda$ da \eqref{eq:cay1} e \eqref{eq:cay2}; cioè 
\begin{equation}
\label{eq:cay3}
A\varphi+B\psi=0
\end{equation}
Cioè $F$ può riguardarsi come luogo dei punti d'incontro delle curve\dots\\

\noindent
Se si prende $p=1$ e $p=2$, è facile l'applicazione del Teorema. Per $p>2$ non è facile determinare il gruppo base (Es., per $p=3$, se prendiamo $9$ punti sulla curva questi non sono base; se ne prendiamo $8$ resta da determinare un nuovo punto
%----------------------------pag.44-----------------------------------------%
base che in genere non starà sulla curva; anche il $9^o$ deve stare sulla curva). Tuttavia si sa da Chasles  che si possono sempre trovare $p^2$ punti base d'un fascio di curve d'ordine $p<n$ su $F_n$. Quindi risulta che dato $n$, esso si può spezzare comunque  $n=p+p'$ ed $F$ si può sempre generare mediante due fasci proiettivi d'ordine $p$ e $p'$.\\
Così ogni conica può generarsi mediante due fasci proiettivi di rette (Steiner); una cubica mediante due fasci proiettivi, uno di coniche ed uno di rette. Così si studiano le proprietà senza ricorrere alle coordinate, cioè senza intervento analitico.\\
Steiner aveva in programma di generalizzare le curve d'ordine superiore per via simile, indipendentemente dalle coordinate; ma non lo ha svolto. Fu ripreso più tardi, ma non conviene.

\chapter{Corrispondenze algebriche}

%-------------------------------pag.45-------------------------------------%
\section{Corrispondenze tra due rette - In particolare biunivoche - Principio di corrispondenza su una retta - Computo dei punti doppi di una involuzione semplicemente infinita su una retta.}
\label{sec:9}

\subsection{Corrispondenze tra due rette.}

\begin{figure}[!htbp]
\centering
\includegraphics[scale=0.5]{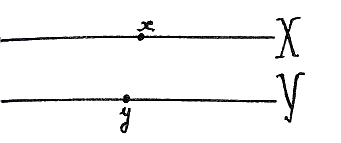}
%\caption{}
\end{figure}

Si può esprimere una corrispondenza $(m,n)$
\[
ax^my^m+bx^my^{n-1}+cx^{m-1}y^n+\dots+k=0
\]
Per un particolare valore di $x$, potrebbe annullarsi il coefficiente di $y^n$, e allora sembrerebbe che ad \underline{$x$} corrispondessero meno di \underline{$n$} punti in $Y$. Ma in realtà i punti mancanti sono andati all'infinito, come si giustifica subito introducendo le coordinate omogenee $x=\frac{x_1}{x_2}, y=\frac{y_1}{y_2}$.\\
Corrispondenza $(m,1)$: $\varphi(x)+y\psi(x)=0$; se tanto $\varphi$ come $\psi$ sono di grado $m$, si ha un'involuzione semplicemente infinita d'ordine $m$ (dipende da un parametro, $y$), e dato $P$ su $X$, ne sono determinate $m-1$.
%dubbio: incerta la lettura del manoscritto%

\subsection{Corrispondenze biunivoche.}
Quando anche $m=1$, si ha la corrispondenza $(1,1)$ o biunivoca:
\[
axy+bx+cy+d=0
\]
Questa è l'equazione della proiettività, sicch\'{e}: una corrispondenza algebrica biunivoca tra due rette è una corrispondenza proiettiva ($1^o$ Principio di Chasles).
Vediamo che in luogo di \emph{algebrica} basta \emph{analitica}. 
Se $x$ è funzione analitica ad un solo valore per tutti i valori della
%-------------------------------pag.46-------------------------------------%
$y$ complessa, e lo stesso per $y$, allora un legame tra $x$ ed $y$ è bilineare, si ha una corrispondenza proiettiva. Infatti: per ipotesi $y$ è funzione della variabile complessa $x$ salvo tutt'al più un numero finito di punti singolari, che possono essere poli o singolarità essenziali. Se un punto è un polo, $y$ diviene $\infty$.
Vediamo se ci possono essere delle singolarità essenziali: si sa dall'analisi che se $y(x)$ ammette una singolarità isolata per $x=a$, allora nell'intorno di $x=a$ si possono trovare infiniti punti $x_1,x_2,\dots, x_n,\dots$ dove $y$ assume uno stesso valore $y'$. \\
Allora ad un valore $y'$ corrispondono infiniti valori di $x$, mentre noi abbiamo supposto anche la $x$ funzione ad un solo valore di $y$. Dunque non vi possono essere singolarità essenziali, cioè $y=y(x)$ è funzione regolare in tutto il piano, ad un sol valore salvo un numero finito di poli. Allora $y$ è funzione razionale di $x$, come $x$ di $y$ la corrispondenza è algebrica biunivoca, quindi bilineare.

\subsection{Geometria proiettiva di curve.}
Chasles ne ha tratto: 
 data una conica, tra i due fasci di centri $S$ ed $S'$ viene a stabilirsi una corrispondenza algebrica (curva = equazione algebrica) biunivoca, quindi proiettiva; e si 
%dubbio: incerta la lettura del manoscritto, probabilmente c'è una correzione%
ha la Generazione proiettiva delle coniche: luogo dei punti d'incidenza di due rette corrispondenti in due fasci proiettivi.\\

\begin{figure}[!htbp]
\centering
\includegraphics[scale=0.6]{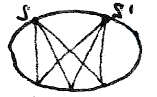}
%\caption{}
\end{figure}
\noindent
Così un punto variabile su una cubica sghemba è determinato
%-----------------------------------pag.47---------------------------%
da due piani variabili in due fasci proiettivi intorno a due corde della cubica. Una cubica sghemba si può riguardare come luogo dei punti d'incontro di $3$ piani corrispondenti in $3$ fasci proiettivi.\\
Un caso notevole di corrispondenza algebrica tra  punteggiate si ha su una stessa retta. Si trova subito per una corrispondenza (m,n) su una retta:
%dubbio: lo riporto come teorema?%
\begin{teo}[$II^o$ Principio di  Corrispondenza (Chasles)]
Una corrispondenza $(m,n)$ su una retta ha $(m+n)$ punti uniti (corrispondenti a sè stessi).
\end{teo}
\noindent
 Principio che sostituisce il Teorema fondamentale dell'Algebra in ragionamenti di carattere geometrico. Naturalmente, per avere con certezza questo risultato, occorre operare in coordinate omogenee.

\subsection{Punti doppi d'una involuzione $\infty^1$ sulla retta.}
\begin{figure}[!htbp]
\centering
\includegraphics[scale=0.6]{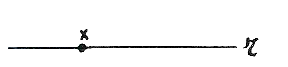}
%\caption{}
\end{figure}
Questa è data da $f(x)+\lambda\varphi(x)=0$; sia d'ordine $n$. Chiamiamo corrispondenti due punti, $x,y$, quando fanno parte di uno stesso gruppo della involuzione. Così immaginando sulla retta $r$ due punteggiate, si ha una corrispondenza $(n-1,n-1)$. Allora si vede: una involuzione semplicemente infinita, di grado $n$, su una retta, ha $2(n-1)$ punti doppi o uniti (due punti d'un gruppo coincidenti).

\emph{Algebricamente}: ricordiamo che se 
\begin{equation}
\label{eq:corrA}
F_n(x)=0
\end{equation}
ha una radice doppia, questa soddisfa la \eqref{eq:corrA} e la $F'_{n-1}=0$; se una radice $i$-upla, la sua derivata è soddisfatta dalla radice multipla secondo $2(i-1)$.
%dubbio: qui nel manoscritto sembra esserci a matita un no%
%-----------------------------------pag.48---------------------------%
In coordinate omogenee, $$x=\frac{x_1}{x_2},\quad F(x_1,x_2)=0,\quad \frac{\partial F}{\partial x_1}=0$$
 Tenendo conto del Teorema d'Eulero ($F=0$ è omogenea) si vede che una radice doppia soddisfa anche $\frac{\partial F}{\partial x_2}=0$.\\
Così se la radice è multipla secondo $i$ la derivata è soddisfatta dalla radice multipla secondo $2(i-1)$.
%dubbio: nel manoscritto questo capoverso è segnato da un no che sembra in matita.%
Sicch\'{e} in coordinate omogenee, condizione perch\'{e} una radice sia doppia è l'annullamento delle derivate sostituendo la radice.
%dubbio: incerta la lettura del manoscritto%

Si abbia dunque l'involuzione: 
\begin{equation}
\label{eq:corr1}
F=f(x_1,x_2)+\lambda \varphi(x_1,x_2)=0
\end{equation}
Per un $\lambda$, la \eqref{eq:corr1} ammette una radice doppia $(x_1,x_2)$. Allora
\begin{equation}
\label{eq:corr2}
\frac{\partial f}{\partial x_1}+\lambda \frac{\partial\varphi}{\partial x_1}=0 , \qquad 
\frac{\partial f}{\partial x_2}+\lambda \frac{\partial\varphi}{\partial x_2}=0  
\end{equation}
Queste \eqref{eq:corr2} coesistono quando $\lambda$ ha un valore da far sì che nel gruppo \eqref{eq:corr1} ci sia un punto doppio. Se le \eqref{eq:corr2} coesistono dev'essere uguale a $0$ il Jacobiano:
\begin{equation}
\label{eq:corr3}
\frac{\partial f}{\partial x_1}\frac{\partial\varphi}{\partial x_2}- \frac{\partial f}{\partial x_2}\frac{\partial\varphi}{\partial x_1}=0
\end{equation}
Questa, di grado $2(n-1)$ dev'essere soddisfatta dalle coordinate d'ogni punto doppio della involuzione, che così sono in numero $2(n-1)$.
%dubbio: qui parte un pezzo segnato con un no che sembra a matita nel manoscritto%
Ecco trovato non solo il numero dei punti doppi, ma le equazioni che risolte darebbero le loro coordinate.\\
Uno dei gruppi di $F$ abbia un punto $i$-plo, in cui i punti coincidono. Allora $\frac{x_1}{x_2}$ è radice $i$-pla di $F=0$, e quindi $(i-1)$-upla delle $\frac{\partial f}{\partial x_1}=0 \quad  \frac{\partial f}{\partial x_2}=0$. Vuol dire: quel punto $i$-uplo conta come $(i-1)$ tra i punti doppi dell'involuzione.
%-----------------------------------pag.49---------------------------%

\section{Corrispondenze algebrica tra i piani - In particolare razionali e birazionali - Trasformazione quadratica - Cenno all'applicazione delle Trasformazioni quadratiche alla risoluzione delle singolarità sulle curve.}
\label{sec:10}

\subsection{Corrispondenze tra due piani $\pi(x,y)$ e $\pi'(x',y')$.}
Si può dare con due, o tre equazioni. 
\begin{equation}
\label{eq:piani1}
f(x,y,x',y')=0 \qquad \varphi=0 \qquad \psi=0
\end{equation}
Se tra le \eqref{eq:piani1} passano opportune relazioni, può darsi che dati valori arbitrari ad $(x,y)$, le \eqref{eq:piani1} ammettano qualche coppia $(x',y')$ comune. \emph{Se la corrispondenza è univoca}, per trovare $P'$ corrispondente a $P$: nelle \eqref{eq:piani1}, elimino $y'$ tra la prima e la seconda, poi tra la prima e la terza; ottengo, con operazioni razionali, due equazioni in $x'$, con una radice comune. Ora si sa dall'Algebra: se due equazioni algebriche in $x'$ hanno una sola soluzione comune, questa si ottiene razionalmente (operazione di M.C.D.):
\begin{equation}
\label{eq:piani2}
x'=\frac{L(x,y)}{N(x,y)} \qquad y'=\frac{M(x,y)}{N(x,y)}
\end{equation}
corrispondenza univoca o razionale in un senso ($L,M,N$ polinomi primi fra loro).
In coordinate omogenee le \eqref{eq:piani2} divengono 
\begin{gather*}
NX'-LZ'=0, \qquad NY'-MZ'=0 \\
 \text{Matrice }
\begin{pmatrix}
N & 0 & -L\\
0 & N & -M
\end{pmatrix}
\end{gather*}
%dubbio: incerta la lettura di Onde%
Onde, risolute rispetto alle incognite $(X',Y'Z')$ danno:
\begin{equation}
\label{eq:piani2b}
\rho X'=L ,\quad \rho Y'=M, \quad \rho Z'=N
\end{equation}
Così a $P$ corrisponde un determinato $P'$, eccetto i $P$ per cui $L=M=N=0$,
%-----------------------------------pag.50---------------------------%
punti base comuni al sistema (lineare $\infty^2$) di curve:
\begin{equation}
\label{eq:piani3}
\lambda L+\mu M+\nu N=0,
\end{equation}
punti \emph{fondamentali} della corrispondenza, che sono in numero finito (poich\'{e} abbiamo escluso il caso che $L,M,N$ abbiano una componente comune) e possono mancare.\\
Viceversa, dato $P'(x',y')$, per avere i punti corrispondenti $P_i$, si ricorre alle \eqref{eq:piani2}:
\begin{equation*}
L-x'N=0, \quad M-y'N=0
\end{equation*}
che sono due curve del sistema \eqref{eq:piani3}; queste, oltre i punti base, hanno generalmente qualche altro punto comune; questi saranno i $P_i$ corrispondenti a $P'.$
%dubbio: la parte successiva è posta tra paretesi quadre nel manoscritto, forse da mettere come osservazione?%
[Eccezione succederebbe se dato $P'$, le \eqref{eq:piani2} avessero comuni i soli punti fondamentali; se si costringessero a passare per un altro $P$, avendo un punto comune oltre i base, avrebbero una parte comune: o le due curve coinciderebbero, ed $L,M,N$ non sarebbero allora linearmente indipendenti (questo caso va escluso, ch\'{e} non ha interesse); o le le due curve si spezzano ed hanno una parte comune, come allora tutte le curve del sistema \eqref{eq:piani3}; allora si ricava che ogni curva della rete
 \eqref{eq:piani3} si spezza in più curve che al variare dei parametri descrivono uno stesso fascio (degenerazione della corrispondenza, Jacobiano nullo [L,M,N]).]\\
%dubbio: qui finisce la parte tra parentesi quadre%
Eccezioni non si presentano supponendo la rete \eqref{eq:piani3} di curve irriducibili.
La corrispondenza algebrica tra $\pi$ e $\pi'$ si può indicare $(m,1)$.\\
Se un $P'$ descrive in $\pi'$ la retta\footnote{Le coordinate omogenee chiamiamole $x,y,z$; non c'è da dar luogo a confusione.} 
\begin{equation}
\label{eq:piani4}
\lambda x'+\mu y' + \nu z'=0
\end{equation}
%-----------------------------------pag.51---------------------------%
il $P$ corrispondente descrive la curva 	\eqref{eq:piani3} di $\pi$. Se la \eqref{eq:piani4} descrive il $\pi'$ rigato, la \eqref{eq:piani3} descrive la rete \eqref{eq:piani3} di $\pi$. Se la \eqref{eq:piani4} descrive un fascio di centro $P'$, \eqref{eq:piani3} descriverà un fascio di curve avente i punti base della rete \eqref{eq:piani3} più altri \underline{$m$} punti base (i corrispondenti a $P'$).
Dunque il modo più generale di porre una corrispondenza $(m,1)$ tra $\pi$ e $\pi'$ è: si prende in $\pi$ una rete qualunque irriducibile $\eqref{eq:piani3}$, e alle curve faccio corrispondere le rette \eqref{eq:piani4} di $\pi'$.
Dunque le formole che ci danno una corrispondenza fra due piani sono le \eqref{eq:piani5b}
\begin{subequations}
\label{eq:piani5}
\begin{align}
&x'=\frac{L}{N}  && y'=\frac{M}{N} \\
&\label{eq:piani5b}\lambda L+\mu M+\nu N=0 \qquad&&\lambda x'+\mu y'+\nu z'=0
\end{align}
\end{subequations}
Esempio: rete di coniche in $\pi$: $\lambda x^2+\mu y^2 + \nu =0$,
rete di coniche in $\pi'$: $\lambda x'+\mu y' + \nu =0$.\\
A $P$ corrisponde un $P'$, a $P'$ corrispondono quattro $P$; se tre sono i punti fondamentali di $\pi$ la corrispondenza è biunivoca.
\begin{teo} 
Ai punti dell'intorno infinitesimo di un punto fondamentale $O$ di $\pi$, corrispondono in $\pi'$ i punti di una certa curva geometrica $\bar{O'}$. Si suol dire che ad $O$ corrisponde la curva $\bar{O'}$.
\end{teo}
\begin{figure}[!htbp]
\centering
\includegraphics[scale=0.5]{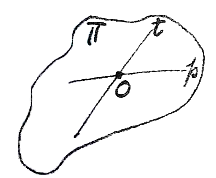}
%\caption{}
\end{figure}
%\begin{proof}
\noindent
Infatti: ad $O$ non corrisponde alcun punto in $\pi'$. Se obblighiamo le curve della rete a passare per un punto infinitamente vicino ad $O$ \emph{nella direzione $p$}, si ha \emph{un fascio} nella rete, con $p$ tangente; a detto punto centro del fascio, corrisponde in $\pi$
%-----------------------------------pag.52---------------------------%
un punto. Analogamente per le altre direzioni, e il Teorema è dimostrato.
%\end{proof}
L'ordine di $\bar{O'}$ corrisponde al numero delle tangenti variabili che le curve della retta hanno sul punto base $O$.\\
\begin{figure}[!htbp]
\centering
\includegraphics[scale=0.5]{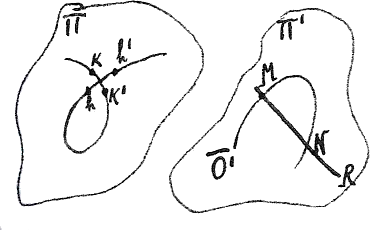}
%\caption{}
\end{figure}
\noindent
Infatti: a $k$ e $k'$ corrisponde un solo punto su $\bar{O'}$ e su $R$ (immagine della curva disegnata sopra), 
%dubbio: qui c'è un riferimento che potrebbe cambiare a seconda dell'impaginazione, sarebbe il caso inserire il riferimento alla numerazione della figura%
e sia $M$; così ad $h$ ed $h'$, $N$. Variando le tangenti variano questi due punti, quindi si ha una conica nel caso disegnato.\\
Se le curve della rete hanno in $O$ tangenti tutte fisse, ad $O$ corrisponde un punto $O'$ in $\pi'$.

\subsection{Corrispondenze biunivoche o birazionali.}
Come la corrispondenza univoca è razionale, così la biunivoca si vede identicamente che è birazionale (Cremoniana). Dunque, anche dato $P(x',y')$ si ha:
\begin{equation}
\label{eq:piani6}
x=\frac{L'}{N'}, \quad y=\frac{M'}{N'}
\end{equation}
La costruzione di una corrispondenza Cremoniana si riduce alla costruzione su uno dei due piani, es. $\pi$, di una \emph{rete} che abbia queste proprietà: due curve si seghino in \emph{un sol punto} (variabile con i parametri) fuori dei punti base (\emph{Rete Omaloidica}). Alle rette di $\pi'$ corrispondono le curve omaloidiche di $\pi$; e viceversa: alle rette di $\pi$ corrispondono le curve $\lambda L'+\mu M'+\nu N'=0$ di $\pi'$, rete omaloidica.
Gli ordini delle curve dei due piani sono uguali, come si vede
%-----------------------------------pag.53---------------------------%
con facili considerazioni.\\

\begin{figure}[!htbp]
\centering
\includegraphics[scale=0.5]{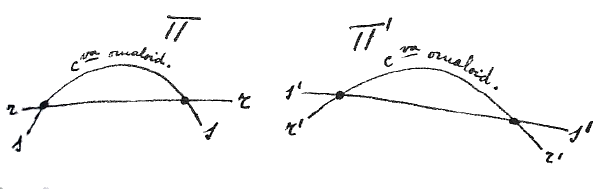}
%\caption{}
\end{figure}
\noindent
I punti base della rete possono avere le molteplicità\footnote{Cremona stabilì essere possibile costruire corrispondenze di tutti gli ordini.} 
%dubbio: Continua.%
$\alpha_1, \alpha_2, \alpha_3, \dots$; la somma dei quadrati sarà il numero delle intersezioni assorbite dai punti base:
\begin{equation}
\label{eq:piani7}
\left. 
\begin{aligned} 
 \alpha_{1}^{2}+\alpha_{2}^{2}+\alpha_{3}^{2}+\dots=n^2-1\\
\text{Ancora: }\quad \frac{\alpha_1(\alpha_1+1)}{2}+\frac{\alpha_2(\alpha_2+1)}{2}+\dots=\frac{n(n+3)}{2}-2
\end{aligned} 
\right\} 
\end{equation}
poich\'{e} vi sono due parametri disponibili. Così ho due legami tra $\alpha_1,\alpha_2,\alpha_3$, ed $n$ noto; non resta che esaminare i sistemi di soluzioni intere e vedere che sono costruibili le reti.
Per $n=1$ si ha la collineazione od omografia tra i due piani. Per $n=3\, $ $\alpha_1=2$ punto base doppio, e quattro punti base semplici. Il caso che ci interessa è $n=2$, cioè a Rette
%dubbio: maiuscola?%
 di un piano corrispondono coniche aventi una sola intersezione variabile, cioè circoscritte ad un triangolo.

\subsection{Trasformazioni quadratiche.}
Alle coniche circoscritte ad un triangolo di $\pi$ corrispondono le rette di $\pi'$, e viceversa.\\
\begin{figure}[!htbp]
\centering
\includegraphics[scale=0.5]{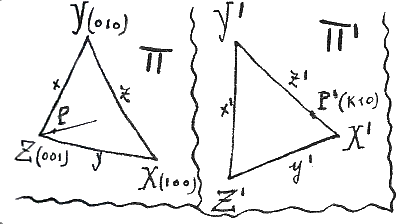}
%\caption{}
\end{figure}
\noindent
I due triangoli assumiamoli come triangoli fondamentali dei due piani. Ad una conica generica pei tre punti di $\pi$
\[
\begin{cases}
\lambda yz+\mu xz+\nu yz=0\\
\alpha x'+\beta y'+\gamma z'=0
\end{cases}
\]
facciamo corrispondere in $\pi'$ la retta
\begin{equation}
\label{eq:piani8}
\begin{cases}
x': y': z'=yz: xz: xy\\
x : y : z = y'z' : x'z' : x'y'
\end{cases}
\end{equation}
%dubbio: incerta l'organizzazione della pagina%
Si ricordino le \eqref{eq:piani2b}: dividendo i tre termini del secondo membro per $xyz$, del primo per $x'y'z'$, ed invertendo, permettono di passare da ogni punto di $\pi'$ a punti di $\pi$.\\
%-----------------------------------pag.54---------------------------%
\noindent
Alla retta di $\pi$
\[
\alpha x+\beta y+\gamma z=0
\]
corrisponde in $\pi'$
\[
 \lambda y' z'+\mu x' z'+ \nu x'y'=0
\]
circoscritta al triangolo fondamentale.
%\[
%\left. 
%\begin{aligned} 
% \text{Alla retta di $\pi$} \quad &\alpha x+\beta y+\gamma z=0\\
%\text{corrisponde in $\pi'$} \quad & \lambda y' z'+\mu x' z'+ \nu x'y'=0
%\end{aligned} 
%\right\} 
%\text{circoscritta al triangolo fondamentale.}
%\]
%dubbio: quale delle due scelte è migliore?%
\begin{teo} 
Ai punti infinitamente vicini a $Z$ in tutte le direzioni uscenti da $Z$ corrispondono i punti del lato $X'Y'$ del triangolo fondamentale del secondo piano.
\end{teo}
%\begin{proof}
\noindent
Dalla prima \eqref{eq:piani8}:
\begin{equation}
\label{eq:piani8b}
\rho x'=yz, \quad \rho y'=xz, \quad \rho z'=xy
\end{equation}
Il punto fondamentale $Z$ non ha punti corrispondenti in $\pi'$. Consideriamo un punto $P$ infinitamente vicino a $Z$ in una data direzione di rapporto direttivo $\frac{y}{x}=k$, cioè di coordinate $P(x\to0,~ y\to0,~1)$. Per trovare $P'$, poniamo nelle \eqref{eq:piani8b} $z=1$, e dividiamole per $x$; divengono:
\[
\sigma x' =\frac{y}{x}, \quad \sigma y'=1, \quad \sigma z'=y;
\]
considerando $P$:
\[
\sigma x'=k, \quad \sigma y'=1, \quad \sigma z'=0;
\]
cioè alla posizione limite di $P$ infinitamente vicino a $Z$ in direzione di rapporto direzionale $k$, corrisponde in $\pi'$ il punto $P'(k,1,0)$, appartenente al lato $X'Y'$. Così l'intorno di $Z$ si dilata sopra la retta $X'Y'$ di $\pi'$. Fra le direzioni delle rette uscenti da $Z$ e i punti di $X'Y'$ vi è una corrispondenza proiettiva. Questo si ripete per l'intorno di ogni punto fondamentale.
%\end{proof}
\begin{teo} 
Alle rette uscenti da $Z$ di $\pi$ corrispondono rette uscenti da $Z'$ in $\pi'$, trascurando la seconda componente delle coniche, $X'Y'$, la quale proviene soltanto dal punto $Z$. Vi è corrispondenza proiettiva fra il fascio di
%-----------------------------------pag.55---------------------------%
rette uscenti da $Z$ e le corrispondenti da $Z'$.
\end{teo}
%\begin{proof}
\noindent
Scrivendo una retta uscente da $Z$, ricordando la seconda \eqref{eq:piani8}, si ha $\nu=0$, quindi l'equazione della conica trasformata si riduce: $z'(\lambda y'+\mu x')=0$, che si spezza precisamente in $z'=0$, $\lambda y'+\mu x'=0$; questa è una retta da $Z'$ i cui punti corrispondono ai punti della retta per $Z$.
%\end{proof}
\begin{teo} 
Una curva $\mathcal{C}_n$ di $\pi$ che non passi per $X,Y,Z$, si muta in una $\mathcal{C'}_{2n}$ di $\pi'$ che passa con molteplicità $n$ per $X',Y',Z'$.
\end{teo}
%\begin{proof}
\noindent
$\mathcal{C}_n$ si muterà in una $\mathcal{C'}_{n'}$ che passerà con una certa molteplicità pei punti fondamentali del secondo piano\footnote{Perchè i punti di contatto di $\mathcal{C}$ con un lato del primo triangolo, corrispondono a punti infinitamente vicini al vertice del secondo triangolo.}. Ad una retta di $\pi'$ corrisponderà una conica di $\pi$ per $X,Y,Z$ (non appartenenti a $\mathcal{C}_n$), che sega la  $\mathcal{C}_n$ in $2n$ punti; quindi $n'=2n$. Un punto di  $\mathcal{C}'_{n'}$ infinitamente vicino ad $X'$ proviene da un punto di  $\mathcal{C}_n$ su $YZ$; ma questi sono $n$, quindi sono $n$ i punti infinitamente vicini ad $X'$ su  $\mathcal{C'}_{n'}$, ossia \underline{$n$} è la molteplicità di  $\mathcal{C'}_{2n}$ in $X'$; come idem in $Y'$ e $Z'$.
%\end{proof}
\begin{teo} [Comprende il precedente]
Una curva  $\mathcal{C}_n$ di $\pi$ passante con molteplicità $\alpha, \beta, \gamma$ per $X,Y,Z$, si trasforma in una  $\mathcal{C}_{2n-\alpha-\beta-\gamma}$ di $\pi'$ che passa con molteplicità $\alpha'=n-\beta-\gamma, \, \beta'=n-\alpha-\gamma, \, \gamma'=n-\alpha-\beta$ per $X', Y', Z'$.
\end{teo}
%\begin{proof}
\noindent
$n'$ si trova come al Teorema precedente, solo bisogna escludere i passaggi della curva $\mathcal{C}$ che cadono nei punti base; resta $n'=2n-\alpha-\beta-\gamma$.\\
Anche le molteplicità si trovano come sopra; $\mathcal{C}_n$ ha con $YZ$
%-----------------------------------pag.56---------------------------%
comuni $n$ punti, di cui però $\beta+\gamma$ cadono in $Y$ e $Z$; quindi
\[
\alpha'=n-\beta-\gamma, \quad \text{ come } \beta'=n-\alpha-\gamma, \quad \gamma'=n-\alpha-\beta.
\]
%\end{proof}
\noindent
Un caso noto che rientra in questa Teoria è fornito nella Trasformazione per raggi vettori reciproci (o inversione) che manda le rette in coniche circoscritte al triangolo formato dal centro d'inversione (Origine delle coordinate) e dai due punti ciclici. I due triangoli fondamentali nei due piani che sono sovrapposti coincidono sopra detto triangolo. Se un cerchio si trasforma in un cerchio non è che un caso particolare che una conica passante per $XY$ (punti ciclici) si muta in una conica per $X'Y'$.

\subsection{Applicazione alla risoluzione delle singolarità di curve.}

Si abbia una curva algebrica con un punto doppio; \\

\begin{figure}[!htbp]
\centering
\includegraphics[scale=0.5]{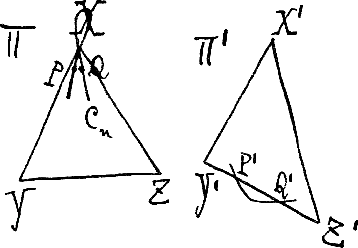}
%\caption{}
\end{figure}
\noindent
applichiamo ad essa una Trasformazione quadratica di cui un punto fondamentale cada in $X$. Nella $\mathcal{C'}$, i due punti infinitamente vicini ad $X$, $P$ e $Q$, in $\mathcal{C}$, si trasformeranno in $P'$ e $Q'$. Cioè il punto doppio ordinario viene risolto in due punti $P'$ e $Q'$; però la trasformazione ha fatto nascere nuovi punti multipli: $X'$ $n$-uplo, $Y'$ e $Z'$ $(n-2)-upli$. Ora $P$ sia \emph{singolare} (a tangenti coincidenti). Allora $P'$ e $Q'$ vengono a coincidere nella trasformata. Questa coincidenza può essere
%-----------------------------------pag.57---------------------------%
per contatto con $Y'Z'$ se $X$ è una cuspide o per punto doppio se $X$ è un contatto almeno quadripunto, come nel Teorema.
\begin{figure}[!htbp]
\centering
\includegraphics[scale=0.5]{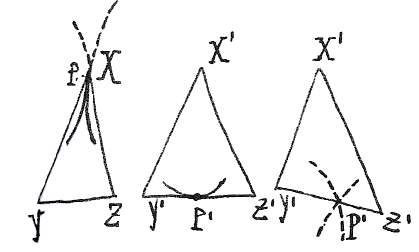}
%\caption{}
\end{figure}
(Nel caso della Cuspide, non si dice essere $X$ origine di due rami, ma origine di un ramo del secondo ordine). Ma anche $P'$ potrebbe essere singolare. Allora occorre applicare una nuova trasformazione quadratica di cui un punto fondamentale cada in $P'$ e gli altri due ad arbitrio.
\begin{teo} 
Dopo un numero finito di trasformazioni quadratiche (Prodotto) si deve arrivare ad una curva trasformata che non ha più punti multipli singolari.\\
Ogni trasformazione fa nascere nuovi punti multipli, ma ordinari.
\end{teo}
\noindent
Se si studiano le proprietà delle curve che non si alterano per trasformazioni birazionali, è lecito sostituire ad una curva con punti multipli singolari una trasformata nel modo anzidetto, con punti multipli ordinari.

\section{Corrispondenze algebriche tra curve, in particolare razionali e birazionali.}
\label{sec:11}
Siano date le curve
\begin{equation}
\label{eq:alg1}
f(x,y)=0 \quad \text{su } \pi \, \text{ ed } \, \varphi(x',y')=0 \quad \text{su } \pi'
\end{equation}
%-----------------------------------pag.58---------------------------%
e due equazioni
\begin{equation}
\label{eq:alg2}
F(x,y,x',y')=0 \qquad \Phi(x,y, x', y')=0
\end{equation}
Dati $x$ ed $y$ le \eqref{eq:alg2} possono essere tali che per ogni coppia $(x,y)$ soddisfacente alla $f$, le \eqref{eq:alg2} e la $\varphi$ abbiano un certo numero di soluzioni comuni, $\nu$; in modo analogo ad ogni $(x',y')$ di $\varphi$ corrispondano $\mu$ punti di $f$. Allora le \eqref{eq:alg2} definiscono una corrispondenza $(\mu, \nu)$ tra le \eqref{eq:alg1}. (Una sola \eqref{eq:alg2} potrebbe non bastare, mentre si dimostra che due bastano).

\subsection{Corrispondenze razionali o univoche.}
Ad $(x,y)$ corrisponde un solo $(x',y')$, non viceversa; allora $x'$ ed $y'$ si possono ricavare razionalmente da $x,y$. Così tra le \eqref{eq:alg1}, una corrispondenza razionale $(\mu,1)$, si può rappresentare
\begin{equation}
\label{eq:alg2b}
x'=\frac{L(x,y)}{N(x,y)}, \quad y'=\frac{M(x,y)}{N(x,y)}
\end{equation}
Queste \eqref{eq:alg2b} rappresentano una corrispondenza $(m,1)$ tra i due piani in considerazione, ed è evidentemente $\mu\leq m$.\\
Così la corrispondenza razionale di curve è contenuta nella corrispondenza razionale tra i rispettivi piani; non è necessariamente $\mu=m$.
\subsection{Corrispondenze birazionali o biunivoche.}
Dati $x',y'$ soddisfacenti alla $\varphi$, le \eqref{eq:alg2b} e la $f$ abbiano una sola soluzione $(x,y)$ comune; allora questa è ottenibile con operazioni razionali:
\begin{equation}
\label{eq:alg2c}
x=\frac{L'(x',y')}{N'(x',y')}, \quad y=\frac{M'(x',y')}{N'(x',y')}
\end{equation}
%-----------------------------------pag.59---------------------------%
In questo caso la corrispondenza tra le due curve è birazionale; essa è ancora rappresentata dalle \eqref{eq:alg2b}. Si noti di non dire che la corrispondenza birazionale tra due curve è contenuta in una corrispondenza birazionale tra i rispettivi piani; in generale non è così. Perchè le \eqref{eq:alg2b} rappresentano in generale una corrispondenza razionale, e non birazionale tra i due piani. Se dalle \eqref{eq:alg2b} abbiamo ottenute le \eqref{eq:alg2c} è perchè abbiamo tenuto conto anche dell'equazione $f=0$. Cioè le corrispondenze \eqref{eq:alg2b}, \eqref{eq:alg2c} tra $\pi$ e $\pi'$ in genere non sono le inverse l'una dell'altra; chè risolvendo le \eqref{eq:alg2b} rispetto ad $x$ ed $y$, non si ottengono le \eqref{eq:alg2c} se non si tiene conto di $f=0$; ossia le \eqref{eq:alg2b} rappresentano una corrispondenza $(m,1)$, ed in genere $m\neq 1$.
In altri termini le due corrispondenze \eqref{eq:alg2b} e \eqref{eq:alg2c} non sono le inverse, tranne che sulle curve \eqref{eq:alg1}. 
Quindi: data una corrispondenza razionale tra due piani se ne deduce quella tra due curve; viceversa no, a meno che non si trattasse di un caso particolarissimo.\\
Le \eqref{eq:alg2b} mutano la retta $\lambda x'+\mu y'+\nu=0$ di $\pi'$ nella curva di $\pi\,$ $\lambda L+ \mu M+\nu N=0$.\\
Come si vede si hanno due reti corrispondenti.
%-----------------------------------pag.60---------------------------%
Ad un fascio della rete di curve (per $P$) corrisponde un fascio di rette su $\pi'$, che avrà come centro un punto corrispondente $P'$.
Se $P$ è fondamentale non avrà punto corrispondente su $\pi'$; ma possiamo chiamare così il centro del fascio trasformato. Così anche se un punto di $f$ è fondamentale, ad esso si può far corrispondere un punto determinato di $\varphi$. \\
Se ho un punto doppio ordinario su $f$ (punto fondamentale della rete), considerando due punti infinitamente vicini sempre su $f$, si hanno due punti su $\varphi$. Ciò potrebbe pensarsi come una eccezione alla biunivocità. Ma in realtà il punto doppio può pensarsi come la sovrapposizione di due punti appartenenti ai due rami della curva. Se il punto doppio non fosse punto base, ad esso corrisponderebbe un altro punto doppio.\\

\noindent
Ora passiamo allo studio delle proprietà Invarianti per trasformazioni birazionali (\emph{Geometria sopra le curve}). 
\`{E} facile vedere che il \emph{prodotto} di due trasformazioni birazionali è una trasformazione birazionale.

%---------------------------------61------------------------------------------%
\chapter{Geometria sopra le curve}
\section{Serie lineari di gruppi di punti su una curva - Serie lineari semplici e composte.}
\label{sec:12}

\subsection{Serie lineari su una retta (Serie razionale).}
Sia la retta $y=ax+b$ ed un sistema lineare di $\infty^r$ curve d'ordine $n$ (tutte le $\varphi$ di grado $n$):
\[
\lambda_0\varphi_0(x,y)+\lambda_1\varphi_1(x,y)+\dots+\lambda_r\varphi_r(x,y)=0
\]
Al variare dei parametri essenziali, queste segano sulla retta una serie lineare di gruppi di punti, d'ordine $n$ e dimensione $r$, $g^{r}_{n}$.\\
Analiticamente: sia $y$ funzione lineare della variabile complessa $x$: $y=ax+b$, e sia data una funzione razionale (tutte le $\varphi$ di grado $n$):
\[
\lambda_0=\frac{\lambda_1\varphi_1(x,y)+\lambda_2\varphi_2(x,y)+\dots+\lambda_r\varphi_r(x,y)}{\varphi_0(x,y)}
\]
Ad ogni valore di $\lambda_0$ corrispondono (fissati $\lambda_1, \lambda_2, \dots, \lambda_r$ ) sulla retta $y=ax+b$, un \emph{gruppo} di $n$ punti, \emph{di livello} $\lambda_0$; facendo variare $\lambda_1, \lambda_2, \dots, \lambda_r$, si ha naturalmente una serie di gruppo di livello $\lambda_0$.

\subsection{Serie lineari su una curva.}
\begin{equation}
\label{eq:ser1}
f(x,y)=0
\end{equation}
Si potrebbe, come sopra, introdurre una funzione razionale della $f$ (Riemann) $$\lambda=\frac{P(x,y)}{Q(x,y)}$$ ($x,y$ legati dalla \eqref{eq:ser1}). Dato $\lambda$, si ha sulla $f$ il gruppo di livello $\lambda$. Vi sarà il gruppo degli zeri
%---------------------------------62------------------------------------------%
($P=0$), il gruppo dei poli ($Q=0$) o degli infiniti. Variando $\lambda$, si ha una serie lineare di dimensione $1$.
Estendiamo queste considerazioni (dimensione maggiore di $1$) per via Geometrica.\\
Si consideri un sistema lineare di $\infty^r$ curve di ugual ordine:
\begin{equation}
\label{eq:ser2}
\lambda_0\varphi_0(x,y)+\lambda_1\varphi_1(x,y)+\dots+\lambda_r\varphi_r(x,y)=0,
\end{equation}
indipendenti. Le \eqref{eq:ser2} segano sulla 	\eqref{eq:ser1} una serie lineare di gruppi di punti, $g_n$. Nessuna \eqref{eq:ser2} contenga come parte la \eqref{eq:ser1}; nelle \eqref{eq:ser2} vi sono solo $r$ parametri arbitrari. Quindi dati $r$ punti generici su $f$, è determinata la \eqref{eq:ser2} che vi passa (per posizioni particolari, potrebbero passare infinite \eqref{eq:ser2}; ma ciò non può accadere sempre: si verrebbe altrimenti alla conclusione che per un punto generico passano infinite \eqref{eq:ser2}, che invece costituiscono un sistema lineare). In questo caso la $g_n$ ha dimensione $r$ $$g^{r}_{n}$$
Se una \eqref{eq:ser2}, per esempio la $\varphi_r$ contiene la $f$, allora la serie \eqref{eq:ser2} $\to$ \eqref{eq:ser1} si vede che è anche segata sulla \eqref{eq:ser1}  da $\lambda_0\varphi_0(x,y)+\lambda_1\varphi_1(x,y)+\dots+\lambda_{r-1}(x,y)\varphi_{r-1}=0$.
Così, se \underline{$s$} curve linearmente indipendenti contengono la \eqref{eq:ser1}, le \eqref{eq:ser2} segano sulla \eqref{eq:ser1} una $g^{r-s}_{n}$.\\
Se vi sono intersezioni fisse \eqref{eq:ser2}-\eqref{eq:ser1} per definizione se ne può fare astrazione; e viceversa aggiungendo ai gruppi della $g^{r}_{n} \dots$.\\
Esempio: serie segata sulla retta dalle coniche di un fascio; nessuna contiene la retta, si ha $g^{1}_{2}$; si ha il Teorema
%---------------------------------63------------------------------------------%
Desargues sulle coniche. Se dei quattro punti base uno stesse sulla retta, se ne potrebbe fare astrazione e si avrebbe $g_{1}^1$; se ulteriormente la retta passa per un altro punto base, una conica contiene la retta, e si avrebbe un $g^{0}_1$, cioè si ha un solo gruppo di due punti fissi (base).

\`{E} sempre $r<n$.\\ 
Le serie lineari
%dubbio: c'è una virgola tra il soggetto e il verbo che andrebbe tolta% 
costituiscono un concetto della geometria sulle curve; cioè si conservano (trasformazioni lineari) ordine e dimensione.
\subsection{Trasformazioni razionali; serie composte.}
Siano due curve:
\begin{equation}
\label{eq:sis3}
f(x,y)=0 \qquad f'(x',y')=0
\end{equation}
legate da una trasformazione \emph{razionale}:
\begin{equation}
\label{eq:ser4}
x'=\frac{L(x,y)}{N(x,y)} \qquad y'=\frac{M(x,y)}{N(x,y)}
\end{equation}
Ad un punto di $f'$ corrispondono $\nu$ punti di $f$. Su $f'$ si consideri una serie lineare segata da 
\[
\lambda_0\varphi_0(x,y)+\dots+\lambda_r\varphi_r(x,y)=0
\]
Applicando la \eqref{eq:ser4}, $f'$ si trasforma in $f$, e la serie lineare in una nuova serie lineare (per questo basta dunque la razionalità) segata da:
\[
\lambda_0\Phi_0(x,y)+\dots+\lambda_r\Phi_r(x,y)=0
\]
Ad un $P'$ corrispondono $\nu$ $P$; così a $P \to P'$; a $P'\to P,  \dots P_{\nu-1}$. Così tutti i gruppi della serie che contengono $P$ di $f$, contengono i $\nu-1$ punti \emph{coniugati} con $P$. Ogni gruppo della nuova  \emph{Serie (composta)} contiene $n\nu$ punti (\underline{$n$} gruppi, ciascuno di $\nu$ punti coniugati): $g^{r}_{n\nu}$.
%---------------------------------64------------------------------------------%
\section{Trasformazione di una curva contente una $g^{2}_n$ o $g^{3}_n$ semplice in un curva piana o sghemba d'ordine $n$.}
\label{sec:13}

\begin{teo} 
Data su una $f$ una $g^{2}_n$ semplice priva di punti fissi, esiste una trasformazione birazionale $f \to f'_n$, e la $g^{2}_n$ nella $g^{2}_n$ segata su $f'$ dalle rette del piano. 
\end{teo}
%\begin{proof}
Sia la $g^{2}_n$ su $f$ segata da $\lambda L+\mu M+\nu N=0$. Si applichi la trasformazione 
$x'=\frac{L}{N}, \, y'=\frac{M}{N}$; $f \to f'$. La rete si trasforma in $\lambda x'+ \mu y'+\nu=0$. Come si vede questa rete di rette sega su $f'$ la serie trasformata; così $f'$ ha ordine $n$. La trasformazione data è razionale, ma anche birazionale; ch\'{e} se a $P'$corrispondessero più $P$, la serie data sarebbe composta, contrariamente all'ipotesi.\\
Così per la Geometria sopra le Curve ogni $g^{2}_n$ semplice, priva di punti fissi, si può riguardare segata dalle rette del piano su una $f'_n$.\\
%\end{proof}
Esempio: data su una retta la serie $\lambda x^2+\mu x  +\nu=0$ di tutte le coppie di punti, la trasformazione birazionale $x'=x^2, ~ y'=x$ è quella che serve, e che come si vede muta la retta nella conica $x'=y'^2$.\\
Altro esempio: data su una retta la serie $\lambda x^3+\mu x^2+\nu=0$\quad $g^{2}_3$, applicando la trasformazione $x'=x^3, ~ y'=y^2$ raggiungiamo lo scopo; i nuovi gruppi stanno sulla curva $x'^2=y'^3$ (Parabola semicubica).
\begin{teo} 
\label{teo:sec}
Data su una $f$ una $g^{3}_n$ semplice, priva di punti 
%---------------------------------65------------------------------------------%
fissi, esiste una trasformazione birazionale $f \to f'_n$ sghemba, e $g^{2}_n \to g^{2}_n$ segata su $f'$ dai piani dello spazio.
\end{teo}
\noindent
La dimostrazione è analoga alla precedente:
%\begin{proof}
La $g^{2}_n$ data sia segata da $\lambda_0\varphi_0+\dots+\lambda_3\varphi_3=0$. Si applichi la trasformazione $x'=\frac{\varphi_0}{\varphi_3}, \, y'=\frac{\varphi_1}{\varphi_3},\, z'=\frac{\varphi_2}{\varphi_3}$. Allora $f$ si muta nella sghemba $f'_n$, in cui la nuova $g^{3}_n$ è segata da $\lambda_0x'+\dots+\lambda_3=0$, piani dello spazio. Anche qui la trasformazione è birazionale.\\
%\end{proof}
Così abbiamo introdotte le serie lineari sulle curve sghembe. Le stesse considerazioni si possono ripetere per $g^{4}_n$ ecc; si ottengono allora curve di iperspazi.\\
Esempio: data su una retta la serie $\lambda x^3+\lambda_1x^2+\lambda_2x+\lambda_3=0$ di tutti i gruppo di tre punti, la trasformazione birazionale ricorrente è $x'=x^3, \, y'=x^2, \, z'=x$, che rappresentano parametricamente una cubica sghemba, su cui la nuova $g^{3}_{3}$ è segata dai $\pi$ dello spazio.

\section{Serie completa - Determinazione di essa quando ne sia noto un gruppo - Somma o differenza di due serie - Serie residue.}
\label{sec:14}
Data su una curva un $g^{r}_{n}$, si può questa serie ampliare? Cioè costruire sulla curva una $g^{s>r}_{n}$ che contenga la serie data? Se il sistema $\infty^r$ di curve è contenuto in un altro
%---------------------------------66------------------------------------------%
più ampio, es. $\infty^{r+s}$, allora l'ampliamento è evidente. Ma il fatto che il sistema $\infty^r$ di curve non si possa ampliare, non vuol dire che $g^{r}_{n}$ non sia contenuta in una $g^{s>r}_n$; poich\'{e} si potrebbe ampliare un altro sistema lineare di curve col quale si può segare la stessa serie $g^{r}_n$.\\
Esempio: cubica con un punto doppio. Le rette del piano vi segano una $g^{2}_3$; le coniche per $O$ e $P$ una $g^{3}_3$ che contiene la $g^{2}_3$. \\

\begin{figure}[!htbp]
\centering
\includegraphics[scale=0.5]{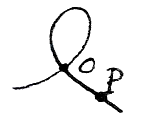}
%\caption{}
\end{figure}
\noindent
Data una $g^{r}_n$, se non esiste una $g^{s>r}_n$ che la contenga, la $g^{r}_n$ si dice \emph{Serie completa}.

\begin{lem}
Se due serie  lineari dello stesso ordine hanno un gruppo comune, sono contenute in una stessa serie lineare di quell'ordine.
\end{lem}

%\begin{proof}
\begin{figure}[!htbp]
\centering
\includegraphics[scale=0.5]{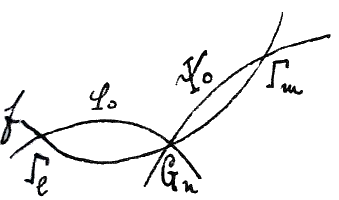}
%\caption{}
\end{figure}
\noindent
$g^{r}_n$ sia segata dalle $\sum \lambda_i\varphi_i$, $g^{s}_n$ sia segata dalle $\sum \lambda'_i\psi_i$.\\
Si consideri la serie:
\begin{equation}
\label{eq:pro1}
\psi_0 \sum \lambda_i\varphi_i+\varphi_0\sum\lambda'_i\psi_i=0
\end{equation}
Una curva \eqref{eq:pro1}, per esempio $\psi_0\varphi_1$ ha su $f$ i punti $\Gamma_l+G_n+\Gamma_m+G'_n$; un'altra curva \eqref{eq:pro1}, per esempio $\varphi_0\psi_1$ ha su $f$ i punti $\Gamma_l+G_n+\Gamma_m+G^{\prime\prime}_n$. 
Per tutte le \eqref{eq:pro1}, trascurando i punti fissi $\Gamma_l+G_n+\Gamma_m$, esse determineranno su $f$ una $g_n$, che come si vede usando opportunamente dei parametri, contiene la $g^{r}_n$ e la $g^{s}_n$.
%\end{proof}
Quindi: ogni serie completa contiene ogni serie di ugual
%---------------------------------67------------------------------------------%
ordine con un gruppo comune. Se due serie complete di ugual ordine hanno un gruppo comune coincidono; e allora:
\begin{teo} 
Dati $n$ punti su $f$, esiste una ed una sola serie lineare completa di quell'ordine che lo contiene (due coinciderebbero).
\end{teo}

%dubbio: questa più che una sub section sembra una definizione%
\subsubsection{Gruppi Equivalenti.}
$A \equiv B$ quando appartengono ad una stessa serie lineare.\\
$$A \equiv B ~\Rightarrow~ A+\Gamma\equiv B+\Gamma, \quad A-K\equiv B-K$$

\subsection{Serie Residua di $g^{r}_n$ rispetto a $\nu<n$ punti.}
Considero tutti i gruppi della $g^{r}_n$ che contengono questi $\nu$ punti; sopprimendoli, essi costituiscono evidentemente una $g^{r-\rho}_{n-\nu}$ ($0\leq \rho \leq \nu$), \emph{Residua} della $g^{r}_n$ rispetto ai $\nu$ punti.
\begin{teo} 
Se $g^{r}_n$ è completa, è completa anche la residua $g^{r-\rho}_{n-\nu}$.
\end{teo}
%\begin{proof}
\noindent
Sia $G_{n-\nu}$ un gruppo di $g^{r-\rho}_{n-\nu}$; basta far vedere che ogni $G'_{n-\nu}\equiv G_{n-\nu}$ fa parte della $g^{r-\rho}_{n-\nu}$.

$G_{n-\nu}+\Gamma_\nu\equiv G'_{n-\nu}+\Gamma_\nu$; anche questo secondo gruppo fa parte della serie completa, quindi $G'_{n-\nu}$ è un gruppo della serie residua di $\Gamma$ rispetto alla $g^{r}_n$, cioè di $g^{r-\rho}_{n-\nu}$.\\
%\end{proof}

\subsection{Somma di due serie complete.}
Date due serie lineari complete, tutti i gruppi ottenuti unendo una della prima con una della seconda in tutti modi possibili, appartengono ad una stessa serie lineare che resa completa si dice Somma delle due; ciascuna di queste è residua dei gruppi dell'altra rispettivamente alla seria somma

%dubbio: c'è un'annotazione in matita poco leggibile. Forse: vedi esempio alle ultime pagine del paragrafo 28, 2¡ quaderno%
%---------------------------------68------------------------------------------%
\section{Costruzione d'una serie completa mediante curve aggiunte - Teorema del Resto.}
\label{sec:15}
\emph{Curve aggiunte} ad una $f_n=0$ data senza punti multipli singolari, si chiamano tutte le curve che passano con molteplicità $(\alpha-1)$ almeno per ogni punto multiplo secondo $\alpha$ per $f$.

\begin{teo} 
La totalità delle aggiunte d'uno stesso ordine, $\varphi_k$, alla $f=0$, segano su questa, fuori dei punti multipli e base, una serie completa.
\end{teo}
%\begin{proof}
\begin{figure}[!htbp]
\centering
\includegraphics[scale=0.5]{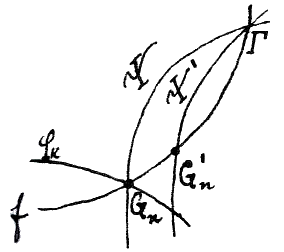}
%\caption{}
\end{figure}
\noindent
Basta dimostrare che ogni $G'\equiv G$ può segarsi con una $\varphi$.
$G'\equiv G$ vuol dire che sono segati da una stessa serie lineare, es. da $\psi$ e $\psi'$, con un gruppo base $\Gamma$ che comprende anche i punti multipli (se non li comprendesse, basterebbe addizionare alle $\psi$ altre curve per essi), pei quali le $\psi$ passino $\beta$ volte, la $f$ ~$\alpha$, le $\varphi$~ $(\alpha-1)$. Ricordando i Teoremi di  N\"{o}ether si vede subito ($\varphi\psi'\equiv Af+B\psi$) che $B$ è d'ordine $k$, aggiunta a $f$; e sega questa in $G'$.
%\end{proof}

\subsection{Costruzione della serie lineare completa di cui fa parte $G_n$.}
Per $G_n$ conduco a $f$ una aggiunta di ordine conveniente, $\varphi$, che segherà oltre $G_n$ un altro gruppo $\Gamma$. Tutte le aggiunte dello stesso ordine di $\varphi$ ad $f$, per $\Gamma$, vi segano fuori di $\Gamma$ e dei punti multipli la $g_n$ completa di cui fa parte $G_n$. Viceversa le aggiunte per $G$ \dots. Cioè:
%---------------------------------69------------------------------------------%
\begin{teo} [Teorema del Resto]
Le serie per $\Gamma$ e $G$ sono mutualmente residue.
\end{teo}
\noindent
\textbf{Esempio:} consideriamo una cubica con un punto doppio.
\begin{figure}[!htbp]
\centering
\includegraphics[scale=0.5]{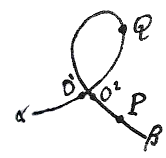}
%\caption{}
\end{figure}
Curve non costrette a passarvi (es. non aggiunte), perchè vi passino sono costrette ad una sola condizione; il passaggio per $O'$ vuol dire anche il passaggio per $O^2$ (il punto doppio si consideri come composto di due punti sui due rami). Le coniche aggiunte per tre punti due in $O$, dànno una $g^3_3$ completa. Vogliamo quella di cui fanno parte $O',P, Q$. Si conduce per $O'PQ$ una aggiunta (conica) che passi due volte per $O'$ (ramo $\alpha$), es. tangente al ramo $\alpha$; per la intersezione residua e per $O'$ si conducono le coniche aggiunte; esse segano la serie di cui fa parte $O'$ e non $O^2$.\\ 
Così, se si avesse un punto triplo: curve che non vi passano, dànno una serie per cui dànno una sola condizione pel passaggio; curve che vi passano semplicemente, perch\'{e} vi ripassino sono dunque costrette a due condizioni; curve aggiunte dànno tre condizioni indipendenti. Sicch\'{e} si ha l'opportunità sopra per scegliere i gruppi di punti.

\section{Genere di una curva - Sua relazione coll'ordine e la dimensione di una serie completa - Sua invarianza per trasformazioni birazionali.}
\label{sec:16}
%---------------------------------70------------------------------------------%
\subsection{Genere.}
\emph{Genere} di $f_m$ con soli punti multipli ordinari:
\begin{equation}
\label{eq:gen1}
p=\frac{(m-1)(m-2)}{2}-\sum\frac{\alpha(\alpha-1)}{2}
\end{equation}
I punti multipli singolari si possono riguardare punti multipli ordinari cui sono avvicinati indipendentemente altri punti multipli.
\begin{es}
Una cuspide conta come un punto doppio ordinario: $$p=\frac{(m-1)(m-2)}{2}-1;$$ il Tacnodo come due punti doppi infinitamente vicini: $$p=\frac{(m-1)(m-2)}{2}-2.$$
\end{es}

\subsection{Sua relazione con $n_i$ ed $r_i$ d'una completa $g_{n_i}^{r_i}$ segata su $f$ dalle aggiunte $\varphi_{m-3+i}$.}
\begin{equation}
\label{eq:gen2}
n_i=m(m-3+i)-\sum\alpha(\alpha-1)=2p-2+mi
\end{equation}
Dimensione del sistema delle $\varphi$:
\begin{equation*}
\rho_i=\frac{(m-3+i)(m+i)}{2}-\sum\frac{\alpha(\alpha-1)}{2}=p-1+mi+\frac{i(i-3)}{2}
\end{equation*}
Questa sarebbe $r_i$ se nessuna $\varphi$ contenesse la $f$; ma proprio ciò avviene appena sia $i \geq 3$. Se una $\varphi_{m-3+i}$ contiene la $f$, la curva residua sarà una $\mathcal{C}_{i-3}$ non soggetta a nessuna condizione. $f$ con una \emph{qualunque} $\mathcal{C}_{i-3}$ dànno complessivamente una aggiunta. Dunque le $\varphi$ linearmente indipendenti che contengono la $f$, sono tante quante le curve linearmente indipendenti d'ordine $i-3$ non soggette a nessuna condizione, cioè $\frac{(i-3)i}{2}+1$; quindi
\begin{equation}
\label{eq:gen3}
r_i=p-2+mi
\end{equation}
da cui per $i\geq 3$
\begin{equation}
\label{eq:gen4}
n_i-r_i=p.
\end{equation}

%dubbio: titolo della subsection?%
\subsection{Col mezzo delle aggiunte.}
Col mezzo delle aggiunte d'ordine abbastanza alto, possiamo segare su $f_m$ una $g^{r}_{n<n_i}$ completa qualsiasi: della $g^{r_i}_{n_i}$ si consideri la residua rispetto a $\nu$ punti; si avrà una
%---------------------------------71------------------------------------------%
$g^{r=r_i-\rho\geq r_i-\nu}_{n=n_i-\nu}$; quindi 
\begin{equation}
\label{eq:gen5}
 n-r\leq p
\end{equation}
Il genere è il massimo valore raggiunto dalla differenza tra $n$ ed $r$ d'una $g^{r}_n$ completa esistente sulla curva.

\begin{teo}[Teorema di Riemann]
Il genere è invariante per trasformazioni birazionali.
\end{teo}
%\begin{proof}
\noindent
Infatti una completa $g^{r}_n$ viene trasformata in una completa $g^{r}_{n}$; quindi non varia $n-r$, e neanche $p$ che ne è il massimo valore.\\
$p=p'$ è condizione necessaria, non sufficiente, per corrispondenze birazionali.
%\end{proof}
\section{Il genere è sempre $\geq 0$ - Curve razionali - Curve di genere $1$.}
\label{sec:17}
$p\geq n-r\geq 0$; se $f_m$ ha soli punti doppi, l'essere $p\geq0$ significa che il numero dei punti doppi è $\frac{(m-1)(m-2)}{2}$ al massimo. Questa proposizione sarebbe da aggiungere alle formole di Pl\"{u}cker.

\subsection{Curve razionali.}
Sono le curve riferibili birazionalmente con una retta. Per esse vale il
\begin{teo} [Teorema di L\"{u}roth]
Se le coordinate $(x,y)$ di un punto generico della curva si possono riguardare funzioni razionali di un parametro $t$, allora o $t$ è funzione razionale di $x,y$, oppure si può introdurre un nuovo parametro $t'$, funzione razionale di $t$, in modo che $x,y$ siano funzioni razionali di $t'$, e $t'$ di $x$ ed $y$.
\end{teo}
%---------------------------------72------------------------------------------%
\noindent
Esempio di curva razionali sono le coniche. Oppure le curve d'ordine $n$ con un punto multiplo d'ordine $(n-1)$; assumendolo come origine, equazione di una curva: $\varphi_{n-1}+\varphi_n=0$ (le $\varphi$ sono polinomi omogenei in $x, y$ di gradi $n-1$ ed $n$); tirando per l'origine una retta $y=tx$, intersezioni son la curva: $x^{n-1}P(t)+x^nQ(t)$, cioè $x=-\frac{P}{Q}$, e quindi $y=-\frac{tP}{Q}$. La Teoria della proiettività sulla retta si trasporta subito a tutte le curve razionali; per le coniche è nota.

\begin{teo} 
Condizione necessaria e sufficiente perchè una curva sia razionale è che sia $p=0$.
\end{teo}
%\begin{proof}
\emph{La condizione è sufficiente}: infatti, per ogni completa $g^{r}_n$ sarà $n=r$, cioè si hanno serie complete $g^n_n$. Così $\varphi-\lambda \psi$ segherà su $f$ fuori dei punti fissi una $g^1_{\underline{1}}$; le coordinate del punto variabile si trovano razionalmente da $f=0$, $\varphi+\lambda \psi=0$
$$x=R(\lambda) \qquad y=\rho(\lambda) \qquad \lambda=\frac{\varphi}{\psi}.$$
Viceversa la condizione \emph{è anche necessaria}: da queste ultime $\varphi-\lambda\psi=0$ sega una $g^1_{\underline{1}}$; si avrà in ogni caso $g^n_{\underline{n}}$ o $p=0$.\\
%\end{proof}
%dubbio: c'è un'annotazione a matita g^{p-1}_{tp-r}\equiv g^0_0%
\subsection{Curve di genere $1$.}
$n-r\leq 1$ \quad$r\geq n-1$ ($r\neq n$, curva razionale altrimenti) quindi $r=n-1$. Su curve di genere $1$ non esistono $g^n_n$ complete; evidentemente non si hanno punti fissi, es. $g^2_3$; e sono anche serie semplici: ch\'{e} se il passaggio per $P$ volesse anche quello per $Q$, fissato ad esempio $P$ (nella $g^2_3$) 
%---------------------------------73------------------------------------------%
si avrebbe residua una $g^1_2$ con $Q$ fisso, cioè una $g^1_1$.\\
Dunque la $g^2_3$ è semplice, priva di punti fissi, e permette la trasformazione di $f$ in $f'_3$ in cui la serie trasformata è segata dalle rette del piano. Così ogni curva di genere $1$ si può trasformare birazionalmente in una cubica piana senza punti doppi.\\
Dunque, tipo di curva con $p=1$ è la cubica piana. Si noti, due cubiche piane
%dubbio: nel manoscritto abbreviato piane con $\pi$%
 generiche non sono riferibili birazionalmente, che $p=p'$ è condizione necessaria, non sufficiente, per corrispondenza birazionale. Così le curve con $p=1$ non formano un'unica famiglia per trasformazioni birazionali; si possono prendere per modelli $\infty^1$ cubiche, ognuna dà famiglie distinte, ciascuna con un determinato valore (reale o complesso) di un birapporto (v. \S \ref{sec:22}).

\section{Serie speciali e non speciali; serie canonica; Ogni serie speciale può essere segata da aggiunte d'ordine $\leq m-3$ sopra la $f_m$, e viceversa.}
\label{sec:18}

Prendiamo ora $p=2$; $n-r\leq2$; quindi fissato $n$ si hanno due possibilità per $r$: $r=n-2$, $r=n-1$ (non $n=r$, curve razionali); cioè $n-r=1$, $n-r=2$.\\
Per esempio: quartica con punto doppio; le coniche aggiunte per tre punti vi segano una $g^1_3$; le rette pel punto doppio una $g^1_2$.\\
Vuol dire che appena $p>1$, esistono serie per cui
%---------------------------------74------------------------------------------%
$n-r<p$ (\emph{Serie speciali})\footnote{Questa classificazione ha carattere  invariante per trasformazioni birazionali.}, e serie con $n-r=p$ (\emph{non speciali}).
\begin{teo} 
Le aggiunte $\varphi_{m-3}$ a $f_m$ vi segano fuori dei punti multipli una serie speciale (\emph{Canonica}).
\end{teo}
%\begin{proof}
\noindent
Ordina (come al \S~ \ref{sec:16},
%dubbio: incerta la lettura del riferimento al paragrafo nel manoscritto 15 o 16?%
 tenendo conto che qui $i=0$):
\begin{equation}
\label{eq:can1}
n=m(m-3)-\sum \alpha (\alpha-1)=2p-2
\end{equation}
Dimensione (nessuna $\varphi$ contiene la $f$):
\begin{equation}
\label{eq:can2}
r=\frac{m(m-3)}{2}-\sum \frac{\alpha(\alpha-1)}{2}+\varepsilon=p-1+\varepsilon
\end{equation}
(c'è $\varepsilon$ perch\'{e} le condizioni imposte dai punti base possono non essere indipendenti). Allora le $\varphi$ segano su $f$ la $g^{p-1+\varepsilon}_{2p-2}$ speciale (\emph{canonica}).
%\end{proof}

\begin{teo} 
Ogni serie contenuta nella canonica è speciale.
\end{teo}
%\begin{proof}
\noindent
Infatti per la canonica $n-r<p$; serie residua rispettivamente a $\nu$ punti, basta togliere ad $n\to \nu$, ad $r$ un numero minore di $\nu$; si avrà quindi una $n_1-r_1<p$.\\
%\end{proof}

\begin{lem}
Dati $n$ punti su $f_m$, se esiste una aggiunta $\varphi_{m-3}$ che passi per $(n-1)$ di essi e non per l'ennesimo, questo risulta fisso nella completa $g_n$ determinata da quel gruppo.
\end{lem}
\begin{figure}[!htbp]
\centering
\includegraphics[scale=0.5]{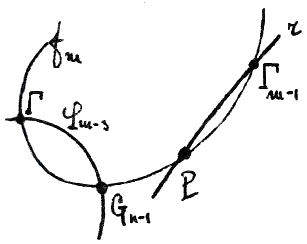}
%\caption{}
\end{figure}
%\begin{proof}
\noindent
La $\varphi_{m-3}$ con la retta $r$ per $P$, dànno una $\varphi_{m-2}$; per avere la nostra completa $g_n$, basta considerare le $\varphi_{m-2}$ per $\Gamma$ e $\Gamma_{m-1}$; ognuna di esse, avendo $m-1$ punti allineati, si spezza necessariamente sulla retta $r$ per $P$, ed in una residua $\varphi_{m-3}$.\\
Quindi $P$ risulta fisso in questa $g_n$. Così se la $g_n$ non ha
%---------------------------------75------------------------------------------%
punti fissi, ma $\varphi_{m-3}$ per $(n-1)$ passa per l'ennesimo. Cioè gli $n$ punti presentano ad essa $\varphi_{m-3}$ $(n-1)$ condizioni indipendenti, e non $n$.
%\end{proof}

\begin{teo} 
Ogni serie speciale è contenuta nella canonica, cioè per ogni suo gruppo si può condurre almeno una $\varphi_{m-3}$.
\end{teo}
\begin{figure}[!htbp]
\centering
\includegraphics[scale=0.5]{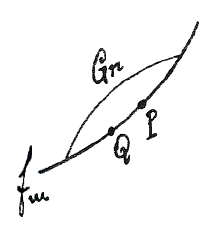}
%\caption{}
\end{figure}
%\begin{proof}
\noindent
Sia dunque $n-r<p$. Prendo un gruppo, $n$ punti ($G_n$) pei quali voglio dimostrare passa la $\varphi_{m-3}$. Nel gruppo prendo $P$ non comune a tutti i gruppi; serie residua $g^{r-1}_{n-1}$ speciale; vero il Teorema precedente per essa sarà vero per $g^r_n$, altrimenti $P$ sarebbe fisso. Così via, in $G_{n-1}$ si consideri $Q$ non comune a tutti i gruppi della $g^{r-1}_{n-1}$, ecc. \dots Ci riduciamo a vedere se il Teorema è vero per $g^0_{n-r}$, cioè per un solo gruppo: $n-r\leq p~ \Rightarrow ~n-r\leq p-1 $;
ora, come si vede dalla \eqref{eq:can2} le $\varphi_{m-3}$ dipendono da almeno $p-1$ parametri; dunque per $n-2$ punti (il gruppo), passa almeno una aggiunta $\varphi_{m-3}$.
%\end{proof}

\subsection{Altre definizioni.}
Serie speciali, possono segarsi con aggiunte $\varphi_{n\leq m-3} $. Serie speciali non possono segarsi solo con aggiunte $\varphi_{m-3}$. Ma così non appare il carattere invariantivo per trasformazione birazionale, ch\'{e} curve aggiunte non si trasformano in curve aggiunte. Appare dicendo: serie speciali contenute nella canonica; non speciali, non contenute.
%dubbio: incerta l'interpretazione dell'abbreviazione cont.%

\`{E} chiaro ora che: se $f_m$ ed $f'_{m'}$ si corrispondono birazionalmente, la serie canonica su $f_m$ corrisponde alla canonica
%---------------------------------76------------------------------------------%
su $f'_m$; non è però che le aggiunte si trasformino in aggiunte, ch\'{e} la trasformazione si opera sulla curva, e non sul piano.
\begin{teo} 
La serie canonica è $g^{p-1}_{2p-2}$.
\end{teo}
\noindent
Infatti, possa essere $g^p_{2p-2}$. Aggiungiamo $P$ a tutti i gruppi, si avrebbe una $g^p_{2p-1}$, speciale, quindi contenuta nella canonica, il ch\'{e} è impossibile.\\
 La serie canonica è caratterizzata dall'invarianza per trasformazioni birazionali: $f_m\to f'_{m'}$, ogni serie speciale $\to$ serie speciale; quindi la canonica si trasforma in una serie speciale che contiene tutte le altre, cioè canonica.

\section{Invariantività della serie canonica - Suoi caratteri - Caso della serie canonica composta - Curve iperellittiche.}
\label{sec:19}
Abbiamo trattato le prime due questioni al \S ~ precedente.
Esempio: Quartica generale: genere $3$; serie canonica $g^2_4$, segata dalle rette del piano; le speciali possono segarsi con rette di un fascio, generico ($g^1_4$) o per la curva ($g^1_3$). Vediamo su detta quartica una serie non speciale; per $4$ punti generici, si ha una completa $g_4$ non speciale ($4$ punti non stanno su una retta); 
%dubbio: incerta l'interpretazione dell'abbreviazione dimens.% 
segandola ad esempio con coniche aggiunte, $g^1_4$, completa non speciale.
\begin{teo} 
La serie canonica non ammette punti fissi comuni
%---------------------------------77------------------------------------------%
a tutti i suoi gruppi; cioè le $\varphi_{m-3}$ non hanno altri punti comuni fuori dei punti multipli.
\end{teo}
%\begin{proof}
\noindent
Se la $g^{p-1}_{2p-2}$ avesse $P$ fisso, aggiungendo alla residua $g^{p-1}_{2p-3}$ a tutti i gruppi $Q$ si avrebbe una $g^{p-1}_{2p-2}$; questa dovrebbe coincidere colla canonica, cioè questa avrebbe fisso anche $Q$; ciò è assurdo, ogni punto della curva sarebbe fisso per la serie canonica.
%\end{proof} 

\begin{teo} 
Condizione necessaria e sufficiente perch\'{e} una curva con $p>2$ abbia la serie canonica composta, è che la curva contenga una $g^1_2$ speciale.
\end{teo}
%\begin{proof} 
\begin{enumerate}[label=\emph{\roman*)}]
\item \emph{Se $g^{p-1}_{2p-2}$ è composta}, i gruppi per $P_1$, passano almeno per $P'_1$; ecc.; se $g^{p-1}_{2p-2}$ è composta, i gruppi per $P_{p-2}$, passano almeno per $P'_{p-2}$.\\ 
I gruppi che hanno comuni $P_1, P_2, \dots, P_{p-2}$, hanno così comuni $2(p-2)$ punti; sopprimendoli ho una $g^1_2$ speciale. (I detti gruppi non possono avere altri punti comuni, sarebbe ordine $\leq 0$).
\item \emph{Se la curva contiene una $g^1_2$ speciale}, questa sarà segata da un fascio di curve, parametro $\lambda$; ad ogni $\lambda$ due punti; con una retta di parametro $\lambda$, si ha allora una corrispondenza $(2,1)$.\\

\begin{figure}[!htbp]
\centering
\includegraphics[scale=0.5]{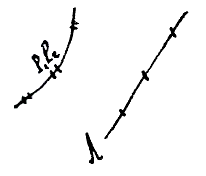}
%\caption{}
\end{figure}
Se allora sulla retta considero una $g^{p-1}_{p-1}$, sulla curva le corrisponderà $g^{p-1}_{2p-2}$ composta, ch\'{e} il gruppo contenente $P_1$ deve contenere $P_2$.

\end{enumerate}

%\end{proof}

\subsection{Curve iperellittiche.}
Nelle curve con $p=2$, la serie canonica è già $g^1_2$ speciale (curve ellittiche).
%---------------------------------78------------------------------------------%
Curve iperellittiche sono appunto quelle che hanno una $g^1_2$ speciale; per esse la serie canonica è composta, e solo per esse.

\subsubsection{Curva canonica di $f$.}
Abbiamo visto che data su una $f_m$ una $g^r_n$ semplice priva di punti fissi , si può trasformare birazionalmente, la $f_m$ in $f'_n$ dello spazio $S_r$ su cui gli iperpiani ($S_{r-1}$) segano la $g^r_n$ trasformata della primitiva. Applicando alla $g^{p-1}_{2p-2}$: ogni $f_m$ di genere $p$ non iperellittica si può trasformare birazionalmente in una $f'_{2p-2}$ di $S_{p-1}$ in cui gli iperpiani $S_{p-2}$ segano la serie canonica corrispondente. $f'$ è la curva canonica di $f$.

\section{Teorema di Riemann-Roch.}
\label{sec:20}
Diamo un risultato:\\
un gruppo di $n$ punti indipendenti (generici) d'una $f_m$ di genere $p$, se $n\leq p$ non ha gruppi equivalenti; se $n>p$ il gruppo appartiene ad una serie completa non speciale $g^{n-p}_n$.

\begin{teo}[Teorema di Roch]
Ogni gruppo $G_n$ d'una completa $g^r_n$ speciale presenta esattamente $n-r$ condizioni alle aggiunte $\varphi_{m-3}$ costrette a contenerlo.
\end{teo}
 %\begin{proof}
 \noindent
Si consideri $P_1$ di $G_n$, non comune a tutti i gruppi; le $\varphi_{m-3}$ per $G_{n-1}$ dovranno passare per $P_{1}$; in $G_{n-1}$ si consideri $P_2$, ecc. \dots\\
\[
G_n=P_1+P_2+\dots+P_r+G_{n-r}
\]
%---------------------------------79------------------------------------------%
Le $\varphi_{m-3}$ per $G_{n-1}$ dovranno passare per $G_n$ intero, cioè pel passaggio per $G_n$ occorrono al massimo $n-r$ condizioni; cioè le $\varphi_{m-3}$ per $G_n$ hanno la dimensione
\begin{equation}
\label{eq:rr1}
\rho=p-1-(n-r)+\varepsilon
\end{equation}
Serie residua di $g^r_n$ rispetto alla canonica è $g^{\rho}_{\nu=2p-2-n}$.\\
Relazioni tra i caratteri della $g^r_n$ e della $g^\rho_\nu$:
\[
\rho-r=p-1-n+\varepsilon \qquad n+\nu=2p-2
\]
Ora $g^r_n$ è anche residua di $g^\rho_\nu$ rispetto alla canonica:
\[
r-\rho=p-1-n+\varepsilon';
\]
sommando ($\varepsilon \geq 0 \leq \varepsilon'$):
\[
0=2p-2-(n+\nu)+\varepsilon+\varepsilon' ~ \Rightarrow~ \varepsilon=\varepsilon'=0
\]
Quindi nella \eqref{eq:rr1} $\varepsilon=0$, le condizioni sono esattamente $n-r$.\\
%\end{proof}

\begin{teo} [Teorema di Riemann]
Il Teorema è valido in un certo senso per le $g^r_n$ non speciali.
\end{teo}
%\begin{proof}
\noindent
Ch\'{e} dire che $G_n$ presenta $n-r=p=\frac{(m-1)(m-2)}{2}$ condizioni alle $\varphi_{m-3}$ costrette a contenerlo, essendo le $\varphi$ dipendenti da $\frac{(m-1)(m-2)}{2}-1<p$ parametri equivale a dire che nessuna $\varphi$ passa pel gruppo.
%\end{proof}

\begin{teo} [Teorema di reciprocità (Brill e di N\"{o}ther)]
Una completa speciale $g^r_n$ ammette come residua rispetto alla serie canonica una serie $g^{p-1-n+2}_{2p-2-n}$.
\end{teo}

\emph{Esempio:} consideriamo una curva sghemba, intersezione di
%---------------------------------80------------------------------------------%
una quadrica rigata con una superficie cubica: è una $\mathcal{C}_6$;
%dubbio: si confonde il 5 con il 6%
 il genere si valuta come $p$ di una piana in corrispondenza birazionale con essa. Proiettando la $\mathcal{C}_6$ da una suo punto su un piano,
 %dubbio: sciogliere $\pi$ in piano%?
  si ha una $\mathcal{C}_5$. Ogni retta di un quadrica è segata in $3$ punti dalla $\mathcal{C}_6$;
%dubbio: qui il 5 è messo ad apice, ma probabilmente è una svista%
per il punto di proiezione $S$ passano due rette della quadrica, quindi $\mathcal{C}_5$ ha due punti doppi, e quindi il genere è $4$. La serie canonica sarà la $g^3_6$. Vediamo di trovare due serie mutuamente residue rispetto alla serie canonica. Prendiamo una retta del primo sistema della quadrica ($3$ punti su $\mathcal{C}_6$); facendole girare intorno un piano si ha una $g^1_3$ speciale; ogni gruppo d'una serie canonica che contiene due punti della $g^1_3$ contiene il terzo; quindi il tre punti devono essere allineati, e sono su una retta del secondo sistema. Così si ha l'altra $g^1_3$ segata dalle rette del primo sistema della quadrica. Queste sono due serie mutuamente residue rispetto alla serie canonica, e corrispondono a quelle del Teorema di reciprocità. Se la quadrica fosse un cono 
%dubbio: incerta l'interpretazione del manoscritto%
c'è un solo sistema di rette (generatrici) e le $g^1_3$ coincidono.

\section{Numero dei punti doppi d'una $g^1_n$; cenno sulla formola che dà il numero dei punti tripli d'una $g^2_n$, ecc.}
\label{sec:21}
\begin{teo} 
Una curva possedente due serie $g'_m$ $g'_n$ che non 
%---------------------------------81------------------------------------------%
abbiano infinite coppie comuni si può sempre trasformare birazionalmente in una curva $\mathcal{C}_{m+n}$ dotata di un punto multiplo secondo $m$ ed uno secondo $n$. Le rette uscenti dall'uno o dall'altro di essi segano sulla nuova curva le serie trasformate o corrispondenti alle due primitive.
\end{teo}
%\begin{proof}
\noindent
Su $f$ siano $g'_m$ segata da 
\begin{equation}
\label{eq:rralpha}
\lambda=\frac{\varphi_1}{\varphi_2},
\end{equation}
$ g^1_n$ da 
\begin{equation}
\label{eq:rrbeta}
\mu=\frac{\psi_1}{\psi_2}
\end{equation}
\noindent
Ad $(x,y)$ corrisponde un solo valore di $(\lambda, \mu)$,
%dubbio: questa parte è cancellata ma sembra in matita...%
cioè due curve dei due fasci, cioè due gruppi delle due serie; se queste non hanno legami particolari, in genere le due curve $\lambda$ e $\mu$ per $P(x,y)$ non hanno in generale altri punti comuni (altrimenti, se avessero in comune anche $Q$, le due serie sarebbero composte mediante infinite coppie comuni).
%dubbio: qui la parentesi si apre ma non si chiude%
 Così in generale anche a $(\lambda, \mu) \to$ un solo $P(x,y)$. Al variare di $P$, $P'$ descrive $F$; $f$ ed $F$ si corrispondono birazionalmente.
Ad un valore di $\lambda$ corrispondono $m$ valori di $\mu$; sicch\'{e} nella $F$, $\mu$ comparirà al grado $m$; così $\lambda$ al grado $n$: allora ad $f(x,y)=0$ corrisponde $F(\lambda^n, \mu^m)=0$ si ottiene dalle \eqref{eq:rralpha} \eqref{eq:rrbeta} eliminando $x$ ed $y$. $F$ è di grado $m+n$.\\

\begin{figure}[!htbp]
\centering
\includegraphics[scale=0.4]{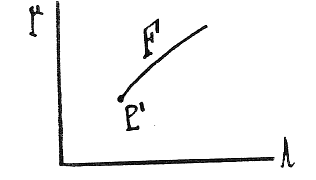}
%\caption{}
\end{figure}
\noindent
Ad un gruppo della data $g^1_m$ su $f$, corrisponde un valore di $\lambda$; quindi le rette $\lambda=cost$ segheranno su $F$ la $g^1_m$ corrispondente. Così le $\mu=cost$, segheranno la $g^1_n$ trasformata.
%---------------------------------82------------------------------------------%
Ora $F$ è d'ordine $m+n$, le $\lambda=cost$ segano $m$ punti variabili, vuol dire le altre intersezioni saranno fisse al centro del fascio; cioè l'infinito di $\mu$ sarà multiplo secondo $n$ per la curva $F$. Così ancora l'infinito di $\lambda$ sarà multiplo secondo $m$ per la curva $F$.
%\end{proof}

\begin{teo} 
Data su $f$ una $g^{\underline{1}}$ la differenza tra il numero dei punti doppi della serie, $\mu$, e il doppio dell'ordine della serie, $2m$, non varia sostituendo alla serie un'altra serie semplicemente infinita.
\end{teo}
%\begin{proof}
\noindent
Dimostriamo il Teorema per la $f$ con $g^1_m$ e $g^1_n$; intanto detta differenza non varia per trasformazioni birazionali, non variando n\'{e} $m$ n\'{e} $\mu$. Dunque trasformiamo la $f$ in $F$ come al Teorema precedente, punti multipli $S^{(n)}$ e $T^{(m)}$.\\
\begin{figure}[!htbp]
\centering
\includegraphics[scale=0.5]{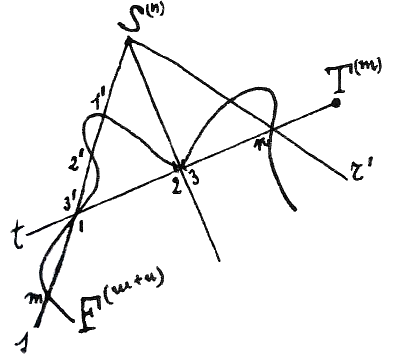}
%\caption{}
\end{figure}
 Posso stabilire tra due fasci per $S$ una corrispondenza algebrica così:\\
alla retta $s$ del primo fascio che passa pel punto $1$ di $t$, far corrispondere nel secondo le $m(n-1)$ rette che da $S$ vanno agli $m(n-1)$ punti d'incontro di $T1',T2', \dots, Tm$ coi punti della curva. E ad $r'$ del secondo fascio faccio corrispondere le $m(n-1)$ rette del primo\dots\\
Si ha la corrispondenza dei due fasci per $S$: [$m(n-1),m(n-1)$]. Per il Teorema di Chasles le rette unite saranno $2m(n-1)$. \emph{Possono cadere \emph{o}} nei punti di tangenza per $T$, come $2-3$;
%---------------------------------83------------------------------------------%
allora $t$ dà un gruppo della $g^1_n$ che ha due punti coincidenti; allora ivi c'è un punto doppio \emph{della serie}; i contatti delle rette da $T$ sono tanti quanti i punti doppi della $g^1_n$, e siano $\nu$. \emph{O} nei punti multipli della curva, e siano $d$, e ciascuno conti per $k$ coincidenze; allora deve essere
\[
2m(n-1)=\nu+kd;
\]
analogamente scambiando le due serie $g_m$ e $g_n$:
\[
2n(m-1)=\mu+kd;
\]
sottraendo membro a membro:
\begin{equation}
\label{eq:dopp1}
\mu-2m=\nu-2n
\end{equation}
%\end{proof}

\subsection{Punti doppi d'una $g^1_n$.}
Dunque $I=\nu-2n$ è invariante. Ma $\nu$ si può interpretare come classe di $F$; allora $\nu=n(n-1)-\sum \alpha(\alpha-1)$ e risulta
\begin{equation}
\label{eq:dopp2}
I=\nu-2n=2p-2; \quad \nu=2(n+p-1).
\end{equation}
Con la \eqref{eq:dopp2} \emph{posso definire} il genere della curva.\\
Così si può trovare: \emph{numero dei punti ($r+1$)-upli d'una $g^r_n$}:
\begin{equation}
\label{eq:dopp3}
\nu_r=(r+1)(n+rp-r)
\end{equation}

\section{Curve piane dei primi generi - Riduzione a tipi.}
\label{sec:22}
\subsection{Curve con $p=0$.}
La curva è razionale. Quindi la retta è il tipo delle curve con $p=0$, ch\'{e} si possono trasformare birazionalmente in una retta.

\subsection{Curve con $p=1$.}
Vedi anche \S ~\ref{sec:17}. Tali curve hanno serie speciali.
%---------------------------------84------------------------------------------%
Tre punti ad arbitrio danno la completa $g^2_3$ non speciale, semplice e priva di punti fissi; ricordare allora il \S ~ \ref{sec:13}. La cubica piana è il tipo delle curve con $p=1$.

\begin{teo} [Teorema di Salmon]
Da un punto d'una cubica piana si possono condurre quattro rette a toccar la cubica \emph{altrove} (formole di Poncelet); se il punto varia sulla curva, non varia il birapporto delle quattro rette.
\end{teo}
\begin{teo} 
Condizione necessaria e sufficiente perchè due cubiche piane si possano porre in corrispondenza birazionale è che i due birapporti siano uguali fra loro.
\end{teo}
\noindent
L'equazione d'ogni cubica piana, con scelta particolare delle coordinate proiettive, cioè con proiezione opportuna della cubica, in coordinate cartesiane, si può scrivere
\[
y^2=x^3+bx^2+cx+d
\]
Allora le coordinate $X,Y$ d'un punto generico d'una curva con $p=1$, si possono esprimere come funzioni razionali d'un parametro, e della radice d'un polinomio di terzo grado nel parametro:
\begin{teo} [Teorema di Clebsch]
$X=R(x, \sqrt P_3(x)), Y=S(x, \sqrt P_3(x))$. 
\end{teo}
\noindent
Sull'asse $x$ ($y=0$) la curva ha tre intersezioni; questi tre punti più il punto all'infinito di $x$ danno quattro punti il cui birapporto è il birapporto (modulo) della curva.

\subsection{Curve con $p=2$.}
Hanno la serie canonica $g^1_2$, e non altre serie speciali. Quattro punti generici (ma due si corrispondono nella $g^1_2$) dànno luogo ad una completa $g^2_4$ semplice priva di punti fissi, poich\'{e} deve essere $4-2=p=2$. Allora ricordare il \S~\ref{sec:13}. La curva si può trasformare in una
%---------------------------------85------------------------------------------%
quartica piana, che per avere $p=2$ deve avere un punto doppio. Anche qui due curve di $n=4$, perch\'{e} siano riferibili birazionalmente, le sei tangenti dai due punti doppi (esclusi) devono formare due gruppi proiettivi; si hanno sei condizioni.

\subsection{Curve con $p=3$.}
Serie canonica $g^2_4$, che può essere semplice, o composta mediante due $g^1_2$. Se semplice si può trasformare la curva in una curva piana del quarto ordine senza punti multipli. La quartica piana è il tipo più semplice cui si può ridurre una curva con $p=3$ non iperellittica. Al solito: se due quartiche sono riferibili birazionalmente, sono trasformabili fra loro mediante collineazione.

\subsubsection{Caso iperellittico.}
Prendo sulla curva cinque punti generici, ho una completa $g^2_5$ semplice e priva di punti fissi, e permette la trasformazione in una quintica su cui la corrispondenza $g^2_5$ sarà segata dalle rette del piano, e che vediamo ha un punto multiplo. Unisco $A$ con $A'$ di un gruppo; ho altri tre punti allineati, $B_1, B_2, B_3$. Serie residua di questi rispetto alla $g^2_5$, è la $g^1_2$, cioè infinite rette di un fascio per $B_1, B_2, B_3$; cioè $B_1, B_2, B_3$ si devono riunire (in senso proiettivo) sulla curva. Così la quintica deve avere un punto triplo. Effettivamente si ha allora la serie canonica composta.
%dubbio: presente un'annotazione a matita poco leggibile, forse: per g^1_2%

\subsection{Curve con $p=4$.}
Serie canonica $g^3_6$, semplicemente composta. Se la serie è semplice si può (la curva) trasformare birazionalmente 
%dubbio: presente un'annotazione a matita poco leggibile, forse: con una g^1_2%
in una sestica sghemba, che si può prendere come tipo. Vediamo che questa $\mathcal{C}^6$
%---------------------------------85------------------------------------------%
è intersezione completa di una quadrica con una cubica.\\
Infatti: la $\mathcal{C}^6$ è di genere $4$ nello spazio ordinario. Consideriamo la serie segatasi dalle $\infty^9$ quadriche dello spazio; dev'essere $12-r\leq p$, quindi $r$ può essere al massimo $8$, vuol dire vi è almeno una quadrica per la $\mathcal{C}^6$; così si vede anche che vi sono $5$ sup.
%dubbio: incerta l'interpretazione dell'abbreviazione sup. forse superiori%
del terzo ordine linearmente indipendenti per la $\mathcal{C}^6$, cioè $\infty^4$ cubiche per essa. Così c'è qualche cubica non contenente la quadrica che passa per la $\mathcal{C}^6$. E così è trovato ciò che volevamo.\\
Come tipo d'una curva iperellittica con $p=4$ conviene prendere la sestica con un punto quadruplo.

%-------------------------------pag. 87 -----------------------------------%
\chapter{Curve sghembe}

\section{Vari modi di rappresentare analiticamente una curva sghemba - Genere della curva sghemba - Serie lineari su di essa.}
\label{sec:23}

\subsection{$I^{a}$ Definizione.}
Se le curve sghembe si definiscono come intersezione di due superficie
%dubbio: incerta l'interpretazione di superf.%
 $f$ e $\varphi$, si vengono ad escludere quelle che non sono intersezioni complete. Allora alle $f$ e $\varphi$ basta aggiungere la proiezione piana della sghemba, per esempio $\psi(x,y)=0$. In genere, tre equazioni:
\begin{equation}
\label{eq:def1}
f(x,y,z)=0 \quad \varphi(x,y,z) \quad \psi(x,y,z)=0,
\end{equation}
se esiste una successione continua di punti che le soddisfa, questi costituiscono una Curva algebrica sghemba. Siccome le \eqref{eq:def1} potrebbero rappresentare anche punti isolati estranei, per escluderli si dimostra che quattro equazioni bastano.

\subsection{$II^a$ Definizione di Cayley.}
Eliminando $z$ nelle \eqref{eq:def1} si hanno due equazioni o curve piane su $(x,y)$. La parte comune (fattore comune) è la proiezione della sghemba, la cui equazione si trova dunque con operazioni razionali: $F(x,y)=0$. Una coppia $(x,y)$ che soddisfi la $F=0$, sostituita nelle \eqref{eq:def1}, queste devono fornire per $z$ un'unica soluzione comune; allora $z$ razionalmente
%-------------------------------pag. 88 -----------------------------------%
\begin{equation}
\label{eq:def2}
F(x,y)=0, \quad z=\frac{P(x,y)}{Q(x,y)}
\end{equation}
La prima è un cono, o cilindro (analiticamente) proiettante la curva, di ordine $\underline{n}$ (ordine di $F$ e di $\mathcal{C}$). Si vede subito che variando il centro di proiezione, $\mathcal{C}$ resta di grado $n$, cioè un piano pel nuovo centro sega lo stesso $n$ punti su $\mathcal{C}$. Così ogni piano dello spazio.
La seconda
\begin{equation}
\label{eq:def3}
P_\nu-z \,Q_{\nu-1}=0
\end{equation}
è una superficie d'ordine ad esempio $\nu$; $z$ è al primo grado, ad ogni $(x,y)$ un solo valore di $z$, cioè ogni parallela a $z$ ($x=cost$, $y=cost$) sega la \eqref{eq:def3} in un sol punto proprio. Ma ogni retta deve incontrare la \eqref{eq:def3} in $\nu$ punti; dunque nell'origine $S^\infty$ ci sarà per la superficie \eqref{eq:def3} un punto multiplo di molteplicità $(\nu-1)$; ogni piano per $S$ dà una curva d'ordine $\nu$ con in $S$ un punto multiplo di molteplicità $(\nu-1)$; le $\nu-1$ tangenti principali in $S$, al variare del piano sezionante descrivono un cono di vertice $S$ e grado $(\nu-1)$. 
La \eqref{eq:def3} è un \emph{Monoide} di vertice $S$.
\begin{teo} 
Il monoide contiene $\nu(\nu-1)$ rette che escono dal vertice e giacciono per intero su di esso. Sono le rette comuni ai due cilindri $P_\nu=0$ \quad $Q_{\nu-1}=0$.
\end{teo}
Ogni curva algebrica sghemba si può rappresentare come intersezione delle \eqref{eq:def2}, cioè di un cono che la proietta da un punto generico dello spazio con un monoide avente ivi il vertice. 
Però la intersezione
%-------------------------------pag. 89 -----------------------------------%
completa di queste due superficie si compone, oltre che della curva, di $n(\nu-1)$ rette costituenti la intersezione del cono proiettante ($F=0$) col cono tangente al monoide nel vertice ($Q=0$).
Infatti la intersezione delle \eqref{eq:def2} è di grado $n\nu$, mentre la nostra $\mathcal{C}$ è di grado $n$; vi sarà un'intersezione residua d'ordine $n(\nu-1)$.
\begin{figure}[!htbp]
\centering
\includegraphics[scale=0.5]{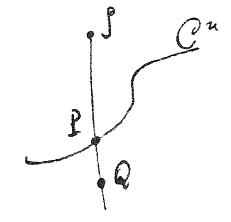}
%\caption{}
\end{figure}
\noindent
Se $Q$ appartiene a questa residua, $SQ$ è retta del cono ($F$), deve segare $\mathcal{C}$ in $P$; ma allora $SQ$ appartiene anche al monoide (perch\'{e} ha $\nu+1$ intersezioni con esso). Così la intersezione residua si compone di $n(\nu-1)$ rette. Queste generatrici, stanno su $Q_{\nu-1}=0$ (ch\'{e} danno $z=\infty$), stanno su $F_n=0$, quindi costituiscono \emph{tutte} le intersezioni dei due coni. ($Q_{\nu-1}=0$ è un cono, le cui generatrici dànno, [\eqref{eq:def2}], $z=\infty$; cioè esse hanno $\nu$ intersezioni col monoide, d'ordine $\nu$, in $S^\infty$; e $Q=0$ è il cono tangente al monoide nel vertice).

Tra le rette proiettanti $\mathcal{C}$ da $S$, o rette del cono proiettante, vi sono certo delle generatrici doppie; queste stanno anche sul monoide (vi hanno $\nu+1$ intersezioni), e sono generalmente semplici per esso. Queste generatrici, fanno parte evidentemente delle $n(\nu+1)$ rette viste sopra; esse stanno sul monoide \eqref{eq:def3}, stanno su $Q=0$, staranno su $P=0$. Vuol dire, nella proiezione piana, $P$ e $Q$ sono aggiunte ad $F$.

L'intersezione residua si compone esclusivamente di rette,
%-------------------------------pag. 90 -----------------------------------%
in numero $l$. Allora ordine della curva sghemba, $n\nu-l=N$. Allora si vede subito che una curva sghemba d'ordine $N$, da ogni superficie d'ordine $\rho$ che non la contenga vien segata in $N\rho$ punti.

\subsection{Altra definizione.}
Ogni curva algebrica sghemba si può riguardare, nel modo più perfetto, come luogo di un punto le cui coordinate variano come funzioni razionali di due parametri legati da una relazione algebrica:
\begin{equation}
\label{eq:def4}
F(\xi,\eta)=0\qquad x=\frac{L(\xi,\eta)}{P} \qquad y=\frac{M}{P} \qquad z=\frac{N}{P}
\end{equation}
In particolare, se la corrispondenza della $F$ (piana), con la $\mathcal{C}$ è per proiezione, se il centro di proiezione è all'infinito di $z$, se $x=\xi, y=\eta$, si ritrova la rappresentazione monoidale.\\
Si può ritrovare il secondo Teorema del \S~ \ref{sec:13}.
%dubbio: aggiungere il riferimento al teorema secondo la numerazione di Latex? \ref{teo:sec}%

Se $F$ contenente una $g^3_n$ si trasforma in  una $F'$, è facile vedere che $F$ ed $F'$ sono collineari.\footnote{Esempio di rappresentazione monoidale. Curva del terzo ordine, intersezione d'un monoide del secondo ordine con un cono del terzo ordine; intersezione residua tre rette, due coincidenti.}
%dubbio:  non so dove riferire la nota%

\section{Condizione perch\'{e} una curva piana d'ordine $n$ e genere $p$ sia proiezione d'una curva sghemba dello stesso ordine.}
\label{sec:24}
%dubbio: nel titolo a pag. 90 c'è indicato l'abbreviazione $\pi$ per piana che invece nell'indice è sciolta%
%...ogni volta che c'è \pi và sciolto in piana??%

\begin{teo} 
Se $F_\pi$ è proiezione di $\mathcal{C}$ sghemba, la $g^2_n$ segata dalle rette del
%-------------------------------pag. 91 -----------------------------------%
piano su $F$ è contenuta in una $g^3_n$.
\end{teo}
\begin{figure}[!htbp]
\centering
\includegraphics[scale=0.5]{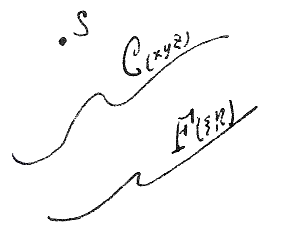}
%\caption{}
\end{figure}
\noindent
La rappresentazione monoidale di $\mathcal{C}$ può scriversi:
\[
F(\xi,\eta)=0 \qquad x=\frac{\xi Q}{Q} \qquad y=\frac{\eta Q}{Q} \qquad z=\frac{P}{Q}
\]
Ai punti che i piani dello spazio
\[
\lambda x+\mu y+\nu z+\pi=0
\]
segano su $\mathcal{C}$, corrispondono i punti che su $F$ segano le curve (fatta la trasformazione birazionale indicata dalla proiezione):
\[
(\lambda \xi+\mu \eta+\pi)Q_{\nu-1}+\nu P_\nu=0
\]
Ed entro questo sistema c'è, come si vede, quello delle rette piane ($\nu=0$, e trascurate le intersezioni fisse di $F$ con $Q$).

\begin{teo} 
Se la $g^3_n$ segata su $F$ dalle rette del piano è contenuta in una $g^3_n$, la $F$ è proiezione d'una curva sghemba. 
\end{teo}

\begin{figure}[!htbp]
\centering
\includegraphics[scale=0.5]{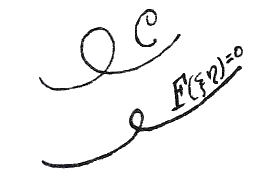}
%\caption{}
\end{figure}
\noindent
Per costruire la $g^3_n$ basta considerare un gruppo  della $g^2_n$, farvi passare una aggiunta d'ordine conveniente $\nu$, e considerare tutte le aggiunte $\varphi_\nu$ pel gruppo residuo; esse segheranno su $F$ la serie voluta $g^{3 \text{ almeno}}_n$ completa, di cui fa parte la $g^2_n$. Tra le $\varphi_\nu$ prendiamo $(\lambda \xi + \mu \eta +\pi)Q_{\nu-1}$ che segano su $F$ la $g^2_n$. Un'altra aggiunta sia $P_\nu$; il sistema completo delle aggiunte sarà:
\begin{equation}
\label{eq:pro1}
(\lambda \xi+\mu \eta +\pi)Q_{\nu-1}+\nu P_{\nu}=0
\end{equation}
che sega su $F$ la nostra $g^3_n$; allora dal \S ~\ref{sec:13}
%-------------------------------pag. 92 -----------------------------------%
la $F$ si può trasformare birazionalmente in una sghemba, basta porre (si metta $\pi Q$ in ultimo) 
\[
x=\xi \qquad y=\eta \qquad z=\frac{P}{Q}
\]
che è appunto la rappresentazione monoidale d'una sghemba (proiezione). \`{E} il procedimento inverso del Teorema precedente. 
Dunque: condizione necessaria e sufficiente perch\'{e} una curva piana $F_n$ sia proiezione d'una sghemba $\mathcal{C}_n$ è che la $g^2_n$ dalle rette del piano sia incompleta.
%dubbio: nel manoscritto d. \pi, intende della curva piana o del piano?%

\begin{teo} 
Condizione necessaria e sufficiente perch\'{e} una $\mathcal{C}_n$ di genere $p$ sia proiezione d'una sghemba $\mathcal{C}_n$, se $n-p\leq2$ è che i punti multipli della $F$ presentino condizioni non tutte indipendenti alle aggiunte $\varphi_{n-4}$. Se $n-p>2$ sempre $F$ è proiezione d'una $\mathcal{C}$ sghemba.
\end{teo}

\begin{enumerate}
\item \underline{Se $n-p>2$}; la $g^2_n$ dovrebbe essere incompleta, e difatti è così; ch\'{e} per una completa $r\geq n-p>2$
\item  \underline{Se $n-p\leq2$}; \emph{la condizione è necessaria}. Se $F$ è proiezione di $\mathcal{C}$, allora esiste una $g^3_n\supset g^2_n$;
%dubbio: è un maggiore?%
è speciale,
%dubbio: sembra aggiunto a matita anche la g^2_n%
quindi per $n$ punti allineati della $F$ (un gruppo della $g^2_n$ segata dalle rette di $\pi$) si può condurre almeno un'aggiunta $\varphi_{n-4}$. $g^2_n$ è contenuta dunque in una $g^{3+\varepsilon}_n$ che è completa speciale. Pel Teorema di Riemann-Roch, ogni $G_n$ di $g^{3+\varepsilon}_n$ presenta $n-3-\varepsilon$ condizioni alle $\varphi_{n-3}$. Dunque, dimensione effettiva delle $\varphi_{n-3}$ per $G_n$, o dimensione effettiva delle $\varphi_{n-4}$:
%-------------------------------pag. 93 -----------------------------------%
\[
\rho_e=p-1-(n-3-\varepsilon)=p-n+2+\varepsilon
\]
Dimensione virtuale (punti doppi $d=\frac{(n-1)(n-2)}{2}-p$), se le condizioni imposte dai punti base fossero indipendenti:
\[
\rho_v=\frac{(n-4)(n-1)}{2}-d=p-n+1
\]
Cioè $\rho_e$ supera $\rho_v$ almeno di $1$; cioè i punti doppi per la $F$ presentano condizioni non tutte indipendenti alle $\varphi_{n-4}$.
\end{enumerate}

\section{Formule di Pl\"{u}cker-Cayley per le curve sghembe.}
\label{sec:25}
Caratteri di Cayley delle curve sghembe:

\begin{description}
\item[$n$ ordine] numero dei punti della curva in un piano generico.
\item[$\nu$ classe] numero dei piani osculatori in un punto generico.
\item[$r=\rho$ rango] numero delle tangenti che incontrano una retta generica, cioè ordine della sviluppabile formata dalle tangenti alla sghemba.
\item[$d$ punti doppi apparenti] coppie di punti allineati con un punto fisso.
%dubbio: incerta l'interpretazione dell'abbreviazione appar.%
\item[$\delta$] coppie di piani osculatori secantisi su un piano fisso.
\item[$t$ ~$\pi$ tangenti doppi] coppie di tangenti da un punto
\item[$\tau$] coppie di tangenti che si segano su un piano.
\item[$K$ punti stazionari,] cuspidi per
%dubbio: incerta l'interpretazione di p. forse per%
curve reali, in cui s'incontrano $4$ piani osculatori successivi.
\item[$\chi$ ~piani osculatori stazionari,] che hanno un contatto quadripunto anzich\'{e} tripunto. 
\end{description}
%-------------------------------pag. 94 -----------------------------------%
Proiettando su un piano questa sghemba, la proiezione ha i seguenti caratteri di Pl\"{u}cker:
\[
\begin{matrix}
n & d & K\\
r & t & \nu\\
\end{matrix}
\]
Quindi possiamo scrivere le formole di Pl\"{u}cker:
\[
\left.
\begin{matrix}
r=n(n-1)-2d-3K\\
n=r(r-1)-2t-3\nu\\
\nu=3n(n-2)-6d-8K
\end{matrix}
\quad
\right\vert
\begin{matrix}
\quad r=\nu(\nu-1)-2\delta-3\chi \quad\\
\quad \nu=r(r-1)-2\tau-3n \quad\\
\quad n=3\nu(\nu-2)-6\delta-8\chi \quad \\
\end{matrix}
\left\vert
\quad
\begin{matrix}
p=\frac{(n-1)(n-2)}{2}-d-K\\
p=\frac{(\nu-1)(\nu-2)}{2}-\delta-\chi\\
\end{matrix}
\right.
\]
\noindent
Non si possono esprimere tutti i caratteri pl\"{u}ckeriani in funzione d'un solo, $n$; già da $n=4$ vi sono più tipi di sghembe di ugual ordine, ma genere diverso, ecc.\\
Supponendo $K=0$ si possono esprimere tutti i caratteri in funzione di due di essi, es. $n$ e $p$:
\begin{align*}
d&=\frac{(n-1)(n-2)}{2}-p\\
r&=2(n+p-1)\\
\nu&=3(n+2p-2)\\
\chi&=\dots\dots, \text{ecc.}
\end{align*}

\section{Caratteri della intersezione completa di due superficie, $f_\mu, \varphi_\nu$. - Le superficie d'ordine $\mu+\nu-4$ segano sulla curva la serie canonica completa; teorema sulla rappresentazione analitica di una superficie che passi per l'intersezione di due altre.}
\label{sec:26}

\subsection{Caratteri della intersezione $\mathcal{C}$ di $f_\mu$ e $\varphi_\nu$.}
\begin{equation}
\label{eq:int1}
n=\mu\nu
\end{equation}
%-------------------------------pag. 95 -----------------------------------%
Cerchiamo $r$: tangenti di $\mathcal{C}$ in un suo punto generico $(x~y~z~t)$:
\[
\left.
\begin{matrix}
X \frac{\partial f}{\partial x}+Y \frac{\partial f}{\partial y}+ Z \frac{\partial f}{\partial z} + T \frac{\partial f}{\partial t}=0\\
X \frac{\partial \varphi}{\partial x}+Y \frac{\partial \varphi}{\partial y}+ Z \frac{\partial \varphi}{\partial z} + T \frac{\partial \varphi}{\partial t}=0\\
\end{matrix}
\right\}
\quad
\text{Retta generica:}
\quad
\left\{
\begin{matrix}
aX+bY+cZ+dT=0\\
a'X+b'Y+c'Z+d'T=0\\
\end{matrix}
\right.
\]
Condizione perch\'{e} le quattro equazioni abbiano una soluzione $(X~Y~Z~T)$ comune:
\[
\begin{vmatrix}
\frac{\partial f}{\partial x} & \frac{\partial f}{\partial y} & \frac{\partial f}{\partial z} & \frac{\partial f}{\partial t}\\
\frac{\partial \varphi}{\partial x} & \frac{\partial \varphi}{\partial y} & \frac{\partial \varphi}{\partial z} & \frac{\partial \varphi}{\partial t}\\
a & b& c & d\\
a' & b' & c' & d'\\
\end{vmatrix}=0
\]
Equazione di grado $\mu+\nu -2$ in (x~y~z~t).\\
Soluzioni comune a questa ed alla $\mathcal{C}_{\mu\nu}$ sono:
\begin{equation}
\label{int:2}
r=\mu\nu(\mu+\nu-2)
\end{equation}
Se $K=0$, dal \S~precedente:
\begin{gather}
\label{eq:int3}
2p-2=\mu\nu(\mu+\nu-4),~ \text{ o }~ p=\frac{\mu\nu}{2}(\mu+\nu-4)+1\\
d=\frac{\mu\nu}{2}(\mu-1)(\nu-1) \quad\text{ ecc.}
\end{gather}
Dalla \eqref{eq:int3} le superficie d'ordine $\mu+\nu-4$, segano su $\mathcal{C}$ una $g_{2p-2}$.
\begin{teo} 
Se una superficie $\psi$ passa per la completa intersezione semplice di $f_\mu$ e $\varphi_\nu$, si può porre $\psi=Af+B\varphi$, coi tre termini di ugual grado.
\end{teo}
\noindent
Segando $f, \varphi, \psi$ col piano $z=0$, si hanno curve piane e si può applicare il Teorema di  N\"{o}ether:
\[
\psi(x ~y~0)=\alpha(x~y)f(x~y~0)+\beta(x~y)\varphi(x~y~0)
\]
Scrivo:
\begin{equation}
\label{eq:int4}
\Phi(x~y~z)=\psi(x~y~z)-\left[ \alpha f(x~y~z)+ \beta\varphi(x~y~z) \right]
\end{equation}
polinomio identicamente nullo per $z=0$, quindi divisibile per $z$:
\begin{equation}
\label{eq:int5}
\Phi(x~y~z)=z \psi'(x~y~z)
\end{equation}
$\psi'$ passa per le intersezioni di $f$ e $\varphi$ (perch\'{e} ogni punto di
%-------------------------------pag. 96 -----------------------------------%
questa curva annulla $f, \varphi, \psi$), ed è di ordine più basso di uno di $\Phi$. Ora il nostro Teorema sia dimostrato per la superficie d'ordine $(n-1)$: $\psi'(x~y~z)=A'f+B'\varphi$.
Allora sostituendo nella \eqref{eq:int5} e poi nella \eqref{eq:int4}
\[
\psi=(\alpha+zA')f+(\beta+zB')\varphi
\]
Cioè è dimostrato il Teorema per le superficie d'ordine $n$.\\
Resta a dimostrare il Teorema vero per le superficie $\psi_\mu$ d'ordine abbastanza basso, che debbono passare per la intersezione di $f_\mu$ e $\varphi_\nu$.
\begin{enumerate}[label=\emph{\arabic* $^o$)}]
\item Sia $\mu<\nu$; $f$ e $\varphi$ s'intersecano in una $\mathcal{C}_{\mu\nu}$. Devo scrivere una nuova superficie d'ordine $\mu$ per la curva, ma vi passa solo la $f$ se no avrebbe ordine $\mu^2$; sicch\'{e} $\psi$ deve coincidere con $f$: $\psi=\alpha f+0\varphi$ ($\alpha$ costante).
%dubbio: la resa del cost scritto sotto come pedice non so se è ottimale (provato anche con \atop ma sposta alpha in alto%
\item Sia $\mu=\nu$; allora abbiamo $f_\mu$ e $\varphi_\mu$, e $\psi_{\mu^2}$ per la loro intersezione. L'intersezione appartiene al fascio $f=0$ ~$\varphi=0$ (si vede analogamente al procedimento noto in geometria piana);
%dubbio: indicato geom. \pi%
allora ($\alpha$ e $\beta$ costanti):
\[
\psi=\alpha f+\beta \varphi
\]
\end{enumerate}

\begin{teo} [Teorema di N\"{o}ther]
Sull'intersezione completa di $f_\mu$ e $\varphi_\nu$, le superficie d'ordine $\mu+\nu-4$ segano la serie canonica completa.
\end{teo}
\noindent
(Serie canonica su una curva sghemba è la serie che corrisponde alla serie canonica di una curva piana in corrispondenza birazionale con la sghemba). La serie canonica è l'unica di caratteri $g^{p-1}_{2p-2}$. Si può vedere così: sia $F$ una trasformata birazionale piana della $\mathcal{C}$. Punti di tangenza da $P$ (intersezioni colla prima polare) sono in numero: 
%-------------------------------pag. 97 -----------------------------------%
$n(n-1)-2d$.
\begin{figure}[!htbp]
\centering
\includegraphics[scale=0.5]{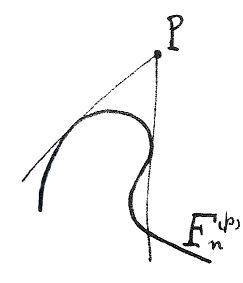}
%\caption{}
\end{figure}
La prima polare è una particolare aggiunta $\varphi_{n-1}$; tra queste aggiunte $\varphi_{n-1}$ ve ne sono che si spezzano in $\varphi_{n-3}$
%dubbio: s o 3?%
 e due rette; il gruppo dei punti di contatto $G_r$ appartiene alla serie rigata dalle $\varphi_{n-1}$; considerando un'altra $\varphi_{n-1}$ che si spezzi e che dà un altro gruppo,
\[
G_r\equiv G_{2p-2}+2G_n.
\]
Può enunciarsi: data una $g^1_n$ su una curva, il gruppo dei punti della $g^1_n$ (tangenze), $G_r$, $\equiv$ 
%dubbio: forse meglio sciogliere il simbolo in equivale%
un gruppo canonico più due gruppi della $g^1_n$ data. Sotto questa forma è una proprietà della Geometria sulle curve, quindi si può trasportare alle curve sghembe: $G_r\equiv 1 \text{ gruppo canonico} + 2 \text{ gruppi segati dai nostri piani (corrispondenti alla rette):}$
\begin{equation}
\label{eq:no1}
G_r\equiv G_{2p-2}+2G_n
\end{equation}
Ora $G_r$ è segato sulla $\mathcal{C}$ da una $\psi_{\mu+\nu-2}$, cioè equivale ad ogni gruppo segato da una superficie d'ordine $\mu+\nu-2$. Tra queste superficie posso pensarne una che si spezzi in una $\psi_{\mu+\nu-4}$, $\Gamma_{2p-2}$, più due gruppi segati da piani
\begin{equation}
\label{eq:no2}
G_r\equiv \Gamma_{2p-2}+2G_n;
\end{equation}
quindi 
\[
G_{2p-2}+\cancel{2G_n} \equiv \Gamma_{2p-2}+\cancel{2G_n}
\]
Cioè $1 \text{ gruppo canonico}\equiv $ un gruppo d'una $\psi_{\mu\nu-4}$. Cioè le $\psi_{\mu+\nu-4}$ segano gruppi canonici, cioè la serie canonica.

Vediamo la dimensione della $g_{2p-2}$ segata dalle
%-------------------------------pag. 98 -----------------------------------%
superficie d'ordine $\mu+ \nu-4$. Queste dipendono da parametri:
\[
\rho=\frac{(\mu+\nu-3)(\mu+\nu-2)(\mu+\nu-1)}{6}-1
\]
Questa sarebbe  la dimensione della serie, senonch\'{e} vi sono delle superficie che contengono come parte la $\mathcal{C}_{\mu\nu}$, completa intersezione di $f$ e $\varphi$, cioè che passano per $\mathcal{C}_{\mu\nu}$. Per il Teorema precedente esse si possono scrivere:
\[
\psi_{\mu+\nu-4}=A_{\nu-4}f_\mu+B_{\mu-4}\varphi_\nu
\]
[Qui non si presenta l'incertezza della Geometria piana, se fosse $Af+B\varphi\equiv A'f+B'\varphi$, conseguirebbe $(A-A')f\equiv (B'-B)\varphi$; $\varphi$ dovrebbe dividere $A-A'$, ma guardando i gradi non è possibile; sicch\'{e} non dobbiamo preoccuparci della possibile multipla rappresentazione]
%dubbio: forse meglio mettere una nota a piè di pagina?%
e la loro dimensione è il numero ($\sigma$) dei parametri di questa espressione:
\[
\sigma=\frac{(\nu-3)(\nu-2)(\nu-1)}{6}+\frac{(\mu-3)(\mu -2)(\mu-1)}{6}-1
\]
Quindi dimensione della serie $\rho - \sigma-1=p-1$.
%dubbio: 6 o sigma?%
Allora la serie è $g^{p-1}_{2p-2}$ canonica completa.

\section{Relazioni tra i caratteri di due curve costituenti insieme la intersezione completa di due superficie; come la serie canonica completa possa segarsi sopra una delle due curve.}
\label{sec:27}
Siano le due curve $\mathcal{C}^{p_1\atop r_1}_{n_1}$  $\mathcal{C}^{p_2\atop r_2}_{n_2}$ costituenti la intersezione di $f_\mu$ e $\varphi_\nu$.
\begin{equation}
\label{eq:cur1}
n_1+n_2=\mu\nu
\end{equation}
$\mathcal{C}_1$ e $\mathcal{C}_2$ hanno certamente \underline{$i$}$>0$ punti comuni (non lo dimostriamo).
%-------------------------------pag. 99 -----------------------------------%
Come al \S ~precedente, cerchiamo $r_1$; i punti di $\mathcal{C}_1$ da cui partono tangenti che incontrano una retta assegnata, soddisfano l'equazione:
\[
\begin{vmatrix}
\frac{\partial f}{\partial x} & \frac{\partial f}{\partial y} & \frac{\partial f}{\partial z} & \frac{\partial f}{\partial t}\\
\frac{\partial \varphi}{\partial x} & \frac{\partial \varphi}{\partial y} & \frac{\partial \varphi}{\partial z} & \frac{\partial \varphi}{\partial t}\\
a & b & c & d\\
a' & b' & c' & d'\\
\end{vmatrix}=0
\]
\begin{wrapfigure}{l}{0pt}
\includegraphics[scale=0.5]{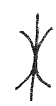}
%\caption{}
\end{wrapfigure}
%dubbio: vale la pena inserire questa figura?%

Quindi sarebbe $r_1=n_1(\mu+\nu-2)$; ma analogamente a considerazione piana, negli \underline{$i$} punti si hanno contatti di $f$ e $\varphi$, quindi coincidono i piani tangenti, ed essi \underline{$i$} punti sono anche comuni a $\mathcal{C}$, ed al determinante sopra; quindi 
\begin{subequations}
\label{eq:cur2}
\begin{align}
\label{eq:cur2a} r_1&=n_1(\mu+\nu-2)-i\\
\label{eq:cur2b} r_2&=n_2(\mu+\nu-2)-i
\end{align}
\end{subequations}
Quindi, ricordando
\begin{subequations}
\label{eq:cur3}
\begin{align}
\label{eq:cur3a} 2p_1-2=n_1(\mu+\nu-4)-i; \\
\label{eq:cur3b} 2p_2 -2=n_2(\mu+\nu-4)-i.
\end{align}
\end{subequations}
Così bastano $n_1$ e $p_1$ per determinare i caratteri di $\mathcal{C}_2$.
\[
p_1+p_2=\frac{\mu\nu}{2}(\mu+\nu-4)-i+2
\]
Il genere di $\mathcal{C}_1\mathcal{C}_2$ supposta momentaneamente irriducibile sarebbe $p=\frac{\mu\nu}{2}(mu+\nu-4)+1$, quindi
\begin{equation}
\label{eq:cur4}
p=p_1+p_2+i-1
\end{equation}
Ora $f$ sia fissa e $\varphi$ variabile; varierà $\mathcal{C}_1\mathcal{C}_2$, e per valori particolari dei parametri si spezzerà; se una curva generalmente irriducibile variabile in un sistema continuo per valori determinati dai parametri si spezza, per Genere della spezzata, definito dalla \eqref{eq:cur4}, è il genere della irriducibile infinitamente vicina.\\
Per lo spezzamento in più parti $p=$ somma generi  $+$ somma punti
%-------------------------------pag. 100 -----------------------------------%
 d'incontro 2 a 2 delle componenti $-$ differenza tra numero componenti ed unità.
\begin{teo} 
Su $\mathcal{C}^{n_1p_1}_{1}$ parziale intersezione di $f_\mu$ e $\varphi_\nu$, la serie canonica completa è segata dalle $\psi_{\mu+\nu-4}$ che passano per $\mathcal{C}^{n_2p_2}_{2}$.
\end{teo}
Intanto le $\psi$ per $\mathcal{C}_2$ passando per gli \underline{$i$} punti comuni a $\mathcal{C}_1$ e $\mathcal{C}_2$, segano su $\mathcal{C}_1$ una $g_{2p_1-2}$.

\subsubsection{Dimensione delle $\psi$ per $\mathcal{C}_2$.}
Dimensione delle $\psi$:
\[
\frac{(\mu+\nu-3)(\mu+\nu-1)\dots}{6}-1.
\]
Le $\psi$ (libere) su $\mathcal{C}_2$ determinano $g^{p_2-2+i-\delta}_{2p_2-2+i}$ (non speciale, $i>1$, $\delta$ se la serie non è completa). Pel passaggio delle $\psi$ per $\mathcal{C}_2$ occorrono $p_2-1+i-\delta$ (\emph{dimensione +1}) condizioni. Quindi, dimensione delle $\psi$ per $\mathcal{C}_2$:
\[
\rho=\frac{(\mu+\nu-3)\dots}{6}-(p_2-1+i-\delta)
\] 
Questa sarebbe la dimensione della $g_{2p_1-2}$ segata dalle $\psi$ per $\mathcal{C}_2$ su $\mathcal{C}_1$ se nessuna di queste $\psi$ contenesse $\mathcal{C}_1$;

\subsubsection{Dimensione delle $\psi$ per $\mathcal{C}_2\times \mathcal{C}_1$.}
La $g_{2p_1-2}$ detta possiamo indicarla 
\[
g^{p_1-1-\varepsilon}_{2p_1-\varepsilon}
\]
perch\'{e} sia la serie canonica completa, deve essere $\varepsilon=0$.\\
Perch\'{e} queste $\psi$ per $\mathcal{C}_2$ dipendenti da $\rho$ parametri passino per $\mathcal{C}_1$ occorrono $p_1-\epsilon$ condizioni, quindi la dimensione delle $\psi$ per $\mathcal{C}_1\times \mathcal{C}_2$ è:
\[
\tau=\rho-(p_1-\varepsilon)=\frac{(\mu+\nu-3)\dots}{6}-1-(p_1+p_2+i-1-\delta-\varepsilon)=\frac{(\mu+\nu-3)\dots}{6}-1-p+\delta+\epsilon.
\]
Vediamo $\tau$ altrimenti:\\
Le $\psi$ dello spazio segano su $\mathcal{C}_1\times \mathcal{C}_2$ (intersezione completa) la $g^{p-1}_{2p-2}$;
%-------------------------------pag. 101 -----------------------------------%
perch\'{e} una di esse contenga $\mathcal{C}_1\times \mathcal{C}_2$ occorrono $p$ condizioni; sicch\'{e} il sistema delle $\psi$ per $\mathcal{C}_1\times \mathcal{C}_2$ ha la dimensione 
\[
\tau=\frac{(\mu+\nu-3) \dots}{6}-1-p
\]
Quindi si vede $$\delta=\epsilon=0,$$ e le $\psi$ per $\mathcal{C}_2$ segano su $f$ la $g^{p-1}_{2p-2}$ canonica completa.

\section{Postulazione di una sghemba $\mathcal{C}^{n,p}$ rispetto alle superficie $f_m$ - Serie segata dalle superficie d'ordine $\mu+\nu-s$ sulla intersezione di $f_\mu$ e $\varphi_\nu$.}
\label{sec:28}

Postulazione di una sghemba $\mathcal{C}^{n,p}$ rispetto alle superficie d'ordine $m$, $f_m$, è il numero delle condizioni perch\'{e} $f_m$ passi per $\mathcal{C}^{n,p}$. Tutte le $f_m$ dello spazio segano su $\mathcal{C}$ una $g^x_{mn}$; ci basta trovare $x$.
Se la $g_{mn}$ fosse completa non speciale, allora $x=mn-p$, ed \underline{$mn-p+1$} valore normale della Postulazione. Ma la $g_{mn}$ potrebbe essere speciale (completa), allora si avrebbe $g^{mn-p+s}_{mn}$, e se la serie fosse ancora incompleta si avrebbe $g^{mn-p+s-\delta}_{mn}$.\\
Si vede che appena $mn>2p-2$ la serie diviene non speciale, $s=0$; quindi \emph{per $m$ abbastanza alto la dimensione è $\leq mn-p$, e la postulazione $\leq mn-p+1$}. Cioè il sistema delle $f_m$ per $\mathcal{C}$ è regolare o sovrabbondante.
\begin{teo} 
Per $m$ sufficientemente alto anche $\delta=0$, o meglio: se $\mathcal{C}_{\mu\nu}$ 
%-------------------------------pag. 102 -----------------------------------%
è intersezione di $f_\mu$ e $\varphi_\nu$, le superficie d'ordine $m\geq \mu+\nu-3$ vi segano una serie completa non speciale.
\end{teo}
%dubbio: grassetto per rendere la sottolineatura%
\begin{description}
\item[$1^o$ Caso - $\mathcal{C}_n$ si intersezione totale di $f_\mu$ e $\varphi_\nu$.]
$n=\mu\nu$ ~$p-1=\frac{\mu\nu}{2}(\mu+\nu-4)$.
Postulazione di $\mathcal{C}=$ Dimensione delle $\psi_m - $Dimensione delle $\psi$  per $\mathcal{C}$.
Per $\mathcal{C}$:
\[
\psi_m=A_{m-\mu}f_\mu+B_{m-\nu}\varphi_{\nu}
\] 
Per avere la dimensione delle $\psi$ per $\mathcal{C}$ sembrerebbe bastasse contare i coefficienti di $A$ e $B$; ma per $m$ abbastanza alto si troverebbe così infinite volte la stessa $\psi$. Si trova che in realtà la dimensione della $\psi$ per $\mathcal{C}$ è:
\[
\frac{(m+1)(m+2)(m+3)}{\sigma}-1-(mn-p+1)
\]
donde la postulazione ha il valore normale $(mn-p+1)$.
\item[$2^o$ Caso - $\mathcal{C}'_n$ sia intersezione parziale di $f_\mu$ e $\varphi_\nu$.]
Le $\psi_{m\geq \mu+\nu-3}$ $$mn'>(\mu+\nu-4)n_1>2p_1-2$$ segano su $\mathcal{C}'$ una $g_{mn_1}$ certamente non speciale 
%dubbio: la quadra non si chiude...%
\[
g^{mn_1-p_1-\delta}_{mn_1} \qquad \delta\geq0 \qquad [\delta_1 \text{ non sappiamo se è completa o no}]
\] 
Cioè $\mathcal{C}'$ impone alle $\psi$ per essa $mn_1-p_1-\delta+1$ condizioni.\\
Le $	\psi$ per $\mathcal{C}'$ segano su $\mathcal{C}^{\prime\prime}$ una $g_{mn_2-i}$, anche certo non speciale, ma non sappiamo se completa o no ($\varepsilon$)
\[
g^{mn_2-i-p_2-\varepsilon}_{mn_2-i}
\]
Numero totale delle condizioni imposte alle $\psi$ per passare per $\mathcal{C}_1\mathcal{C}_2$:
\begin{equation}
\label{eq:post}
mn-\underline{p}+1-(\delta+\varepsilon)
\end{equation}
($n$ e $p$ ordine e genere di $\mathcal{C}'\times \mathcal{C}^{\prime\prime}$). Pel caso precedente stesso numero è
\[
mn-p+1
\]
%-------------------------------pag. 103 -----------------------------------%
donde 
\[
\delta=\epsilon=0
\]
Cioè la Postulazione ha il valore normale $mn-p+1$.\\
 Cioè la serie segata su $\mathcal{C}'$ dalle $\psi_{m\geq\mu+\nu-3}$ è completa non speciale.
\end{description}
%dubbio: qui c'è un asterisco che non capisco a dove si riferisca%
La rappresentazione si può fare in infiniti modi; forma più generale: 
\begin{equation}
\label{eq:xi}
\psi_m\equiv (A_{m-\mu}+H_{m-\mu-\nu}\varphi_\nu)f_\mu+(B_{m-\nu}-H_{m-\mu-\nu}f_\mu)\varphi_\nu
\end{equation}
Quanti i parametri arbitrari in $\psi$?\\
Posso profittare dei coefficienti di $H$ per esempio per ridurre a zero $m-\mu-\nu$ coefficienti del moltepl.
%dubbio: incerta l'interpretazione dell'abbreviazione moltepl. %
di $f_\mu$ in \eqref{eq:xi}. Allora:\\
numero dei coefficienti di $\psi \equiv$ Coefficienti di $A$ + Coefficienti di $B$ - coefficienti $H$-1.
Risulta

\section{Serie segate su una sghemba $\mathcal{C}_n$ dalle $\psi_{m\geq n-2}$ e dalle $\psi_{m\geq\frac{n-3}{2}}$.}
\label{sec:29}
\begin{teo} 
Su una sghemba $\mathcal{C}_n$, le superficie d'ordine $\geq n-2$ segano serie complete non speciali.
\end{teo}

\begin{figure}[!htbp]
\centering
\includegraphics[scale=0.5]{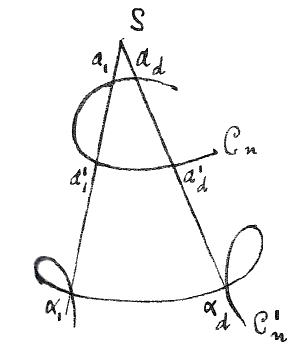}
%\caption{}
\end{figure}
\noindent
Si proietti da $S$; la $\mathcal{C}'_n$ avrà \emph{certamente $d$} punti doppi. Su $\mathcal{C}'$ considero le aggiunte (piane) d'ordine $m\geq n-2$; vi segano una serie completa non speciale $g^{mn-2d-p}_{mn-2d}$; cioè i coni corrispondenti, d'ordine $mn\geq n-2$ segano su 
 %-------------------------------pag. 104 -----------------------------------%
$\mathcal{C}_n$ una $g^{mn-2d-p}_{mn-2d}$ completa non speciale.\\
La stessa può esser segata da tutte le superficie d'ordine $m\geq n-2$ per $a_1a'_1, \dots a_da'_d$ (questa contiene la prima completa, quindi è completa e coincide con essa). Cioè $\mathcal{C}_n$ impone $mn-2d-p+1$ condizioni alle $\Phi_{n-2+k}$ per essa (e pei $2d$ punti doppi apparenti).\\
Liberando le $\Phi$ dalle condizioni di passaggio pei $2d$ punti doppi, \emph{se queste fossero} \emph{ indipendenti}, la $\mathcal{C}_n$ imporrebbe $mn-p+1$ condizioni.\\
Sono indipendenti: basta vedere che possiamo condurre una $\Phi$ per $2d-1$ qualsiasi ($a_1a'_1 - a_d$) e non pel rimanente ($a'_d$).\\
%dubbio: forse (a_1a'_1 \dots a_d)%
Sul piano di $\mathcal{C}'$ costruisco una curva d'ordine $m-1\geq n-3$, che sia aggiunta a $\mathcal{C}'$ in tutti i punti doppi salvo $\alpha_d$ (Per una aggiunta d'ordine $\geq n-3$ i punti doppi d'una curva piana impongono condizioni indipendenti), e proiettiamola da $S$; otteniamo una superficie d'ordine $\underbrace{m-1}$ per $a_1a'_1-a_{d-1}a'_{d-1}$; vi si aggiunga un piano qualsiasi per $a_d$ e non per $a'_d$; così si ha una $\Phi_m$ per $2d-1$ punti e non pel rimanente.
\begin{lem}
Se la serie speciale $g^r_m$ è contenuta in un'altra speciale $g^R_M$, ogni gruppo della prima impone al più $m-r$ condizioni ad un gruppo della seconda costretto a contenerlo, cioè ad una curva della seconda serie.
\end{lem}
\begin{teo} 
\label{teo:sghemba}
Su una sghemba $\mathcal{C}_n$ le $\psi_{m\geq\frac{n-3}{2}}$ segano certo serie non speciali.
\end{teo}
\noindent
Su $\mathcal{C}_n$, se la $g^3_n$ segata dai piani dello spazio è non speciale, sarà non speciale
 %-------------------------------pag. 105 -----------------------------------%
ogni $g_{mn}$ segata da superficie $\psi_m$ (es. quadriche, ecc.) ch\'{e} se la $g_{mn}$ fosse speciale, sarebbe speciale la $g^3_n$ contenuta in essa. Ma se la $g^3_n$ è speciale, che valore occorre per $m$ perch\'{e} le $\psi_m$ seghino una serie non speciale?\\
Se $\psi_m$ sega una $g_{mn}$ speciale, contiene la $g^3_n$ dei piani dello spazio, pel lemma: ogni $\psi_m$ per $n-3$ punti di un gruppo piano, passa pei rimanenti $3$. Perch\'{e} ciò non sia la $g_{mn}$ deve essere non speciale (se fosse speciale \dots). Una sezione piana di $\mathcal{C}$, vi sega $n$ punti, mai tre allineati (non lo dimostriamo); per questi, consideriamo $n$ piani, pei punti (tra gli $n$) $1,2,3,4,\dots, 2m-1, 2m$. Questi costituiscono una $\psi_m$ per $2m$ punti d'una sezione piana, che non passa pei rimanenti $n-2m$; voglio che $\psi_m$ contenga tutti i punti della sezione piana salvo al più tre, dunque deve essere 
\[
n-2m\leq 3,\quad \text{ cioè } \quad m\geq \frac{n-3}{2}
\]
Così la $g_{mn}$ è non speciale, ch\'{e} se fosse speciale la $\psi$ per $n-3$ punti della sezione piana passerebbero pei rimanenti tre.

\section{Massimo genere d'una sghemba $\mathcal{C}_n$ priva di punti multipli.}
\label{sec:30}
\subsubsection{Numero condizioni ($\nu_K$) da imporre ad una $\psi_K$ perch\'{e} contenga un gruppo $\pi$ $G_n$ di $\mathcal{C}_n$.}
Conduciamo un piano generico per $\underline{1 2}$, uno per \underline{3 4}, \dots uno
 %-------------------------------pag. 106 -----------------------------------%
per $2K-1,2K$; ho una $\psi_K$ che passa per $2K$ punti e non pei rimanenti (mai tre allineati) $n-2K$; dunque $2K$ condizioni è troppo poco; \emph{almeno}
\[
\nu_K\geq 2K+1
\]
ma deve anche essere $n>\nu_K\geq 2K+1$, donde massimo valore di $K$
\begin{equation}
\label{eq:max1}
\chi=\frac{n-2}{2} \text{ per $n$ pari} \qquad \text{o} \qquad \chi=\frac{n-3}{2} \text{ per $n$ dispari.}
\end{equation}
\begin{equation}
\label{eq:I}
\nu_K\geq 2K+1 \qquad K=1,2,\dots \chi
\end{equation}

%dubbio: la successiva dovrebbe essere per rimanere coerenti una subsubsection%
\noindent
\emph{Le $\psi_K$ segano su $\mathcal{C}_n$ una $g^{r_K}_{nK}$; cerchiamo $r_K$}; imposte le condizioni $\nu_K$ di passaggio per $n$ punti piani, resta una $g^{r_K-\nu_K}_{(K-1)n}$. Tra queste $\psi_K$ per $G_n$ ve ne sono che si spezzano nel piano di $G_n$ e una $\Phi_{K-1}$; così la $g^{r_K-\nu_K}_{(K-1)n}$ contiene o coincide colla $g^{r_{K-1}}_{(K-1)n}$ segata delle superficie d'ordine $K-1$. Quindi $r_K-\nu_K\geq r_{K-1}$ cioè $r_K-r_{K-1}\geq 2K+1$.
Dando a $K$ successivi valori, posto $r_0=0$, si ha: $r_1-r_0\geq3, r_2-r_1\geq5, \dots r_\chi-r_{\chi-1}\geq 2\chi+1$; sommando:
\[
r_\chi\geq \chi(\chi+2)
\]
Pel Teorema alla pagina precedente \ref{teo:sghemba}
%dubbio: inserito il riferimento al teorema%
 \underline{$g^{r_\chi}_{\chi n}$} è non speciale, dunque 
\[
\chi n-r_\chi=p \quad \Rightarrow \quad  \chi_n-\chi(\chi-2)\geq p
\]
e sostituendo le \eqref{eq:max1}:
\begin{equation}
\label{eq:max2}
\left. 
\begin{aligned} 
p \leq  (\frac{n-1}{2})^2 \quad  \text{ per $n$ pari }\\
p\leq \frac{(n-1)(n-3)}{4} \quad\text{ per $n$ dispari }
\end{aligned} 
\right\} 
\end{equation}
 %-------------------------------pag. 107 -----------------------------------%
\begin{teo} 
\label{teo:quad}
Le curve di massimo $p$ rispetto ad un dato ordine, appartengono a quadriche. 
\end{teo}
\noindent
Le quadriche dello spazio segano su $\mathcal{C}_n$ una $g^{r_2}_{2n}$, dove $r_2$ è $9$ se nessuna quadrica contiene $\mathcal{C}_n$, altrimenti $g^8_{2n}$.\\
Essendo raggiunto il massimo genere, sopra vale il segno uguale: $r_1=3, \, r_2-r_1=5$
\[
\underline{r_2=8}
\]
cioè $\mathcal{C}$ sta su una quadrica.
\begin{teo} 
Idem come sopra,
%dubbio: come sopra?%
e segano le rette di entrambi i sistemi in $\frac{n}{2}$ punti ($n$ pari) o di uno in $\frac{n-1}{2}$ e l'altro $\frac{n+1}{2}$ ($n$ dispari).
\end{teo}
\noindent
Seghino le rette di un sistema in $\mu$ punti , dell'altro in $n-\mu$. Si proietti $\mathcal{C}$ da un punto della quadrica su un piano; la proiezione avrà un punto multiplo d'ordine $\mu$, ed uno $(n-\mu)$, e se $\mathcal{C}$ non ha altri punti multipli, la proiezione neanche ne avrà altri. Allora siamo indotti a valutare il genere d'una curva piana; 
\[
p=\frac{(n-1)(n-2)}{2}-\frac{\mu(\mu-1)}{2}-\frac{(n-\mu)(n-\mu-1)}{2}=(\mu-1)(n-\mu-1)
\]
somma di fattori costante è massima per fattori uguali.
Allora si vede:
\begin{gather*}
p \text{ massimo per } \mu=\frac{n}{2} \quad \text{($n$ pari)}\\
p \text{ massimo per } \mu=\frac{n-1}{2}, \quad n-\mu=\frac{n+1}{2} \quad \text{($n$ dispari)}
\end{gather*}

%dubbio: più che una subsection dovrei inserirlo nell'ambiente teorema%
\subsubsection{Costruzione.}
\emph{Le curve $\mathcal{C}_n$ di $p$ massimo sono intersezioni d'una quadrica con una $f_{\frac{n}{2}}$ ($n$ pari) o $f_{\frac{n+1}{2}}$ ($n$ dispari)}.\\

\begin{description}
\item[$n$ pari] $n=2\mu \quad \mathcal{C}_{2\mu} \quad p=(\mu-1)^2$\\
\emph{$\mathcal{C}_{2\mu}$ sta su una quadrica}, pel Teorema antiprecedente \ref{teo:quad}.
%dubbio: inserire il riferimento al teorema? %
Le $f_\mu$ segano
 %-------------------------------pag. 108 -----------------------------------%
su $\mathcal{C}$ una $g^{2\mu^2-(\mu-1)^2-\delta}_{2\mu^2}
=g^{\mu^2+2\mu-1-\delta}_{2\mu^2}$ non speciale ($\delta \geq 0$, ch\'{e} \dots). Perch\'{e} una $f_\mu$ contenga la $\mathcal{C}$, basta imporre $\mu(\mu+2)-\delta$ condizioni. Ma $f_\mu$ potrebbe spezzarsi \emph{sempre} nella quadrica ed una $f'_{\mu-2}$? Una $f_\mu$ per $\mathcal{C}$ dipende da 
\[
\frac{(\mu+1)(\mu+2)(\mu+3)}{\sigma}-1-[\mu(\mu+2)-\delta]=\frac{\mu(\mu+1)(\mu-1)}{\sigma}+\delta
\]
parametri; cioè dimensione delle $f_\mu$ per $\mathcal{C}$ è almeno $\frac{\mu(\mu-1)(\mu+1)}{\sigma}$. Se tutte le $f$ si spezzassero nella quadrica ed una $f'_{\mu-2}$, la dimensione del sistema sarebbe quella delle $f'_{\mu-2}$, cioè $\frac{\mu(\mu-1)(\mu+1)}{\sigma}-1$, più bassa di quella che abbiamo trovato. Cioè c'è almeno qualche $f_\mu$ per la $\mathcal{C}$ che non si spezza nella quadrica ed una $f'_{\mu-2}$. Ragionamento analogo per $n$ dispari.
 \end{description}

\begin{teo} [Teorema di Halphen-Valentiner]
Una curva appartenente ad $f_\mu$, di genere massimo compatibile coll'ordine  della curva e coll'ordine della superficie, o è intersezione totale di $f_\mu$ e $\varphi_\nu$; o parziale, e allora la residua è una curva piana.
\end{teo}

 %-------------------------------pag. 109 -----------------------------------%
\section{Numero dei parametri da cui dipendono le curve sghembe.}
\label{sec:31}
\subsubsection{Il numero dei parametri da cui dipendono le $\mathcal{C}_{n}$ razionali è $4n$. - ($p=0$).}
I piani dello spazio segano su $\mathcal{C}_{n(p=0)}$ una $g^{3}_{n}$. Sulla retta corrispondente (le curve razionali corrispondono birazionalmente ad una retta), la corrispondente $g^{3}_{n}$ sarà rappresentata da
\begin{equation}
\label{eq:par1}
\lambda_{0}f_{0}(\xi)+\dots+\lambda_{3}f_{3}(\xi)=0
\end{equation}
Viceversa dalla $g^{3}_{n}$ \eqref{eq:par1} sulla retta, a questa corrisponderà la curva:
\begin{equation}
\label{eq:par2}
x:y:z:t=f_{0}:f_{1}:f_{2}:f_{3}
\end{equation}
su cui la serie corrispondente è segata dai piani dello spazio.\\
A prima vista, data la serie \eqref{eq:par1}, la curva \eqref{eq:par2} pare completamente determinata; ma è determinata dando le coordinate proiettive; variandole (\emph{collineazione nello spazio}) varia la curva \eqref{eq:par2}. \emph{Quanti parametri introduce questa libertà di scelta tra una curva e la collineazione?} (Questo numero lo dovremo poi aggiungere ai parametri da cui dipende una $g^{3}_{n}$ su una retta). Una collineazione nello spazio è determinata quando a $5$ punti faccio corrispondere $5$ punti qualsiasi nel secondo spazio; cioè è determinata da $15$ parametri.\\
Quindi, data la $g^{3}_{n}$ \eqref{eq:par1} sulla retta, restano determinate in corrispondenza $\infty^{15}$ curve. \emph{Ma} supponiamo data la curva e
 %-------------------------------pag. 110 -----------------------------------%
messa in corrispondenza birazionale colla retta, posso sulla retta variare la rappresentazione, variare la proiettività senza che la curva cambi. Una proiettività sulla retta dipende da $3$ parametri, così una curva viene rappresentata $\infty^{3}$ volte; resta
\[
15-3=12.
\]
Allora ci basta calcolare $\rho$, \emph{numero dei parametri da cui dipende una $g^{3}_{n}$ su una retta}, quindi 
\[
u=\rho+12
\]
una $g^{3}_{n}\equiv\infty^{3}$ gruppi di $n$ punti $\equiv$ punto dello spazio ad $n$ dimensioni costretto nello spazio a $3$, che dipende da $n-3$ parametri. La retta in $S_{3}$ dipende da $4$ parametri ($\equiv 4$ punti) dunque
\begin{equation}
\label{eq:par3}
u=4(n-3)+12=4n
\end{equation}
Anche per le curve di genere $1$ vale lo stesso risultato, e anche, sotto larghe ipotesi, per curve di genere superiore.

\subsubsection{Curve di genere $1$.}
Una $\mathcal{C}^{p=1}_{n}$ si può trasformare in una $F_{3}$ piana; vi sono $\infty^{1}$ $F_{3}$ (cubiche) modello, con tutti i valori possibili di un certo birapporto (\emph{modulo}). Una retta ammette $\infty^{3}$  trasformazioni birazionali in s\'{e} stessa, mentre una cubica piana $\infty^{1}$ tali che un punto si può trasformare in qualunque altro. Così un punto di $\mathcal{C}_{n}$ si può mutare in un punto che si vuole di $F_{3}$.\\
 %-------------------------------pag. 111 -----------------------------------%
Data $\mathcal{C}^{p=1}_{n}$ i piani dello spazio vi segano una $g^{3}_{n}$; viceversa data una $g^{3}_{n}$ di $F_{3}$ esiste una trasformazione birazionale tale che ad essa corrisponde la $g^{3}_{n}$ segata su $\mathcal{C}$ dai piani dello spazio; come al caso $p=0$ la $\mathcal{C}$ è determinata a meno di una collineazione, cioè sono determinate $\infty^{15}$ $\mathcal{C}$. Ma supponiamo
%dubbio: nel manoscritto riportato supponendo poi corretto a matita in supponiamo%
data $\mathcal{C}$ e messa in corrispondenza birazionale con $F_{3}$.\\
Per avere il numero dei parametri, a $15$ aggiungere $1$ (la cubica piana in corrispondenza a $\mathcal{C}$ dipenda da $1$ valore del birapporto), poi togliere $1$ (vi sono $\infty^{1}$ trasformazioni birazionali della $F_{3}$ su s\'{e} stessa, e variandole non varia conseguentemente la $\mathcal{C}$), quindi aggiungere $\rho$ (\emph{numero dei parametri da cui dipende una $g^{3}_{n}$ sulla $F_{3}$}).
\[
u=\rho+15
\]
Resta a calcolare $\rho$. \emph{Una $g^{3}_{n}$ in una curva ellittica è sempre contenuta in una completa $g^{n-p}_{n}$}
%dubbio: confusione tra il simbolo u e il simbolo n%
(v. Premessa al Teorema di Riemann-Roch). Quante sono queste complete $g^{n-1}_{n}$? Cioè da quanti parametri dipendono? Essa possiede sempre un determinato gruppo che conti $n-1$ punti assegnati sulla curva, ne conterrà un ulteriore $x$; viceversa, se fisso ad arbitrario $x$ sulla curva, questo cogli $n-1$ dà $n$ punti che determinano la $g^{n-1}_{n}$ completa cui esso appartiene. Così sono tante le $g^{n-1}_{n}$ quanti i modi con cui posso prendere $x$ sulle curve,
 %-------------------------------pag. 112 -----------------------------------%
cioè $\infty^{1}$. E allora, \emph{le complete $g^{n-1}_{n}$ dipendono da un solo parametro} \S~ \ref{sec:13}: $\rho=1+\dots$ Ora quante sono le $g^{3}_{n}$ contenute entro le $g^{n-1}_{n}$? \`{E} il numero dei parametri da cui dipende un punto da cui dipende un punto di $S_{n-1}$ costretto in $S_{3}$, cioè dalla lezione precedente
%dubbio: inserire riferimento latex?%
$4(n-4)$. Allora $\rho=1+4(n-4)=4n-15$
\begin{equation}
\label{eq:par3a}
u=4n
\end{equation}

\subsubsection{Curve di genere $>1$.}
Il risultato $4n$ vale certamente, qualunque sia $p$, per le \emph{curve sghembe non speciali} (la $g^{3}_{n}$ segata su $\mathcal{C}$ dai piani dello spazio è non speciale) per $n-3\geq p$ (il risultato vale anche per $n\geq \frac{3}{4}(p+4)$, Brill e N\"{o}ether; Severi ne ha risolte le obbiezioni).

$\mathcal{C}^{p}_{n}$ si può rappresentare in una $F^{p}$ piana; ora \emph{perch\'{e} due curve con $p>1$ siano trasformazioni birazionali l'una dell'altra, occorrono $(3p-3)$ condizioni (Riemann). Una $F^{p>1}$ piana non ammette infinite trasformazioni birazionali in s\'{e} stessa (Schwarz)}, quindi non ci sono infinite rappresentazioni possibili della $\mathcal{C}$ sulla $f$.
Allora, per costruire una $\mathcal{C}^{p}$ si parte da una $F^{p}$ $$u=3p-3+\dots;$$
Poi fissare sulla curva una completa non speciale $g^{3}_{n}$
%dubbio: nel manoscritto presente un'annotazione a matita incerta, forse o n-p?% 
che importerà $\rho$ parametri: $u=3p-3+\rho+\dots$;\\
La $\mathcal{C}$ è determinata a meno di una collineazione, restano $15$ parametri disponibili $u=3p+12+\rho$. Resta a
 %-------------------------------pag. 113 -----------------------------------%
calcolare $\rho$. Dobbiamo fissare sulla curva modello una $g^{3}_{n}$ non speciale; essa è sempre contenuta in una $g^{n-p}_{n}$ completa. \emph{Quante sono le $g^{n-p}_{n}$ complete?} Si possono fissare sulla curva $n-p$ punti ad arbitrio; un gruppo della nostra serie è determinato ad essi; restano ancora $p$ punti; viceversa dati questi ad arbitrio, se li associo agli $n-p$ fissati una volta per tutte, ho $n$ punti che determinano la completa $g^{n-p}_{n}$ cui appartiene. Questa dunque dipende da $p$ parametri; allora $\rho=p+4(n-p-3)=4n-3p-12$
\begin{equation}
\label{eq:par3b}
u=4n
\end{equation}

\tableofcontents
\end{document}